\documentclass[final, a4paper, 12pt]{article}

\newcommand{\version}{4.1}

\usepackage{amsmath, amsthm, amssymb}
\usepackage{xspace}
\usepackage{comment}
\usepackage{paralist}

\author{Antongiulio Fornasiero%
  \thanks{Institut f\"ur Mathematische Logik,
    Einsteinstr.~62, 48149 M\"unster, Germany.\protect\\
\Email{antongiulio.
fornasiero@googlemail.com}
   } 
  \and Tamara Servi%
  \thanks{Centro de Matem\'atica e Aplica\c c\~oes Fundamentais,
    Av.\ Prof.\ Gama Pinto~2, 1649-003 Lisboa, Portugal.\protect\\
\Email{tamara.
servi@googlemail.com}
}}

\newcommand{\Email}[1]{\href{mailto:#1}{\tt<#1>}}

\def\shorttitle{Baire structures  v. \version}
\title{Definably complete Baire structures and Pfaffian closure\\
\medskip
\normalsize{Full article} \quad
\small{Version \version \quad file \tt\jobname.tex}}

\usepackage{newcent}                                                      

\usepackage{fancyhdr}
\makeatletter
\renewcommand\@makefnmark{{\normalfont({\scriptsize\@textsuperscript\@thefnmark})}}
\makeatother

\usepackage{tocloft}
\newlength{\tocwidth}
\cftsetpnumwidth{1.5em}
\cftsetrmarg{1.6em}

\usepackage[linktocpage=true, pdfborder={0 0 0} ]{hyperref}

\hypersetup{%
  pdfauthor={A. Fornasiero \& T. Servi},
  pdftitle={Definably complete Baire structures and Pfaffian closure},
  pdfsubject={Definably complete Baire structures and o-minimality of Pfaffian closure},
  pdfkeywords={Pfaffian functions; Pfaffian closure; definably complete
    structures; Baire spaces; o-minimality}
}
\newcommand{\ie}{i.e\mbox{.}\xspace}
\newcommand{\Ie}{I.e\mbox{.}\xspace}
\newcommand{\eg}{e.g\mbox{.}\xspace}

\newcommand{\cf}{cf\mbox{.}\xspace}

\DeclareMathOperator{\reg}{{reg}}
\newcommand{\Pa}[1]{\bigl( #1 \bigr)}
\newcommand{\set}[1]{\{#1\}}
\newcommand{\bigset}[1]{\bigl\{#1\bigr\}}
\newcommand{\abs}[1]{\lvert#1\rvert}
\newcommand{\norm}[1]{\lVert #1 \rVert}
\newcommand{\card}[1]{\lvert#1\rvert}
\newcommand{\Nat}{\mathbb{N}}
\newcommand{\Real}{\mathbb{R}}
\newcommand{\Rz}{\Real_0}
\newcommand{\K}{\mathbb{K}}
\newcommand{\pair}[1]{( #1 )}
\newcommand{\pairf}[1]{\langle #1 \rangle}
\newcommand{\R}{\mathbb R} \newcommand{\N}{\mathbb N} \newcommand{\Q}{\mathbb Q}
\newcommand{\F}{\mathbb F}

\newcommand{\imp}[1]{\textbf{#1}}
\newcommand{\dconnected}{definably connected\xspace}

\newcommand{\inter}[1]{\mathring{#1}}
\newcommand{\wideinter}[1]{\mathrm{int}(#1)}
\newcommand{\intert}[2]{\mathrm{int}_{#1}(#2)}
\newcommand{\cl}[1]{\overline{#1}}
\newcommand{\clt}[2]{\mathrm{cl}_{#1}(#2)}

\DeclareMathOperator{\vb}{vb}
\DeclareMathOperator{\cc}{cc}
\newcommand{\Bool}{\mathfrak{B}}

\newcommand{\ball}{B}

\newcommand{\Fs}{\mathcal{F}_{\sigma}}
\newcommand{\Gd}{\mathcal{G}_{\delta}}
\newcommand{\Afam}{\mathcal{A}}
\DeclareMathOperator{\bd}{bd}
\newcommand{\rest}{\mathbin\upharpoonright}
\newcommand{\Continuous}{\mathcal C}
\newcommand{\Czero}{\mathcal{C}^0}
\newcommand{\Cone}{\mathcal{C}^1}
\newcommand{\Ctwo}{\mathcal{C}^2}
\newcommand{\Cinf}{\mathcal{C}^\infty}
\newcommand{\CN}{\Continuous^N}
\newcommand{\CM}{\Continuous^M}
\newcommand{\Disc}{\mathcal D}
\newcommand{\Df}{\Disc_f}
\newcommand{\jhat}{\hat \jmath}

\newcommand{\DCN}{$\text{\rm DAC}^N$\xspace}
\newcommand{\Wlog}{W.l.o.g\mbox{.}\xspace}
\newcommand{\wloG}{w.l.o.g\mbox{.}\xspace}
\newcommand{\et}{\ \&\ }
\newcommand{\ao}{a.o\mbox{.}\xspace}

\newcommand{\dcompact}{d-compact\xspace}

\newcommand{\sdiff}{\mathop{\Delta}}
\newcommand{\meq}{\sim}
\DeclareMathOperator{\acc}{acc}

\newcommand{\Ffam}{\mathcal F}
\newcommand{\Gfam}{\mathcal G}

\newcommand{\PK}{\mathcal{P}(\Kt, \K)}
\newcommand{\RK}{\mathrm{Rolle}(\Kt,\K)}
\newcommand{\RKt}{\widetilde{\RK}}
\newcommand{\Kt}{\K_0}
\newcommand{\Kz}{\K_0}

\newtheorem{lem}{Lemma}[section]
\newtheorem{thm}[lem]{Theorem}
\newtheorem{cor}[lem]{Corollary}

\newtheorem{prop}[lem]{Proposition}
\newtheorem*{inductivehypothesis}{Inductive Hypothesis}

\theoremstyle{remark}

\newtheorem{claim}{Claim}
\newtheorem*{claim*}{Claim}
\newtheorem{oclaim}{Claim}

\theoremstyle{definition}
\newtheorem{dfn}[lem]{Definition}
\newtheorem{notation}[lem]{Notation}
\newtheorem{rem}[lem]{Remark}
\newtheorem{final remark}[lem]{Final remark}
\newtheorem{example}[lem]{Example}
\newtheorem{examples}[lem]{Examples}
\newtheorem{exa}[lem]{Example}
\newtheorem{open problem}[lem]{Open problem}
\newtheorem{open problems}[lem]{Open problems}

\newtheorem{proviso}[lem]{Proviso}
\newtheorem*{proviso*}{Proviso}
\swapnumbers

\newcommand{\case}[1]{\par\textsc{Case} #1.\xspace}

\makeatletter
\renewenvironment{proof}[1][\proofname]{\par
  \pushQED{\qed}%
  \normalfont \topsep6\p@\@plus6\p@\relax
  \setcounter{claim}{0}
  \trivlist
  \item[\hskip\labelsep
        \itshape
    #1\@addpunct{.}]\ignorespaces
}{%
  \popQED\endtrivlist\@endpefalse
}

\providecommand{\rom}{\textup}

\newcommand{\allsmall}{\forall^s}

\newcommand{\s}{{\mathcal S}}
\newcommand\tts{\widetilde{\mathcal{S}}}

\def\hyph{-\penalty0\hskip0pt\relax}

\newcommand{\intro}[1]{\textbf{#1}}

\newcommand{\pf}{correspondence\xspace}
\newcommand{\pfs}{correspondences\xspace}

\newcommand{\adm}{admissible\xspace}
\newcommand{\Adm}{Admissible\xspace}

\newcommand{\ov}[1]{\overline{#1}}
\newcommand{\ovs}[1]{\bar{#1}}
\newcommand{\x}{\ovs{x}}
\newcommand{\w}{\ovs{w}}
\newcommand{\z}{\ovs{z}}

\newcommand{\y}{\ovs{y}}
\newcommand{\av}{\ovs a}
\newcommand{\bv}{\ovs b}

\newcommand{\eps}{\varepsilon}
\newcommand{\veps}{{\ovs\varepsilon}}
\newcommand{\cll}[1]{\overline{#1}}
\DeclareMathOperator{\im}{Im}\DeclareMathOperator{\dom}{dom}
\newcommand{\de}{\mathrm{d}} 
\newcommand{\cof}{\textup{\textsf{Cofin}}}

\newcommand{\xseq}{(x_k)_{k < \cof(\K)}}\newcommand{\epseq}{(\eps_k)_{k < \cof(\K)}}

\DeclareMathOperator{\dist}{dist}

\newcommand{\cB}[1]{\overline{B(#1)}}

\newcommand{\Rt}{\mathfrak R}
\newcommand{\Lang}{\mathcal L}
\newcommand{\Th}{\mathrm{Th}}
\newcommand{\Tp}{T_0}
\newcommand{\Lp}{\Lang_0}
\newcommand{\RR}{\mathcal{P}(\Rt)}

\newcommand{\Cv}{\ovs C}
\newcommand{\psiv}{\ov \psi}

\newcommand{\tc}{\mathrm{t.c.}}
\newcommand{\re}{r.e\mbox{.}\xspace}
\newcommand{\RCF}{\ensuremath{\mathrm{RCF}}\xspace}
\newcommand{\Cfam}{\mathcal C}
\newcommand{\Pproof}{\mathbb P}
\newcommand{\Vreg}{V^{\reg}}

\begin{document}

\setlength{\emergencystretch}{1em}
\setlength{\hfuzz}{7pt}
\hbadness=10000

\maketitle

\begin{abstract} 
We consider definably complete Baire expansions of ordered fields: 
every definable subset of the domain of the structure has a supremum
and the domain can not be written as the union of a definable increasing
family of nowhere dense sets. Every expansion of the real field is definably
complete and Baire, and so is every o-minimal expansion of a field.
Moreover, unlike the o-minimal case, the structures considered form an
axiomatizable class. 
In this context we prove the following version of Wilkie's Theorem of the 
Complement: given a definably complete Baire expansion $K$ of an ordered field
with a family of smooth functions, if there are uniform bounds on the
number of definably connected components of quantifier free definable sets,
then $K$ is o-minimal. 
We further generalize the above result, along the line of Speissegger's
theorem, and prove the o-minimality of the relative Pfaffian closure of an
o-minimal structure inside a definably complete Baire structure. 
\end{abstract}

\textit{Key words:}
Pfaffian functions; Pfaffian closure; definably complete structures;
Baire spaces; o-minimality\\
\indent\textsl{MSC2000:} Primary
58A17; 
Secondary
03C64, 
32C05, 
54E52. 

{
\tableofcontents
}

\section{Introduction}\label{section:intro}

We recall that a subset $A$ of a topological space $X$ is said to be meager if
there exists a collection $\set{Y_i:\ i\in\Nat}$ of nowhere dense subsets of
$X$ such that $A\subseteq\bigcup_{i\in\Nat}Y_i$.
The Baire Category Theorem implies that every open subset of~$\Real$ (with the
usual topology) is not meager, \ie $\Real$ is a Baire space.

The notion of Baire space is clearly not first order. 
Here we consider a similar ({\em definable}) notion, which instead is preserved under elementary equivalence, and which coincides with the classical notion over the real numbers (this is made precise in Section~\ref{section: meager}).

The (first order) structures we consider are definably complete expansions of ordered fields.
Definable completeness (see Definition \ref{def:def-complete}) is a weak version of Dedekind completeness, which is preserved under elementary equivalence. It is shown in \cite{miller}, \cite{servi-articolo}, \cite{frat} that, as in the o-minimal case, (a definable version of) most results of elementary real analysis can be proved in every definably complete expansion of an ordered field. However, to obtain less elementary results one would need some more sophisticated machinery, in the direction of Sard's Lemma and Fubini's Theorem. Both of the quoted classical results refer to a notion of smallness (having measure zero), which has no natural translation in our context. We consider instead a topological notion of smallness (being meager), propose a definable version of this notion and carry out a theory of {\em definably complete Baire structures}, i.e. expansions of ordered fields such that every definable subset of the domain has a supremum and the domain can not be written as the union of a definable increasing family of nowhere dense sets. In this context we prove an analogue to Fubini's Theorem (the Kuratowski-Ulam's Theorem \ref{thm:Fubini}) and a very restricted form of an analogue to Sard's Lemma (Theorem~\ref{sard for noeth rings}). Notice that it is not known whether every definably complete structure is definably Baire.

\bigskip

Once we have developed the basic tools for definably complete Baire structures
(Sections \ref{section: meager} to \ref{section:open problems}), our next task
is to give necessary and sufficient conditions, for a definably complete
expansion with $\Cinf$ functions of an ordered field, to be o-minimal.

In \cite{wilkie99}, the author proves his Theorem of the Complement:
given an expansion $\mathcal R$ of the real field with a family of $\Cinf$
functions, if there are bounds (uniform in the parameters) on the number of
connected components of quantifier free definable sets, then $\mathcal R$ is
o-minimal.
In particular, thanks to a well known finiteness result in~\cite{khovanskii},
the structure generated by all real Pfaffian functions is o-minimal
(see \cite{khovanskii} or \cite{wilkie99} for the definition of Pfaffian
functions and examples).
In \cite{KM99}, the authors generalize Wilkie's Theorem of the Complement (by
weakening the smoothness assumption) in a way which allows them to derive the
following result (originally due to Speissegger, see~\cite{speissegger}): the
Pfaffian closure of an o-minimal expansion of the real field is o-minimal.

In Section~\ref{section:complement} we proceed to generalize the o-minimality
results present in~\cite{wilkie99} and~\cite{KM99}, to a situation where the
base field is not necessarily~$\Real$; moreover, we further weaken the
assumption of Wilkie's Theorem of the Complement, by allowing not only
functions, but also \imp{\adm \pfs{}} (roughly, partial multi-valued functions
with finitely many values at each point).
We deduce that, given a definably complete Baire expansion $\K$ of an
ordered field with a family of $\Cinf$ functions, if there are bounds (uniform
in the parameters) on the number of definably connected components of
quantifier free definable sets, then $\K$ is o-minimal (Theorem~\ref{o-minimality theorem}).
In Section~\ref{section:pfaffian}, by using our restricted version~\ref{sard for noeth rings} of Sard's Lemma,
we proceed to prove the analogue to Khovanskii's finiteness result in
the context of definably complete Baire structures (Theorem~\ref{khov1e2}).
We derive the o-minimality of every definably complete Baire expansion of an
ordered field with any family of definable Pfaffian functions
(Theorem~\ref{pfaff is omin}).

Finally, in Section~\ref{section:speissegger}, we prove that the relative
Pfaffian closure of an o-minimal structure $\K_0$ inside a definably complete
Baire expansion $\K$ of $\K_0$ is o-minimal.
This latter result, whose proof is shaped on the one present in \cite{KM99},
can be compared with the main result in \cite{frat} (which can be derived from
ours), where instead Speissegger's method was followed; it is here where our
generalization of Wilkie's Theorem to \adm \pfs is necessary.

In Section~\ref{sec:effective} we use the above results to find effective
bounds for various topological invariants of sets definable in the Pfaffian
closure of the fields of reals, and more generally of recursively axiomatized
o-minimal expansions of~$\Real$.

\bigskip
The results in this article have been submitted for publication.
Since we do not have constraints of space, we opted to give more detailed
proofs, explanations and examples that would be suitable for a published
version.

\subsection{Notation}
Throughout this paper, $\K$~is a (first-order) structure expanding an ordered field.
We use the word ``definable'' as a shorthand for ``definable in $\K$ with parameters from $\K$''.

We denote by $x,y,z,\ldots$ the points in $\K^n$. 
When we want to stress the fact that they are tuples, we write
$\x,\y,\z,\ldots$, where $\x=(x_1,\ldots,x_n)$, etc.

For convenience, on $\K^m$ instead of the usual Euclidean distance we will use
the equivalent distance 
\[
d : (x,y) \mapsto \max_{i = 1, \dotsc, m} \abs{x_i - y_i}.
\]
For every $\delta > 0$ and $x \in \K^m$, we define
\[\begin{aligned}
B^m(x;\delta) &:= \set{y \in \K^m: d(x,y) < \delta},\\
\cl B^m(x;\delta) &:= \set{y \in \K^m: d(x,y) \leq \delta},
\end{aligned}\]
the open and closed ``balls'' of center $x$ and ``radius''~$\delta$; 
we will drop the superscript $m$ if it is clear from the context.

\begin{notation}
Let $X \subseteq Y \subseteq \K^n$, with $Y$ definable.
We write $\clt Y X$
\rom(or simply $\cl X$ if $Y$ is clear from the context\rom)
for the topological closure of $X$ in~$Y$, $\intert Y X$
\rom(or simply $\inter X$\rom)  for the interior part of $X$ in~$Y$,
$\bd_Y(X) := \cl{X} \setminus \inter{X}$ for the boundary of~$X$
\rom(in~$Y$\rom) ,
and $\partial_Y X := \cl X \setminus X$ for the frontier of~$X$
\rom(in~$Y$\rom).
\end{notation}

\begin{notation} We define $\Pi^{m+n}_n:\K^{m+n}\to\K^m$ as the projection onto the first $m$ coordinates. If $A\subset\K^{m+n}$ and $x\in\K^m$, we denote by $A_x$ the fibre of $A$ over $x$, i.e. the set $\set{y\in\K^n:\ (x,y)\in A}$.
\end{notation}

\begin{notation}
Let $\K_+:=\set{x\in\K:\ x>0}$.
\end{notation}

\begin{dfn}\label{def:Rtilde}
Let $\tilde \Real$ be the structure on the reals numbers, with a predicate for
every subset of~$\Real^n$ \rom(it will be used for examples\rom).
\end{dfn}

\subsection{Definably complete structures}
\begin{dfn}\label{def:def-complete} An expansion $\K$ of an ordered field is called {\em definably complete} if every definable subset of $\K$ has a supremum in $\K\cup\set{\pm\infty}$.
\end{dfn}

Generalities on definably complete structures can be found in \cite{servi-articolo}, \cite[\S 2]{DMS08} and \cite{miller}.

\begin{proviso*}
For the remainder of the article,
$\K$~will always be a definably complete structure.
\end{proviso*}

\begin{dfn}
$X \subseteq \K^m$ is {\em definably compact} \rom(\dcompact for short\rom)
if it is definable, closed in $\K^m$ and bounded.
\end{dfn}

\begin{proviso}\label{cofinal set}
We order $\K^m$ lexicographically.
In this subsection we will denote by $N$ a definable subset of $\K^m$ which is
cofinal in the lexicographic ordering.
\end{proviso}

\begin{lem}[Miller]\label{lem:1}
$X$ is definably compact iff for every $\Pa{Y(y)}_{y \in N}$ definable
decreasing family of closed non empty subsets of~$X$,
we have $\bigcap_y Y(y) \neq \emptyset$.
\end{lem}

\begin{dfn}Let $f: N \to \K^n$ be definable.
Define $\acc_{y\to \infty} f(y)$ (and write for simplicity~$\acc f$)
to be the set of accumulation points of~$f$;
that is, $x \in \acc f$ iff
\[
(\forall r \in \K^m) (\forall \varepsilon \in \K_+) (\exists y > r)\
y \in N \ \& \ d(f(y), x) < \varepsilon.
\]
\end{dfn}

\begin{lem}\label{lem:2}
If $X$ is definably compact, then for all definable~$N$
\rom(satisfying \ref{cofinal set}\rom)
and for all $f: N \to X$ definable we have
$\acc f \neq \emptyset$.
\end{lem}
It is not clear if the converse of the above lemma is true.

\begin{dfn}
Let $(A(y))_{y \in N}$ be a definable family of non-empty
subsets of~$\K^n$.
Define $\acc_{y\to \infty} A(y)$ (and write for simplicity $\acc A$)
to be the set of accumulation points of~$A$;
that is, $x \in \acc A$ iff
$(\forall r \in \K^m)$ $(\forall \varepsilon \in \K_+)$ $(\exists y > r)$
$y \in N$ and $d(A(y), x) < \varepsilon$.

Note that $\acc A = \bigcap_y \cl{\Pa{ \bigcup_{z \geq y} A(z)}}$.
\end{dfn}

\begin{rem}
Let $(A(t))_{0 < t \in \K}$ be a definable family of subsets of~$\K^m$,
and $G := \bigcup_{t > 0} \set t \times A(t)$.
Then, $\acc_{t \to 0} A(t) = \Pa{\clt{\K^{m + 1}} G}_0 :=
\set{x \in \K^m: (0,x) \in \cl G}$.
\end{rem}

\begin{lem}\label{accumulation point for compact sets}
$X$ is definably compact iff for all $A$ definable family of non-empty subsets of $X$ we have $X\cap\acc A \neq \emptyset$.
\end{lem}
\begin{proof}
First assume that $X$ is \dcompact. Let $Y(y) := \cl{\Pa{\bigcup_{z \geq y} A(y)}}$. Then $(X\cap Y(y))$ is a definable decreasing family of closed subsets of $X$.
By Lemma~\ref{lem:1}, $\bigcap_y Y(y) \neq \emptyset$, and we are done.

Conversely, assume that $X$ is not \dcompact.
By Lemma~\ref{lem:1},
there exists a definable decreasing family $Y:=(Y(y))_{y \in \K^m}$ of closed
subsets of~$X$ such that  $\bigcap_y Y(y) = \emptyset$.
However, since $Y$ is decreasing, $X \cap \acc Y = \bigcap_y Y(y)$, and we are done.
\end{proof}

\begin{proof}[Proof of Lemma~\ref{lem:2}]
Define $A(y) := \set{f(y)}$.
By Lemma~\ref{accumulation point for compact sets}, $\acc A$ is non-empty.
Note that $\acc A = \acc f$.
\end{proof}

\begin{lem}\label{lem:exchange}
Let $X \subseteq \K^n$ be definably compact.
Fix $\varepsilon \in \K_+$.
Let $(A(t))_{t \in N}$ be a definable family of subsets of~$\K^n$.
The following are equivalent
\begin{enumerate}
\item $\forall x \in X$ $\forall t \in N$ large enough
$X \cap B(x;\varepsilon) \subseteq A(t)$;
\item $\forall t \in N$ large enough $X \subseteq A(t)$.
\end{enumerate}
\end{lem}
\begin{proof}
That (2) implies (1) is clear.

Conversely, assume that (1) is true.
Suppose, for contradiction, that (2) is false.
Let $D(t) := X \setminus A(t)$.
Let $N' := \set{t \in N: D(t) \neq \emptyset}$.
Since (2) is false, $N'$ is cofinal in~$N$.
Let $C := \acc_{t \in N', t \to \infty} D(t)$.
By Lemma~\ref{accumulation point for compact sets}, $C \neq \emptyset$; let $x \in C$.
By (1), if $t$ is large enough,
then $X \cap B(x;\varepsilon) \subseteq A(t)$.
Choose $t \in N'$ such that $X \cap B(x;\varepsilon) \subseteq A(t)$
and $d(x, D(t)) < \varepsilon$.
Let $y \in D(t)$ such that $d(x,y) < \varepsilon$.
Since $y \in D(t)$, we have $y \notin A(t)$.
Since $y \in X \cap  B(x;\varepsilon)$, we have $y \in A(t)$, a contradiction.
\end{proof}

\begin{lem}\label{Lebesgue number}
Let $C\subset\K^n$ be a nonempty \dcompact set, and let $V:=\set{V(t):\ t\in I}$ be a definable open cover of~$C$.
Then, there exists $\delta_0\in\K_+$
\rom(a \emph{Lebesgue number} for $V$ and~$C$\rom)
such that, for every subset $X\subseteq C$ of diameter smaller than $\delta_0$, there exists $t\in I$ such that $X\subseteq V(t)$.
\end{lem}

\begin{proof} Suppose for a contradiction that
\[
(\forall\delta > 0) (\exists y\in C) (\forall t\in I)\
 B(y;\delta)\nsubseteq V(t).
\]
For every $\delta > 0$, define
\[
Y(\delta) := \set{y \in C:  (\forall t\in I)\  B(y; \delta)\nsubseteq V(t)}.
\]
Note that  $\Pa{Y(\delta)}_{\delta>0}$ is a definable family of subsets of~$C$,
increasing as $\delta$ decreases.
Let $y_0$ be an accumulation point for the family $\Pa{Y(\delta)}_{\delta>0}$, as $\delta \to 0$ (which exists by Lemma \ref{accumulation point for compact sets}).

Let $t_0 \in I$ and $\delta_0 > 0$ such that $B(y_0; 2\delta_0) \subseteq V(t_0)$.
Let $\delta_1 \leq \delta_0$ and $y \in Y(\delta_1)$ such that $\abs{y - y_0} < \delta_0$.
Therefore, $B(y; \delta_0) \subseteq B(y_0; 2\delta_0) \subseteq V(t_0)$,
contradicting the fact that $y \in Y(\delta_1)$.
\end{proof}

We will often use without further comment the following result:

\begin{lem}[Miller] Let $f:\K^n\to\K^m$ be a definable continuous function and let $C\subset\K^n$ be \dcompact. Then $f(C)$ is \dcompact.
\end{lem}

\begin{dfn}\label{dfn:manifold}
A $n$-dimensional definable embedded $\CN$ $\K$-manifold $V \subseteq
\K^d$ \rom(which we will simply call $n$-dimensional $\K$-manifold\rom)
is a definable subset $V$ of~$\K^d$, such that for every $x \in V$ there
exists  a definable neighbourhood $U(x)$ of~$x$ \rom(in~$\K^d$\rom),
and a definable $\CN$ diffeomorphism $f_x: U(x) \simeq \K^d$, such that
$U(x) \cap V = f_x^{-1}\Pa{\K^n \times \set 0}$.
\end{dfn}

\begin{rem}\label{manifold}Note that a $\K$-manifold $V$ can always be written
as the intersection of a definable closed set and a definable open set.
In fact, let $\delta : V \to \K_+ \cup \set{+ \infty}$ be
the definable map
\[
\delta(x) := \sup\set{r \in \K_+:
\forall s \in \K_+\ \Pa{ s < r \Rightarrow B(x; s) \cap V \text{ is closed in }
B(x; s)}}.
\]
Let $U := \bigcup_{x \in V} B\Pa{x; \delta(x)/2}$; then, $V = \cl V \cap U$.

Note moreover that the dimension $n$ of a $\K$-manifold $V$ is uniquely
determined by~$V$, because $\K^n$ and $\K^{n'}$ are locally diffeomorphic iff
$n = n'$.
If we consider only $\Czero$ manifolds, it is not clear anymore if the dimension is well defined.
\end{rem}

Finally, recall the following definition.

\begin{dfn}\label{def:connected}
A definable set $X\subset\K^n$ is {\em definably connected} if it can not be expressed as a union of two {\em definable} nonempty disjoint open sets. A subset $C\subseteq X$ is a {\em definably connected component} of $X$ if it is a maximal definably connected subset of~$X$.
\end{dfn}

Note that if $X$ has finitely many definably connected components, then
each component of $X$ is definable. Moreover, if $\K$ expands the real field, every definable and (topologically) connected set is also definably connected. The converse could in general not be true. However it is true if $\K$ is o-minimal. For example, it is true for any expansion of the real field by a Pfaffian chain (see Theorem~\ref{pfaff is omin}).

\section{Meager sets}\label{section: meager}
Let $X \subseteq Y \subseteq \K^n$, with $Y$ definable.
\begin{dfn}
$X$~is {\em nowhere dense} \rom(in~$Y$\rom)
if $\intert Y {\clt Y X} = \emptyset$.
$X$~is {\em definably meager} \rom(in~$Y$\rom) if there exists a
{\em definable increasing} family $(A(t))_{t \in \K}$
of nowhere dense subsets of~$Y$, such that  $X \subseteq \bigcup_t A(t)$.
We will call the family $\Pa{\clt Y {A(t)}}_{t \in \K}$ a {\em witness} of
the fact that $X$ is definably meager.
$X$~is {\em definably residual} \rom(in~$Y$\rom) if $Y \setminus X$ is
definably meager.
\end{dfn}

Notice that, if $(A(t))_{t \in \K}$ is a witness of the fact that $X$ is
meager in $\K^n$, then also the family
\[
\Pa{\cl B^n(0; \abs t) \cap A(t)}_{t \in \K}
\]
is a witness, hence we may always assume that each $A(t)$ is \dcompact.

Notice also that we do {\em not} require that a definably meager set is
definable.

The subsets of $Y$, with the operations $\Delta$ (symmetric
difference) and~$\cap$, form a commutative ring;
the definably meager subsets of $Y$ form an ideal of this ring.

\begin{dfn}\label{def:Baire}
$Y$~is {\em definably Baire} if every non-empty open definable subset of $Y$
is not definably meager \rom(in~$Y$\rom).
\end{dfn}

Note that if $\K$ has countable cofinality, then $X$ is definably meager
(Baire, respectively) in $\K^n$ if $X$ is meager (Baire, respectively) in the
usual topological sense.
In general, the converse is not true: for instance, if $\K$ is a countable o-minimal structure, then it is definably Baire, but not Baire in the topological sense.
However, the two notions coincide for $\tilde \Real$.
In fact, assume that $X$ is a meager subset of~$\Real^n$;
therefore, $X = \bigcup_{i \in \Nat} Y(i)$,
where each $Y(i)$ is a nowhere dense subset of~$\Real^n$.
For each $t \in \Real$, define $Z(t) := \bigcup_{i \in \Nat, i < t} Y(i)$.
Then, each $Z(t)$ is nowhere dense, the family $\Pa{Z(t)}_{t \in \Real}$ is
increasing and definable in~$\tilde \Real$, and $X = \bigcup_{t \in \K} Z(t)$.

From now on, we will write ``meager'' for ``definably meager'', and ``topologically meager'' for the usual topological notion, and
similarly for ``residual'' and ``Baire''.
Moreover, if $Y$ is clear from the context, we will simply say that $X$ is
nowhere dense (resp.\ ``meager'', ``residual'') instead of
``nowhere dense'' (resp.\ ``definably meager'', ``definably residual'') in~$Y$.

\begin{prop}\label{prop:open}Let $Y$ be definable, and $\emptyset \neq U \subseteq Y$ be definable and open. Then,
 $U$ is meager in~$Y$ iff it is meager in itself.
\end{prop}

\begin{proof}
Suppose $U$ is meager in $Y$ and let $\Pa{Y(t)}_{t \in \K}$ be a witness of this fact.
For every $t \in \K$, define $X(t) := Y(t) \cap U$.
Since $U$ is open, $\intert U{X(t)} = \intert Y {Y(t)} \cap U = \emptyset$.
Hence, $\Pa{X(t)}_{t \in \K}$ is a witness of the fact that $U$
is meager in itself.

Vice versa,
let $\Pa{X(t)}_{t \in \K}$ be a witness of the fact that $U$ is meager in
itself, and $Y(t) := \clt Y{X(t)}$.
We claim that $\intert Y{Y(t)} = \emptyset$.
In fact, $\intert Y{Y(t)}=\intert Y{(\clt Y{X(t)}\cap U)}=\intert U{Y(t)}=\intert U{X(t)}=\emptyset$.
Hence, $\Pa{Y(t)}_{t \in \K}$ is a witness of the fact that $U$
is meager in~$Y$.
\end{proof}

\begin{cor}\label{lem:open}
Let $Y$ be definable, and $\emptyset \neq U \subseteq Y$ be definable and open. Then,
\begin{enumerate}
\item If $U$ is of not meager in itself, then $Y$ is also not meager in itself.
\item If $Y$ is Baire, then $U$ is also Baire.
\end{enumerate}
\end{cor}
\begin{proof}

For (1), if $Y$ is meager in itself, then any subset of $Y$, in particular  $U$, is meager in $Y$. Since it is open, $U$ is also meager in itself.

Regarding (2),
if $Y$ is Baire, let $V \subseteq U$ be non-empty, definable and open in~$U$.
Since $U$ is open in~$Y$, $V$~is also open in~$V$.
Hence, by \ref{prop:open}, $V$ is not meager in
itself, and, again by \ref{prop:open}, $V$~is not meager
in~$U$. Therefore, $U$~is Baire.
\end{proof}

\begin{lem}\label{lem:Baire}
Let $Y \subseteq \K^m$ be definable.
The following are equivalent:
\begin{enumerate}
\item $Y$ is Baire;
\item for all $X \subseteq Y$, if $X$ is meager, then $\inter X = \emptyset$;
\item every $x \in Y$ has a definable neighbourhood which is Baire;
\item every residual subset of $Y$ is dense;
\item every open definable non-empty subset of $Y$ is not meager in
itself;
\item every meager closed definable subset of $Y$ has empty interior.
\end{enumerate}
\end{lem}
\begin{proof}\
\begin{enumerate}
\item[($2 \Rightarrow 1$)] is obvious.
\item[($1 \Rightarrow 3$)] is obvious, because $Y$ itself is a Baire
neighbourhood of each point.
\item[($3 \Rightarrow 4$)] Let $X \subseteq Y$ be meager.
Suppose, for a contradiction, that $U$ is a non-empty definable subset of $X$
open in~$Y$, and let $x \in U$.
Let $V$ be a definable Baire neighbourhood of~$x$, and $W := V \cap U$.
By Proposition~\ref{prop:open}, $W$ is Baire, and therefore it is not meager in~$Y$
(by the same proposition), which is not possible.
\item[($4 \Rightarrow 2$)] Let $X \subseteq Y$ be meager. Hence, $Y \setminus
X$ is dense, and therefore $\inter X = \emptyset$.
\item[($1 \Leftrightarrow 5$)] Use Proposition~\ref{prop:open}.
\item[($1 \Rightarrow 6$)]
Let $C \subseteq Y$ be definable, closed and meager.
If $\inter C \neq \emptyset$, then $\inter C$ is not meager, and thus $C$ is
not meager.
\item[($6 \Rightarrow 1$)]
Let $U \subseteq Y$ be open, definable and meager in~$Y$.
Then, $\cl U$ is also meager, because $\cl U = U \sqcup \bd U$, and $\bd U$ is
nowhere dense.
Therefore, $\cl U$ has empty interior, and therefore $U$ is empty.
\qedhere
\end{enumerate}
\end{proof}

\begin{rem}\label{K is baire iff nonmeager}
$\K^n$~is Baire iff it is not meager in itself.
\end{rem}

\begin{proof}
One implication is obvious.

For the other implication, assume that $\K^n$ is not meager in itself, and
let $U \subseteq \K^n$ be an open definable subset.
If, for a contradiction, $U$ were meager in itself, then we could find an open non empty box $B \subseteq U$.
By Proposition~\ref{prop:open}, $B$~is also meager in itself.
However, $B$ is definably homeomorphic to $\K^n$, because $\K$ expands a field, contradicting the hypothesis.
\end{proof}

The following result is not trivial and will be proved in Section \ref{section: kuratowski}.
\begin{prop}\label{cor:Baire}
If $\K$ is Baire, then for every $m \geq 1$, $\K^m$ is Baire.
\end{prop}

The converse, however, \emph{is} trivial
\begin{rem}
If $\K^m$ is Baire for some $m \geq 1$, then $\K$ is Baire.
\end{rem}

\subsection{Baire structures}

\begin{dfn}
A definably complete structure~$\K$ is a \emph{Baire structure}
if $\K$ is definably Baire as a definable subset of~$\K$ itself,
in the sense of Def.~\ref{def:Baire}.
A~theory $T$ is definably complete and Baire if every model of $T$ is a
definably complete Baire structure.
\end{dfn}

\begin{rem}
The fact that $\K$ is Baire can be expressed by a set of first-order
sentences: therefore, every $\K'$ elementary equivalent to $\K$ also satisfies
the hypothesis.
If moreover the language is recursive, this set of sentences is also recursive.

Notice that an ultra-product of definably complete \rom(resp.\ Baire\rom)
structures is also definably complete \rom(resp.\ Baire\rom);
the same cannot be said for ``o-minimal'' instead of ``definably complete''.
\end{rem}

\begin{examples}
The following are examples of definably complete Baire structures.
\begin{itemize}
\item
Every expansion of $\Real$ 
(because $\Real$ is Dedekind complete and topologically Baire).
\item
Every o-minimal expansion of a field.
In fact, a nowhere dense definable subset of $\K$ is finite,
and definable families of finite sets are uniformly finite;
hence, the union of a definable increasing family of nowhere dense sets is finite, and can not coincide with the whole structure.
\item
Let $\mathcal{B}$ be an o-minimal expansion of a field, let
$A\preccurlyeq\mathcal{B}$ be a dense substructure. Then the structure
$\mathcal{B}_A$, generated by adding a unary predicate symbol for $A$, is
definably Baire. This follows from the fact that if $X\subseteq \mathcal{B}$
is $\mathcal{B}_A$-definable, then its topological closure $\cl{X}$ is
$\mathcal{B}$-definable (see \cite[Theorem 4]{vddries}).
Hence, a closed nowhere dense set is finite, and, since $\mathcal{B}_A$
satisfies the Uniform Finiteness property (see \cite[Corollary 4.5]{vddries}),
the union of a definable increasing family of nowhere dense sets is finite.
More generally, as shown in \cite[\S 3.5]{DMS08}, any definably complete
structure satisfying the Uniform Finiteness condition is definably
Baire (one can even show that if $\K$ is definably complete, and  every
definable closed discrete subset of $\K$ is bounded, then $\K$ is Baire).
\end{itemize}
\end{examples}


\section{\texorpdfstring{$\Fs$-sets}{F-sigma sets}}\label{section:Fs}
We now consider a class of sets for which it is easy to determine whether they are meager or not. This sets have also been studied in \cite{DMS08}, where they are called $D_{\Sigma}$-sets.

\begin{dfn}
Let $X \subseteq Y \subseteq \K^n$, with $Y$ definable.
$X$~is in $\Fs$ in~$Y$ (we will also say ``$X$ is an $\Fs$ subset of $Y$'', or
``$X$ is $\Fs$'',
and drop the reference to $Y$ if it is clear from the context)
if $X$ is the union of a definable increasing family of closed subsets of~$Y$, indexed by~$\K$.
$X$~is in $\Gd$ if its complement is an~$\Fs$.
\end{dfn}

\begin{lem}\label{lem:Fs}
Let $\Afam$  be either the family of $\Fs$ or the family of $\Gd$ subsets of some~$\K^n$, for $n \in \Nat$.
Then, each $A \in \Afam$ is definable.
Moreover, $\Afam$~is closed under finite unions, finite intersections,
Cartesian products, and preimages under definable continuous functions.
Besides, the following are in~$\Afam$
\begin{enumerate}
\item definable closed subsets of $\K^n$;
\item definable open subsets of $\K^n$;
\item finite boolean combinations of definable open subsets of~$\K^n$.
\end{enumerate}
The family of $\Fs$ subsets is also closed under images under definable continuous functions.
\end{lem}

\begin{proof}
See \cite{DMS08}.
Let $A$ and $B$ be in $\Fs$.
The fact that $A \cup B$ and $A \times B$ are also in $\Fs$ is obvious.

Let $A = \bigcup_{t} A(t)$ and $B = \bigcup_{t} B(t)$,
where $(A(t))_{t \in \K}$ and $(B(t))_{t \in \K}$ are two definable increasing families of closed (and, we may assume, \dcompact) sets.

Then, $A \cap B = \bigcup_{t} (A(t) \cap B(t))$, because
$(A(t))_{t \in \K}$ and $(B(t))_{t \in \K}$ are increasing families.
Hence, $A \cap B$ is also in $\Fs$.

If $f: \K^n \to \K^{n'}$ is continuous, then $f(A) = \bigcup_t f(A(t))$.
For every $t \in \K$, $f(A(t))$ is \dcompact, because $A(t)$ is \dcompact, and therefore $f(A)$ is in~$\Fs$.
A similar proof works for preimages.

Let $U \subseteq \K^n$ be open and definable, and $C := \K^n \setminus U$.
For every $r \in \K_{+}$, define $U(r) := \set{ x \in \K^n: d(x, C) \geq r}$.
Note that each $U(r)$ is closed.
Since $U$ is open, $U = \bigcup_{r > 0} U(r)$, and therefore~$U$ is in $\Fs$.

If $D$ is a finite boolean combination of open definable subsets of~$\K^n$, then it is a finite union of sets of the form $C_i \cap U_i$, for some definable sets $C_i$ and~$U_i$, such that each $C_i$ is closed and each $U_i$ is open.
Hence, $D$~is in~$\Fs$.

The corresponding results for $\Gd$ follow immediately by considering the complements.
\end{proof}
It is not true in general that, if $Y \subseteq \K^n$ is definable, $X
\subseteq Y$ is an $\Fs$-subset of~$Y$, and $f: Y \to Y$ is definable and
continuous,  then $f(X)$ is an~$\Fs$.
The point where the above proof breaks down for $Y \neq \K^n$ is the fact that
it is not necessarily true that every $\Fs$ subset of $Y$ is an increasing definable union of \dcompact sets.

Notice that, by Remark \ref{manifold}, every $\K$-manifold is an $\Fs$-set.

\begin{rem}
Let $X \subseteq \K^n$.
$X$~is an $\Fs$ iff $X$ is of the form $\Pi^{n+m}_n(Z)$
for some $Z \subseteq \K^{n + m}$ closed and definable.
\end{rem}
\begin{proof}
The ``if'' direction follows from Lemma~\ref{lem:Fs}.
For the other direction, let $\Pa{X(t)}_{t \in \K}$ be a definable increasing family of closed subsets of~$\K^n$, such that $X = \bigcup_{t \in \K} X(t)$.
Define $Z := \ov{\bigsqcup_{t \in \K} \Pa{X(t) \times \set t}}$.
\end{proof}

Notice that, if $\K$ is o-minimal, then every $X$ definable subset of $\K$ is a finite Boolean combination of definable closed sets (because $X$ is a finite union of cells), and therefore $X$ is an~$\Fs$.

\begin{rem}
If $X \subseteq \K^n$ is meager, then there exists a meager $\Fs$-set containing~$X$.
\end{rem}

\begin{lem}\label{lem:meager-F}
Let $Y$ be definable and Baire, and $D \subseteq Y$.
Assume that $D$ is in~$\Fs$.
Then, $D$ is meager iff $\inter D = \emptyset$.%
\footnote{This is not true for $\Gd$ sets: for instance, the set of irrational numbers in $\Real$ is a $\Gd$ which is not meager (it is even residual), but has empty interior.}
\end{lem}
\begin{proof}
If $\inter D \neq \emptyset$, then, since $Y$ is Baire, $D$~cannot be meager.
Conversely, assume that $D$ is not meager.
If $D$ is in~$\Fs$, then $D = \bigcup_{t} D(t)$, for some definable increasing family of closed subsets.
Since $D$ is not meager, at least one of the~$D(t)$, say~$D({t_0})$, is not meager.
Hence, $\wideinter{D({t_0})} \neq \emptyset$ (otherwise, $D({t_0})$~would be
nowhere dense),  and therefore $\inter{D} \neq \emptyset$.
\end{proof}

Note that if $X \subseteq \Real^n$ is in $\Fs$ for the structure
$\tilde \Real$, and of Lebesgue measure zero,
then $X$ is meager, but the converse is not true.

We now give a local condition which is sufficient to prove that the image of an $\Fs$-set under a continuous definable function is meager.

\begin{prop}\label{prop:local-meager}
Let $C \subseteq \K^m \times \K^n$ be in $\Fs$, $f: C \to \K^d$ be definable
and continuous.
Assume that for every $y \in \Pi^{m+n}_m(C)$ there exists a
neighbourhood $V_y \subseteq \K^m$ of~$y$, such that $f
\Pa{(V_y\times\K^n)\cap C }$ is meager.
Then, $f(C)$~is meager.
\end{prop}
\begin{proof}
If $\K$ is meager in itself, then by Proposition \ref{cor:Baire} there is
nothing to prove.
Thus, we may assume that $\K$ (and hence $\K^d$) is Baire.

We proceed by induction on~$m$.
The case $m = 0$ is clear, because if $m = 0$, then $V_0 = \K^0$.

Assume that we have already proved the conclusion for $m - 1$ (and every~$n$).
We want to prove it for~$m$.
First, we consider the case when $C$ is \dcompact.
\Wlog, $0 \in C$.
Remember that, for every $r > 0$ and $y \in \K^m$,
$\cl B^m(y; r)\subset\K^m$ is the
closed hypercube of side $2r$ and center~$y$;
let $S^m(y; r)$ be its boundary.
Moreover, define
$D(r) := f \Pa{ C \cap (\cl B^m(0; r)\times  \K^n ) }$.

Note that $f(C) = \bigcup_r D(r)$, and that each $D(r)$ is \dcompact.
Therefore, to prove that $f(C)$ is meager,
it suffices to prove that each $D(r)$ has empty interior.
Suppose, for a contradiction, that $f(C)$ is not meager, and let
\[
r_0 := 
\inf \set{r > 0 : \wideinter{D(r)} \neq \emptyset}.
\]
Since the $D(r)$ are closed,
$r_0 = \inf \set{r > 0 : D(r) \text{ is not meager}}$.
We have that $0 < r_0$ by hypothesis, and $r_0 < + \infty$
because $f(C)$ is not meager.

Let $P := \Pi^{n+m}_m(C)$.
Since $P$ is \dcompact, if $\K = \Real$, we could find
$y_{1}, \dotsc, y_k \in P$ such that
$P \subseteq V_{y_1} \cup \dots \cup V_{y_k}$.
In the general situation, we need another argument.
Let $5 \delta_0$ be a Lebesgue number for the open cover $\set{V_y:\ y\in P}$
of $P$ (we may also assume that $\delta_0$ is small in comparison with~$r_0$);
$\delta_0 > 0$ exists by Lemma~\ref{Lebesgue number}.


Note that
\[
\cl B^m(0; r_0 + \delta_0/2)\subseteq \cl B^m(0; r_0 - \delta_0/2) \cup
\bigcup_{y \in   S^m(0; r_0)}\cl B^m(y; \delta_0),
\]
hence
\[
D(r_0 + \delta_0/2) \subseteq D(r_0 - \delta_0/2) \cup \bigcup_{y \in S^m(0;r_0)}
f \Pa{ C \cap( \cl B^m(y; \delta_0)\times\K^n)}.
\]
By definition of~$r_0$, we know that $D(r_0 + \delta_0/2)$ is not meager, while
$D(r_0 - \delta_0/2)$ is meager.
Hence, to obtain a contradiction, it suffices to show that
$\bigcup_{y \in S^m(0;r_0)} f \Pa{ C \cap( \cl B^m(y;\delta_0)\cap\K^n)}$
is meager.

Note that $S^m(0;r_0)$ is the finite union of the faces of the closed
hypercube $\cl B^m(0; r_0)$:
hence, we only need to show that for each face~$S$ of $S^m(0;r_0)$ the set
$D := \bigcup_{y \in S} f \Pa{ C \cap( \cl B^m(y;\delta_0)\times\K^n)}$
is meager.
\Wlog, we can assume that $S$ is the ``top'' face
$\set{y \in \cl B^m(0;r_0): y_m =  r_0}$ and we may identify
$S$ with $\cl B^{m-1}(0;r_0)\times\set{r_0}$.

Define
\[\begin{aligned}
\tilde C  & := C \cap \bigcup_{y \in S} \Pa{ \cl B^m(y;\delta_0)\times\K^n},\\
\tilde f  &:= f \rest {\tilde C}.
\end{aligned}\]

\begin{claim*}
$\tilde C$ and $\tilde f$
satisfy the hypothesis of the proposition, with $n' = n + 1$,
$m' = m - 1$, and $V'_z = B(z; \delta_0)$.
\end{claim*}

$\tilde C$~is \dcompact,
and therefore it is in~$\Fs$.
Let $\tilde P \subseteq \K^{m - 1}$ be the projection of
$tilde C$ onto $\K^{m - 1}$; note that $\tilde P$ is \dcompact.
Fix $z \in \tilde P$; by definition, there exists $t \in [r_0-\delta_0, r_0+\delta_0]$ such that $y := (z, t) \in P$.
Notice that
\[
\tilde C \cap(  V'_z \times\K\times\K^{n}) \subseteq
C \cap ( V'_z \times[r_0-\delta_0,r_0+\delta_0]\times\K^{n} ) \subseteq
C \cap (\cl B^m(y; 2\delta_0)\times\K^n) .
\]

Since $5 \delta_0$ is a Lebesgue number for the cover $\set{V_y:\ y\in P}$
of~$P$, it follows 
that there exists $y'\in P$ such that $\cl B^m(y;2\delta_0)\subset V_{y'}$.
Putting everything together,
we have that $\tilde C \cap(  V'_z \times\K^{n+1})\subset C \cap ( V_{y'}\times\K^n)$
and thus $\tilde f \Pa{ \tilde C \cap ( V'_z \times\K^{n + 1} ) }$ is meager, which proves the claim.

Therefore, by inductive hypothesis, $\tilde f(\tilde C)$ is meager.
However, $D \subseteq \tilde f(\tilde C)$, and we reached a contradiction.

We now treat the general case when $C$ is in~$\Fs$.
Note that $C$ is an increasing union of \dcompact sets~$C(t)$.
For each $t \in \K$, define $D(t) := f(C(t))$: note that each $D(t)$ is \dcompact.
By the \dcompact case, we can conclude that each $D(t)$ is meager, and therefore
nowhere dense.
Thus, $D = \bigcup_t D(t)$ is meager.
\end{proof}

\begin{cor}\label{cor:local-meager}
Let $C \subseteq \K^m$ be in~$\Fs$,
and $f : C \to \K^d$ be definable and continuous.
Assume that for every $x \in C$ there exists $V_x \subseteq C$
neighbourhood of~$x$, such that $f \Pa{ C \cap V_x }$ is meager.
Then, $f(C)$ is meager.
\end{cor}
\begin{proof}
Apply the proposition to the case $n = 0$.
\end{proof}

With a similar method, one can prove the following.
\begin{lem}
Let $C \subseteq \K^m$ be \dcompact and $f: C \to \K^d$ be definable
\rom(but not necessarily continuous\rom).
Assume that for every $x \in C$ there exists $V_x \subseteq C$
neighbourhood of~$x$, such that $f \Pa{ C \cap V_x }$ is nowhere dense.
Then, $f(C)$ is nowhere dense.
\end{lem}

\begin{cor}
Let $W \subseteq \K^m$ be a definable $\K$-manifold,
$C \subseteq W$ be an~$\Fs$ subset of~$W$,
and $f : C \to \K^d$ be definable and continuous.
Assume that for every $x \in C$ there exists $V_x$
neighbourhood of~$x$, such that $f \Pa{C \cap V_x }$ is meager.
Then, $f(C)$ is meager.
\end{cor}
\begin{proof}
Since $W$ is a $\K$-manifold, it is in~$\Fs$.
Since $C$ is $\Fs$ in~$W$, it is also $\Fs$ in~$\K^m$.
Apply the previous corollary.
\end{proof}

\begin{cor}
Let $C \subseteq \K^m$ be an $\Fs$.
If every $x \in C$ has a neighbourhood $V_x$ such that $C \cap V_x$ is
meager, then $C$ is meager.
\end{cor}

Proposition~\ref{prop:local-meager} and the following results are trivial
if $\K$ is o-minimal, since in this case $C \subseteq \K^n$ is meager iff
$\dim(C) < n$.

For the topological notions, we know the following facts to be true:
\begin{enumerate}
\item
Let $C \subseteq \Real^n$ and $f: C \to Y$ (not necessarily continuous).
Assume that, for every $x \in C$, there exists $V_x \subseteq C$ neighbourhood
of~$x$, such that $f(C \cap V_x)$ is topologically meager.
Then, $f(C)$ is topologically meager.
\item
Let $C \subseteq Y$. If every $x \in C$ has a neighbourhood~$V_x$,
such that $V_x \cap C$ is topologically meager, then $C$ is topologically
meager.
\end{enumerate}
The first fact follows from the fact that $\Real^n$ is second countable;
the second from~\cite[Theorem~6.35]{Kelley55}.
We were able to prove the definable versions only under additional
hypothesis (\eg, $C$ in~$\Fs$); however, these results are strong enough
for our applications.


\begin{lem}\label{lem:discontinuity}
Let $f: \K^n \to \K^m$ be definable, and
\[
\Df := \set{\x\in A:\ f\mathrm{\ is\ discontinuous\ at\ }\x}.
\]
If the graph of $f$ is an $\Fs$ set, then $\Df$ is meager.
\end{lem}
\begin{proof}
If, for contradiction, $\Df$~is not meager, then
since it is an $\Fs$, it contains a non-empty open box~$B$.
Therefore, \wloG we can assume that $\Df = \K^n$, and that $\K^n$ is Baire.
Let $\Gamma(f) = \bigcup_t X(t)$, where $\Pa{X(t)}_{t \in \K}$ is a definable
increasing family of \dcompact sets.
Let $Y(t) := \Pi^{n+m}_n \Pa{X(t)}$.
Note that each $Y(t)$ is \dcompact, and $\K^n = \bigcup_t Y(t)$.
Since $\K^n$ is Baire, there exists $t_0$ such that $Y(t)$ contains a
non-empty open box~$B'$.
Let $B''  \subseteq B'$ be a closed box with non-empty interior, and $g := f
\rest B''$.
Note that $\Gamma(g) = X(t_0) \cap (B'' \times \K^m)$; therefore, $\Gamma(g)$
is \dcompact, and so, as in the classical case, $g$ is continuous,
contradicting the fact that $B'' \subseteq \Df$.
\end{proof}

\section{The Kuratowski-Ulam's Theorem}\label{section: kuratowski}

The main result of this section is the following theorem.

\begin{thm}\label{thm:Fubini}
Let $D \subseteq \K^{m+n}$.
For every $x \in \K^m$, let $D_x := \set{y \in \K^n: (x,y) \in D}$ be the corresponding section of~$D$.
Let $T := T^m(D) := \set{x \in \K^m: D_x \text{ is meager in } \K^n}$.

If $D$ is meager \rom(in~$\K^{m+n}$\rom), then $T$ is residual.
\end{thm}

This is a definable version of Kuratowski-Ulam's Theorem
\cite[Theorem~15.1]{Oxtoby80}, which in turn is an analogue of Fubini's Theorem:
they both imply that if $D$ is negligible, then $D_y$ is negligible for almost
every~$y$;  in Kuratowski-Ulam's Theorem negligible means ``meager'', while in Fubini's Theorem negligible means ``of measure zero''.

It is not clear whether in the above theorem $D$ definable implies that $T$ is
definable.
Note that if $\K$ is o-minimal and $D$ is definable, then $T$ is also definable.

As a corollary, we obtain Proposition \ref{cor:Baire}.

\begin{proof}[Proof of Proposition \ref{cor:Baire}]
By induction on~$m$.
The case $m = 1$ is our assumption on~$\K$.
Assume that we already proved that $\K^m$ is Baire: we want to prove that
$\K^{m + 1}$ is Baire.
Suppose not; then $\K^{m + 1}$ is meager in itself.
If we apply Theorem~\ref{thm:Fubini} with $n = 1$, we obtain that either $\K^m$ or $\K$ is meager in itself, a contradiction.
\end{proof}


\begin{dfn}
A definable function $f: Y \to \K$ is {\em lower semi-continuous} if, for every $x \in Y$,
either $x$ is an isolated point of~$Y$, or
\[
\liminf_{ \stackrel{x' \to x}{ x' \in Y}} f(x') \geq f(x).
\]
\end{dfn}

\begin{rem}
Let $C  \subseteq \K^{n+1}$ be \dcompact.
For every $x \in D := \Pi^{n+1}_n(C)$, let $f(x) := \min C_x$.
Then, $f: D \to \K$ is lower semi-continuous.
\end{rem}

\begin{lem}\label{lem:semicontinuous}
Let $Y \subseteq \K^n$ be definable,
$f: Y \to \K$ be lower semi-continuous and definable,
and $\Df \subseteq Y$ be the set of points of discontinuity of~$f$.
Then, $\Df$~is meager \rom(in~$Y$\rom).
\end{lem}
\begin{proof}
See \cite[Lemma 2.8(1)]{DMS08}.
\end{proof}

In the above lemma, if $Y = \K = \Real$, we can not conclude that $\Df$
has Lebesgue measure zero.
In fact, let $C \subseteq \Real$ be closed, with empty interior, and of
positive measure, and  $f$ be the characteristic function of
$\Real \setminus C$.
Then, $\Df = C$, and therefore it is of positive measure.

On the other hand, it is always true that if $f : \K^m \to \K$ is definable,
then $\Df$ is in $\Fs$ (see \cite[Theorem~7.1]{Oxtoby80}).

\bigskip

\begin{proof}[Proof of Theorem~\ref{thm:Fubini}]
If $\K^m$ is meager in itself, then the conclusion is trivially true, because
 then every subset of $\K^m$ is meager.
Hence, we can assume that $\K^m$ is Baire.

\case{1} $n = 1$ and $D$ is \dcompact.\\
Hence, $D$ has empty interior, and
each $D_x$ is also \dcompact.
Therefore, by Lemma~\ref{lem:meager-F},
$T = \set{x \in \K^m: \inter{D_x} = \emptyset}$.
Let $E := \K^m \setminus T$.
We have to prove that $E$ is meager.

For every $\varepsilon > 0$ let
\[
X(\varepsilon) := \set{(x,y) \in \K^m \times \K: \ball^1(y;\varepsilon) \subseteq D_x}.
\]
Let
\[
E(\varepsilon) := \pi(X(\varepsilon)) = \set{x \in \K^m: D_x \text{ contains a ball of radius } \varepsilon}.
\]
Note that $X(\varepsilon)$ is \dcompact, since its complement is the projection of an open set, therefore so is $E(\varepsilon)$.
Note that $E = \bigcup_{\varepsilon > 0} E(\varepsilon)$; hence, to prove that $E$ is meager, it suffices to prove that each $E(\varepsilon)$ is nowhere dense.
Since each $E(\varepsilon)$ is \dcompact, it suffices to prove the following claim.

\begin{claim}
For every $\varepsilon > 0$, $\wideinter{E(\varepsilon)} = \emptyset$ (see
also \cite[Lemma 2.8(2)]{DMS08}).
\end{claim}

Assume, for a contradiction, that there exists a nonempty open box $U \subseteq E(\varepsilon)$.
Define
\[\begin{aligned}
f:\ & U \to \K\\
  & x \mapsto \min\set{y \in \K: (x,y) \in X(\varepsilon)}.
\end{aligned}\]
Note that $f$ is lower semi-continuous and definable.
By Lemma~\ref{lem:semicontinuous}, $f$~is continuous outside a meager
set~$\Df \subseteq U$.
Since $\K^m$ is Baire, $\Df \neq U$, and therefore
there exists $x_0 \in U$ such that $f$ is continuous at~$x_0$.
It is now easy to show that a neighbourhood of $(x_0,f(x_0))$ is contained
in~$D$, contradicting the fact that $\inter{D} = \emptyset$.

\case 2
$n = 1$ and $D$ arbitrary meager subset of $\K^m$.\\
Let $\Pa{D{(p)}}_{p \in \K}$ be an increasing definable family of \dcompact
subsets of $\K^{m + 1}$ with empty interior, such that $D \subseteq
\bigcup_{p} D{(p)}$.
For each $p \in \K$, let $E{(p)} := \set{x \in \K^m: D{(p)}_x \text{ is
    not meager in } \K}$.
By what we have seen above, $E{(p)} = \bigcup_{\varepsilon > 0}
E(p,\varepsilon)$, where $(E(p,\varepsilon))_{\stackrel{ \varepsilon
    \in \K_+}{ p \in \K}}$ is a definable family of subsets of~$\K$,
increasing in $p$ and decreasing in $\varepsilon$, such that each
$E(p,\varepsilon)$ is closed and nowhere dense.
Let
\[
E' := \bigcup_{\varepsilon, p} E(p,\varepsilon)
= \bigcup_p E{(p)}.
\]

\begin{claim}
$\K^m \setminus T \subseteq E'$.
\end{claim}

In fact, let $x \notin T$.
Thus, $D_x$ is not meager.
However, $D_x \subseteq \bigcup_p D{(p)}_x$.
Since $\Pa{D{(p)}_x}_{p \in \K}$ is an increasing definable family of closed
subsets of~$\K$, we obtain that there exists $p_0$ such that $D{(p_0)}_x$ has
non-empty interior.
Thus, $x \in E{(p_0)} \subseteq E'$.

Therefore, it suffices to prove that $E'$ is meager to obtain that $T$ is residual.
However, $E' = \bigcup_{p > 0} E(p,1/p)$, and we are done.

\case 3
$n > 1$ and $D$ arbitrary meager subset of $\K^m$. We argue by
induction on~$n$.\\
Suppose that we have already proved the conclusion for $n$ (and for every~$m$).
We want to prove the conclusion for $n + 1$.
First, we will assume that $D$ is in~$\Fs$.
We want to prove that the set $T := T^m(D) := \set{x \in \K^m: D_x \text{ is meager}}$ is residual.
Define
\[\begin{aligned}
S &:= \K^{m + 1} \setminus T^{m + 1}(D) := \set{(x,y_{n + 1}) \in \K^m \times \K : D_{(x,y_{n + 1})}
\text{ is not meager}},\\
R &:= T^m(S) = \set{x \in \K^m: S_x \text{ is meager}}.
\end{aligned}\]


Notice that (for the moment) we do not know whether $S$  and $R$ are
definable, even assuming that $D$ is in~$\Fs$.

\begin{claim}
$S$ is meager.
\end{claim}

By inductive hypothesis.

\begin{claim}
$R$ is residual.
\end{claim}

By the case $n = 1$ and the previous claim.

\begin{claim}
$R \subseteq T$.
\end{claim}

Fix $x \in \K^m$.
Assume that $x \notin T$.
We have to prove that $x \notin R$.
Define $F := D_x \subseteq \K^{n + 1}$.
Note that $F$ is in~$\Fs$: therefore,
since $x \notin T$, $\inter{F} \neq \emptyset$.
Let $U := U_1 \times U_2$ be a non-empty open box contained in~$F$, $U_1
\subseteq \K$, $U_2 \subseteq \K^m$.
For every $y_{n + 1} \in U_1$, $D_{(x,y_{n+1})} = F_{y_{n+1}} \supseteq U_2$,
and therefore $(x,y_{n+1}) \in S$.
Thus, $U_1 \subseteq S_x$, and $x \notin R$.

Hence, $T$ contains a residual set, and therefore it is residual.

For $D$ arbitrary, let $D' \subseteq \K^{m + n}$ be a meager $\Fs$
containing~$D$.
By the previous case, the corresponding set $T' := T^m(D')$ is residual.
Since $T' \subseteq T$, we are done.
\end{proof}

\section{Almost open sets}\label{section:a.o.}
\begin{proviso*}
In this section we will assume that $\K$ is {definably complete and Baire}.
\end{proviso*}

Let $Y \subseteq \K^m$ be definable.
We have seen that the family of meager subsets of $Y$ is an ideal, hence it defines an equivalence relation on the family of subsets of~$Y$, given by
$X \meq X'$ iff $X \sdiff X'$ is meager.

\begin{rem}
$X \meq X'$ iff there exists $Z$ meager such that $X \sdiff Z = X'$
\end{rem}
\begin{proof}
Set $Z := X \sdiff X'$.
\end{proof}

\begin{dfn}
$X \subseteq Y$ is {\em almost open} \rom(in~$Y$\rom),
or \ao for short, if $X$ is equivalent to a definable open set.%
\footnote{Almost open sets are called ``sets with the property of Baire''
in~\cite{Oxtoby80}.}
\end{dfn}

\begin{lem}\label{lem:almost open}
Let $Y \subseteq \K^m$ be definable, and $A$ and $B$ be \ao subsets
of~$Y$.
Then, $A \cap B$, $A \cup B$ and $Y \setminus A$ are also \ao.
Moreover, $\Fs$ and $\Gd$ subsets of $Y$ are \ao.

Finally, if $Y_1$ and $Y_2$ are definable, and $A_i \subseteq Y_i$ are \ao for $i = 1, 2$, then $A_1 \times A_2$ is \ao in $Y_1 \times Y_2$.
\end{lem}
\begin{proof}
It is trivial to see that $A \cap B$, $A \cup B$ and $A_1 \times A_2$ are \ao.

Let $A = U \sdiff E$, where $U$ is open and definable, and $E$ is meager.
Then, $Y \setminus A = (Y \setminus U) \sdiff E$.
Hence, to prove that $Y \setminus A$ is \ao it suffices to prove that
$C := Y \setminus U$ is \ao.
However, $C = \inter{C} \cup \bd(C)$.
Since $C$ is closed, $\bd(C)$ is nowhere dense, and \textit{a fortiori} meager, and we are done.

Let $(D(t))_{t \in \K}$ be a definable increasing sequence of closed subsets
of~$Y$. We have to prove that $D := \bigcup_t D(t)$ is \ao.
Let $U := \inter{D}$ and $E := D \setminus U$.
It is enough to prove that $E$ is meager.
For every $t$, let $E(t) := E \cap D(t)$.
Note that $\inter{D(t)} \subseteq U$; therefore, $E(t) \subseteq \bd(D(t))$ is
nowhere dense, and we are done.
\end{proof}

Consequently, $X \subseteq Y$ is \ao iff it is equivalent to a definable closed subset of~$Y$.

\begin{rem} Every meager set is \ao, being equivalent to the empty set.
Every residual set is also \ao, being equivalent to the ambient space.
\end{rem}

\begin{cor}
Let $A \subseteq Y$.
The following are equivalent:
\begin{enumerate}
\item $A$ is \ao;
\item $A$ is of the form $E \sdiff F$, for some meager set $E$ and some set
$F$ in $\Fs$;
\item $A$ is of the form $G \sqcup E$, for some $G$ in $\Gd$ and $E$ meager.
\end{enumerate}
\end{cor}
\begin{proof}Cf.\cite[Theorem~4.4]{Oxtoby80}.
$(1 \Leftrightarrow 2)$ and  $(3 \Rightarrow 1)$ are obvious.
For $(1 \Rightarrow 3)$, let $A = U \sdiff E$ for some $U$ open and $E$ meager.
Let $Q$ be a meager set in $\Fs$ containing~$E$, and $G := U \setminus Q$.
Note that $G$ is in $\Gd$, and
\[
U \sdiff E = [(U \setminus Q) \sdiff (U \cap Q)] \sdiff (E \cap Q) =
G \sdiff [(U \sdiff E) \cap Q] = G \sqcup E',
\]
where $E' := (U \sdiff E) \cap Q$ is meager.
\end{proof}



The following is a partial converse of Theorem~\ref{thm:Fubini}.

\begin{lem}\label{lem:Fubini-converse}
Let $D$ be an \ao subset of $\K^{ m+n}$, and $T(D) := \set{x \in \K^m: D_x
  \text{ is meager}}$.
Then, $D$ is meager iff $T(D)$ is residual.
\end{lem}
\begin{proof}
The \lq\lq only if'' direction is Theorem~\ref{thm:Fubini}.
For the other direction, let $U$ be an open set such that $E := D
\sdiff U$ is meager.
By Theorem~\ref{thm:Fubini}, $T(E)$~is residual.
Moreover, since $U_x = D_x \sdiff E_x$, we have $T(U) \supseteq T(D) \cap
T(E)$, and therefore $T(U)$ is also residual.
However, $U$~is open and $\K^n$ is Baire: therefore, $T(U)$ is the complement of the projection of $U$ on~$\K^m$.
Since $U$ is open, $T(U)$ is closed.
Therefore, $T(U)$ is closed and residual; since $\K^m$ is Baire, $T(U) =
\K^m$. Thus, $U$ is empty, and we are done.
\end{proof}

The hypothesis that $D$ is \ao in the above lemma is necessary:
\cite[Theorem~15.5]{Oxtoby80} gives an example of a set $E \subseteq \Real^2$
that is not topologically meager,
such that no three points of $E$ are collinear.


\section{Further results and open problems}\label{section:open problems}

\begin{open problem} It is not known to the authors if there exists a definably complete structure which is not Baire.
\end{open problem}

\begin{proviso*}
For the remainder of this section,
$\K$~is a definably complete Baire structure.
\end{proviso*}

\begin{open problem}
Let $\Pa{Y(t)}_{t \in \K}$ be a definable increasing family of {\em meager} subsets of~$\K^m$,
and let $Y := \bigcup_t Y(t)$. Is $Y$ necessarily meager? In particular, is it
necessarily $Y\not=\K^m$?
\end{open problem}

Notice that, if in addition the $Y(t)$ are closed,
then $Y$ is meager, whereas the same conclusion does not necessarily hold if
the $Y(t)$ are in $\Fs$ (actually, since every meager set is contained in a
meager $\Fs$-set, it is enough to reduce to this situation).
Moreover, the above question has positive answer if $\K$ is o-minimal, because
then each $Y_t$ has (o-minimal) dimension less than $m$, and therefore $Y$ has
dimension less than~$m$.
In fact, if $\K$ is o-minimal, and $Y \subseteq \K^m$ is definable, then $Y$
is meager iff $\dim Y < m$; moreover,
if $\Pa{Y(t)}_{ t \in \K_+}$ is a definable family, decreasing
in~$t$, then $\bigcup_t Y(t) \subseteq \acc_{t \to 0} Y(t)$.
Thus, the following lemma proves what we want.

\begin{lem}
Let $\K$ be an o-minimal structure,
$n \leq m \in \Nat$, and $\Pa{Y(t)}_{t > 0}$ be a definable family of
subsets of~$\K^m$, and $Z := \acc_{t \to 0} Y(t)$.
If, for every $t > 0$, $\dim \Pa{Y(t)} \leq n$, then $\dim(Z) \leq n$.
\end{lem}
\begin{proof}
Define $W := \bigcup_{t>0} Y(t) \times \set t \subseteq \K^{m+1}$.
Note that $Z = (\cl W)_0 := \set{x \in \K^m: (x,0) \in \cl W}$.
Moreover, since $\dim Y(t) \leq n$, we have $\dim W \leq n + 1$.
Since $Z \times \set 0 \subseteq \partial W $, we have
$\dim Z < \dim W \leq n + 1$.
\end{proof}


The following is a partial result for the case of \ao sets.
\begin{lem}\label{lem:meager-union}
Let $Y \subseteq \K^n$ be definable and Baire, $D \subseteq Y$ be \ao \rom(in
$Y$\rom), and $(Y(t))_{t \in \K}$ be a definable increasing family of closed subsets of~$Y$, such that $Y = \bigcup_t Y(t)$.
Then, $D$~is meager in $Y$ iff each $D \cap Y(t)$ is meager \rom(in~$Y$\rom).
\end{lem}
\begin{proof}
The ``only if'' direction is clear.

For the other direction, let $C \subseteq Y$ be closed, such that $E := C
\sdiff D$ is meager.
It suffices to prove that $C$ is meager.
For every $t \in \K$, define
\[\begin{aligned}
C(t) & := C \cap Y(t),\\
D(t) & := D \cap Y(t).
\end{aligned}\]
Then, $D(t) \sdiff C(t) \subseteq E$.
Therefore, since $D(t)$ and $E$ are meager, $C(t)$ is meager and closed.
Since $Y$ is Baire, $C(t)$ is nowhere dense, and thus $C$ is meager.
\end{proof}




\bigskip


\subsection{The Sard property}

Let $C \subseteq \K^n$ be meager, and $f : \K^n \to \K^n$ be definable
and~$\Cone$.
We want to investigate in which circumstances $f(C)$ is meager.
When $\K=\Real$, Sard's Lemma implies that $f(C)$ is meager.
This suggests the following definition.

\begin{dfn}
Fix $d,r,m$ positive natural numbers.
Let $V\subseteq\K^d$ be a $\K$-manifold of dimension~$n$.
Let $f:V\to\K^m$ be a definable $\Continuous^r$ function and $\Delta_f$ be the set
of singular points of $f$. If $\Sigma_f := f(\Delta_f)$ is meager in $\K^m$, then we say
that $f$ has the {\em Sard property}.%
\end{dfn}

\begin{lem}
If $\K = \tilde \Real$ and $f: V \to \K^m$
is as in the above definition, with $r>\max\set{0,n-m}$,
then $f$ has the Sard property.
\end{lem}
\begin{proof}
By Sard's Lemma, $\Sigma_f$ has Lebesgue measure zero, and therefore it has empty interior.
Since $\Sigma_f$ is in~$\Fs$, it is also meager.
\end{proof}

\begin{open problem}\label{open prob sard}
Does every $\Continuous^r$ definable function $f: \K^n\to\K^m$ 
\rom(with $r>\max\set{0,n-m}$\rom) have the Sard property?
\end{open problem}

\begin{rem}
If $\K$ is o-minimal, then every $\Cone$ definable function $f: V \to \K^m$ has
the Sard property~\cite[Theorem~3.5]{BO01}.
\end{rem}



\begin{prop}
Suppose $f:\K^n\to\K^n$ has the Sard property, and let $C\subset\K^n$ be meager. Then $f(C)$ is meager.
\end{prop}
\begin{proof}
We may assume that $C\in\Fs$, since $C$ is contained in a meager $\Fs$-set. Let $\Lambda:=\K^n\setminus\Delta_f$ be the set of regular points of~$f$.
Note that $\Lambda$ is open.

By the Sard property, $f(C\cap\Delta_f)$ is meager.
Hence, it suffices to show that $f(C \cap \Lambda)$ is
meager.
Let $x \in C \cap \Lambda$.
Since $x$ is a regular point for~$f$, by the Implicit Function Theorem there exists a neighbourhood $V$ of $x$
such that $f$ is a diffeomorphism on~$V$; therefore, $f(C \cap V)$~is meager,
and, by Corollary~\ref{cor:local-meager}, $f(C \cap \Lambda)$ is meager.
\end{proof}


The following lemma is a generalization of Lemma~\ref{lem:Fubini-converse}.
\begin{lem}
Let $f : \K^n \to \K^m$ be a $\Cone$ definable function with the Sard property.
Let $\Lambda$~be the set of the regular points of~$f$, and $C \subseteq \K^n$
be almost open. 
For every $t \in \K^m$, let $F_t := f^{-1}(t)$, $C_t := F_t \cap C$, and
$T := \set{t \in \K^m: C_t \text{ is meager in } F_t \text{ or } 
F_t = \emptyset}$. 
Then, $T$~is residual iff $C \cap \Lambda$ is meager.
\end{lem}
\begin{proof}
If $n < m$, $\Lambda = \emptyset$, and we have a tautology. 
Assume that $n \geq m$.

Let $x \in C \cap \Lambda$.
Since $x$ is a regular point for~$f$, there exists $V$ open neighbourhood
of $x$ such that, up to a change of coordinates, $f\rest V$~is the projection
on the first $m$ coordinates $y := (x_1, \dotsc, x_m)$.
For every $y \in T$, the set $C_{y}\cap V$ is meager
in~$\set{y} \times \K^{n - m}$.

Hence, if $T$ is residual, then, by Lemma \ref{lem:Fubini-converse}, $C \cap V$
is meager; therefore, by Corollary \ref{cor:local-meager}, $C \cap \Lambda$ is meager.

Conversely, assume that $C \cap \Lambda$ is meager; we must prove that $T$ is
residual.
Since $f(\Delta_f)$ is meager, it suffices to prove that $T(C \cap \Lambda)$ is
residual.
Therefore, \wloG we can assume that $C \subseteq \Lambda$.
By Kuratowski-Ulam's Theorem \ref{thm:Fubini}, the set
$T(V \cap C) := \set{y \in \K^m:   C_{y} \cap V \text{ is meager in } 
\K^{n - m}}$ is residual.
Therefore, $T$~is residual.
\end{proof}

In Subsections \ref{preliminary} and \ref{subs:sard} we will produce examples of classes of functions in definably complete Baire structures, which have the Sard property.

\section[A theorem of the complement]{A theorem of the complement for a class of definably complete Baire structures}\label{section:complement}

In this section we prove a version of Wilkie's Theorem of the Complement \cite[Theorem 1.9]{wilkie99} which holds not only, as the original theorem, for expansions of the real field, but also for definably complete Baire structures. This result will give a sufficient criterium to establish if a given definably complete Baire structure is in fact o-minimal.

We will assume the reader to have familiarity with \cite{wilkie99} (which, in turn, uses results from \cite{max}) and we will adapt the proofs contained therein to our situation. We will occasionally refer to the treatment of Wilkie's Theorem of the Complement given in \cite{BS04}, when more suitable to our purposes.

We recall a few definitions (corresponding to \cite[Definitions 1.1, 1.2, 1.3, 1.6]{wilkie99}), adapted to our present situation.

\begin{proviso*}
We fix for the rest of the article a {definably complete Baire}
structure~$\K$.
\end{proviso*}

\begin{dfn}\label{gamma}
For $X\subseteq \K^n$ definable, let $cc(X)$ be the number of \dconnected
components of $X$ (Def.~\ref{def:connected}),
and let $\gamma(X)$ be the least $m \in \N$ such that,
for every affine set $L\subseteq \K^n$,
we have $cc(X\cap L)\leq m$,
with the convention that $\gamma(X) = \infty$ if $m$ does not exist.
\end{dfn}

\begin{dfn} \label{weak structure}  Let $\s = \langle \s_n \colon n \in \N^+
\rangle$, where $\s_n$ is
a collection of definable subsets of $\K^n$. We say that $\s$ is a \emph{weak structure} \rom(over~$\K$\rom)
if $\s$ contains all zero-sets
of polynomials with coefficients in $\K$ and is closed under finite
intersection, cartesian product and permutation of the variables.%

$\s$ is \textit{closed} if for every $n$ and $A\in
\s_n$, $A$ is a closed subset of $\K^n$; $\s$ is \textit{semi-closed} if for
every $n$ and $A\in \s_n$, $A$ can be obtained as the projection onto the first
$n$ coordinates of some closed set $B\in \s_{n+k}$, for some suitable~$k$.
$\s$~is \textit{o-minimal} if for every $n$ and $A \in \s_n$
we have $\gamma(A) < \infty$.
\end{dfn}

\begin{dfn} \label{Ch}
Let $\s$ be an o-minimal weak structure \rom(over~$\K$\rom).
The {\em Charbonnel closure} $\tts = \langle \tts_n \colon n
\in \N^+\rangle$ is obtained from $\s$ by closing under the following {\em Charbonnel operations}: finite union, intersection with affine sets, projection and topological closure\footnote{The set of operations defined here gives rise to the same closure as the one originally defined by Charbonnel, see \cite{BS04}.}.
\end{dfn}

We immediately obtain an analogous result to \cite[Lemmas 1.4, 1.5]{wilkie99}:

\begin{thm}\label{charbonnel}
If $\mathcal S$ is a semi-closed o-minimal weak structure, then its Charbonnel closure
$\tts$ is a semi-closed o-minimal weak structure.
\end{thm}

The reader can easily check that the proof of the quoted lemmas contained in \cite[\S 1]{max} does not use specific properties of $\R$, and can be reformulated in any definably complete structure (the Baire property is not needed here). Definable completeness is necessary because the fact that a continuous definable function on a closed bounded definable set assumes maximum is used to bound the $\gamma $ of the topological closure of a set.

\begin{dfn}\label{cinfinito}Let $\s$ be an o-minimal weak structure.
We say that $\s$ is \textit{determined by its smooth functions}
\rom(DSF\rom) if, given a set $A\in\s_n$, there exist $k\in\N$ and a
$\Cinf$-function $f_A:\K^{n+k}\rightarrow\K$ whose graph lies in
$\tts$, such that $A$ is the projection onto the first $n$ coordinates
of the zero-set of $f_A$ (compare with \cite[Definition 1.7]{wilkie99}).
\end{dfn}

The aim is to prove the following version of \cite[Theorem 1.8]{wilkie99}.

\begin{thm}\label{wilkie1.8}
Let $\s$ be a semi-closed o-minimal weak structure \rom(over $\K$\rom),
which is DSF.
Then $\tts$ is closed under complementation and o-minimal.
\end{thm}

The following version of \cite[Theorem 1.9]{wilkie99} will then automatically follow.

\begin{thm}\label{o-minimality theorem} Let $\K$ be a definably complete Baire structure and $\cal F$ be a family of $\K$-definable $\Cinf$ functions. Let $\K_{\cal F}$ be the reduct of $\K$ generated by the field structure and the functions in $\cal F$.  Then $\K_{\cal F}$ is o-minimal if and only if $\gamma(A)<\infty$, for every $A\subset\bigcup_{n\geq 1}\K^n$ {\em quantifier free} $\K_{\cal F}$-definable set.
\end{thm}

In view of Theorem~\ref{charbonnel}, to prove Theorem~\ref{wilkie1.8} it is sufficient to show that, under the hypothesis of the statement, $\tts$ is closed under complementation.

In \cite{KM99}, the authors generalized Wilkie's Theorem of the Complement \cite[Theorem 1.9]{wilkie99} (by weakening the DSF assumption) in a way which allowed them to derive the o-minimality of the Pfaffian closure of an o-minimal expansion of the real field. Inspired by \cite{KM99}, we will also weaken our DSF assumption and prove a more general statement (Theorem~\ref{thm:KM}), from which Theorem \ref{wilkie1.8} will follow as a corollary. The motivation for giving such a general statement will be clear in Section~\ref{section:speissegger}, where we will show an application.

\bigskip

In Subsection~\ref{subsec:adm} we give some results on {\adm \pfs}, which will play a role in the statement of \ref{thm:KM}.
In Subsection \ref{preliminary} we develop some preliminary results
(corresponding to the results in \cite[\S 2]{wilkie99}) about o-minimal weak structures. In Subsection \ref{subs:boundary} we can finally state our result precisely, and we proceed as in \cite[\S 3]{wilkie99} and give the key ingredient of the proof (the Theorem of the Boundary~\ref{boundary-DCN}). Finally, in Subsection \ref{subsec:cell dec} we conclude the proof by adapting Wilkie's Cell Decomposition Theorem (which can be found in \cite[\S 4]{wilkie99}) to our situation.

\subsection{\Adm \pfs}\label{subsec:adm}

To be able to state exactly the result we want to prove, we need to give some
definitions.
All the results in this subsection do not need that $\K$ is Baire, but only
that it is definably complete.

\begin{dfn}
A \pf $f: \K^n \leadsto \K^m$ is a definable partial multi-valued
function from $\K^n$ to $\K^m$.
\end{dfn}

\begin{dfn}
Given $1 \leq N \in \Nat$, a $\CN$ \intro{\adm \pf} is
a \pf $f : \K^n \leadsto \K^m$,
satisfying the following conditions.
Let $F \subset \K^n \times \K^m$ be the graph of~$f$.
\begin{enumerate}
\item $F$ is definable and
has a finite number of \dconnected components;
\item $F$ is a $\CN$ \emph{closed} embedded submanifold of $\K^{n+m}$, of
dimension~$n$;
\item for every $\x \in F$, the normal space $N_{\x}F$ to $F$ at $\x$ is
transversal to the coordinate space $\K^n$;
equivalently, the restriction to $F$ of the projection map $\Pi^{n + m}_n$ is
a local diffeomorphism between $F$ and~$\K^n$.
\end{enumerate}
\end{dfn}

For the remainder of this subsection,
$f: \K^n \leadsto \K^m$ is a $\CN$ \adm \pf, with graph~$F$.

\begin{dfn}
For every $C \subseteq \K^m$, denote by $f^{-1}(C)$ the preimage of $C$
under $f$, that is
$f^{-1}(C) := \set{\x \in \K^n: \exists \y \in C\ \pair{\x, \y} \in F}$.
Define $V(f) := f^{-1}(\set 0)$. Define the domain of $f$ to be $\dom(f) := f^{-1}(\K^m)$.
For every $A \subseteq \K^n$, denote by $f(A):= \set{\y \in \K^m: \exists \x \in A\ \pair{\x, \y} \in F}$, the image of $A$ under~$f$.
For every $\x \in \K^n$, we define $f(\x) := f\Pa{\set \x}$.
\end{dfn}

\begin{examples}\label{exa:adm}\
\begin{enumerate}
\item Every $\CN$ function is an admissible \pf.
\item The \pf $\sqrt{x}$ is not admissible.
\item
Define $g: \Real \leadsto \Real$ be the \pf with graph
$G := \set{\pair{x,y} \in \Real: y = x^2 \vee y = x^2 - 1}$.
$g$~is $\Cinf$ \adm, it is definable in the real field,
but it is not a partial function.
\item
Define $g: \Real \leadsto \Real$, $g(x) := 1/x$, defined for $x \neq 0$.
$g$ is an \adm $\Cinf$ partial function.
The domain of $g$ is not closed, and therefore it is not true that
the preimage of a closed set is closed.
\end{enumerate}
\end{examples}

\begin{lem}\
\begin{enumerate}
\item
For every $C \subseteq \K^m$ 
\dcompact, $f^{-1}(C)$ is closed \rom(in~$\K^n$\rom).
In particular, $V(f)$~is closed.
\item
For every $U \subseteq \K^m$ open and definable, $f^{-1}(U)$ is open.
In particular, $\dom(f)$~is open.
\end{enumerate}
\end{lem}
\begin{proof}
Let $x \in \cll{f^{-1}(C)}$.
We have to prove that $x \in f^{-1}(C)$.
Let $D := \bigl ( \cll{ F \cap (\K^n \times C)}  \bigr)_{x}$.
Notice that $D \subseteq \cll C = C$, and therefore $D = D \cap C$.
Since $x \in \cll{f^{-1}(C)}$, we have that
for every $U$ neighbourhood of $x$ there exists $y \in U$, such that
$f^{-1}(y) \cap C \neq \emptyset$, 
\ie the section
$( F \cap (\K^n \times C) \bigr )_{y} $ is non-empty.
Since $C$ is \dcompact, $D$~is non-empty.
Since $F$ and $C$ are closed,
we have $\cll{ F \cap (\K^n \times C)} = F \cap (\K^n \times C)$,
and therefore
\[
F_{x} \cap C = \bigl (F \cap (\K^n \times C) \bigr)_x = D.
\]
Since $D \neq \emptyset$, we have that $x \in f^{-1}(C)$.
\end{proof}


\begin{rem}
If $F$ is the graph of an \adm $\CN$ \pf, then every \dconnected component of
$F$ is the graph of an \adm $\CN$ \pf.
Conversely, if $F_1$ and $F_2$ are the graphs of 2 \adm $\CN$ \pfs,
$f_i : \K^n \leadsto \K^m$, and $F_1$ and $F_2$ are disjoint,
then $F_1 \cup F_2$ is the graph of an \adm $\CN$ \pf.
\end{rem}

\begin{lem}
Let $g: \K^n \to \K^m$ be a definable partial function.
Then, $g$~is \adm $\CN$ iff:
\begin{enumerate}
\item
the domain of~$g$ is an open set~$U$;
\item
$g: U \to \K^m$ is a $\CN$ function;
\item
for every $\x \in \partial U$,
\[
\lim_{\substack{\y \to \x,\\ \y \in U}} \abs{g(\y)} = + \infty.
\]
\end{enumerate}
\end{lem}

We conjecture that, if $F$ is \dconnected and $\dom(f) = \K^n$, then, $f$ is a (total and single-valued) function.


The reader can check that the following properties of \adm \pfs hold.

\begin{lem}\label{lem:composition-adm}\
\begin{itemize}
\item
Let $\phi: \K^m \to \K^m$ be a $\CN$, definable diffeomorphism.
Then, $\phi \circ f: \K^n \leadsto \K^m$ is $\CN$ and \adm.
\item
Let $\theta: \K^n \to \K^n$ be a $\CN$ definable diffeomorphism.
Then, $f \circ \theta: \K^n \leadsto \K^m$ is $\CN$ and \adm.
\item
Let $\theta: \K^n \to \K^n$ be a $\CN$ definable function.
If $f \circ \theta: \K^n \leadsto \K^m$ has a finite number of definably
connected components, then it is a $\CN$ \adm \pf.
\end{itemize}
\end{lem}

Notice that in the above lemma we can not drop the hypothesis that $\phi$
is a diffeomorphism, and replace it with the hypothesis
that it is a $\CN$ function, and similarly we cannot drop the additional
conditions on~$\theta$.
\begin{itemize}
\item
In fact, if $m = 1$, $n > 1$,  and $\phi(x) = x^2$, it might happen that
the graph of $\phi \circ f$ self-intersects. For example, let $g$ be defined
as in Example~\ref{exa:adm}(3).
Then the graph of $g^2$ self-intersects.
\item   For instance, let  $\K$ be an expansion of $\R$ where the sine function is defined; let $f(x) := 1/x$, and $\theta(t) := \sin t$.
Then, $f \circ \theta = 1/ \sin( t)$ is not \adm.
\end{itemize}

\begin{lem}\label{lem:difference-adm}
Let $m = 1$, and
define $g: \K^{n + 1} \leadsto \K$ as $g(\x,y) := y - f(\x)$.
That is, the graph of $g$ is
$G:= \set{\pair{\x,y,z} \in \K^{n+2}: \pair{\x, z - y} \in F}$.
Then, $g$ is $\CN$ and \adm.
\end{lem}

\begin{lem}\label{lem:extension-adm}
Given $g : \K^n \to \K$ a \rom(total and single-valued\rom)
$\CN$ and definable function,
define the \pf $h := \pairf{f,g}: \K^n \leadsto \K^{m+1}$;
that is, the graph of $h$ is
\[
H := \set{\pair{\x,\y,z} \in \K^{n+m+1}: \pair{\x,\y}\in F \et z =g(\x)}.
\]
Then, $h$ is $\CN$ and \adm.
\end{lem}

\begin{dfn}
For every $\pair{\x,\y} \in F$, it makes sense to define $Df(\x;\y)$, the
differential of $f$ at the point $\pair{\x,\y}$ (the notational difference
with the usual case when $f$ is a function is that here we have to specify at
which $\y \in f(\x)$ we compute~$Df$).  As usual, we say that $\pair{\x,\y}$
is a regular point for $f$ if $Df(\x;\y)$ has maximal rank, otherwise
$\pair{\x,\y}$ is singular.  Similarly, $\y \in \K^n$ is a regular value if,
for every $\x \in f^{-1}(y)$, $\pair{\x,\y}$ is a regular point; otherwise,
$\y$ is a singular value.
\end{dfn}

Moreover, we have a \pf on $\K^n$,
which assign to every point $\x$ the values of $Df(\x;\y)$,
as $\y$ varies in $f(\x)$.
This \pf in general is not \adm, even if $N \geq 2$, because its graph might
self-intersect.
The following lemma addresses this point.

\begin{lem}\label{lem:differential-adm}
Assume that $N \geq 2$.
\begin{itemize}
\item Let $\tilde Df$ be the \pf $\pairf{f, Df}$ on $\K^n$.  That is, the graph
of $\tilde Df$ is
\[
H := \set{\pair{\x,\y, \z}: \pair{\x,\y} \in F \et \z = Df(\x; \y)}.
\]
Then, $\tilde Df$ is $\Continuous^{N-1}$ and \adm.
\item
Assume that $n = m + k$, with $k \geq 1$.
Fix $1 \leq i_1 < \dots < i_k \leq n$.
Then, the \pf
\[
\pairf{f,\det\Pa{ \frac{\partial(f_1, \dotsc, f_k)}
  {\partial( x_{i_1}, \dotsc, x_{i_k} )} }^2}
\]
is \adm.
\end{itemize}
\end{lem}

The two previous lemmas are particular cases of the following:
\begin{lem}\label{lem:general-composition-adm}
Given $1 \leq M \leq N$, let $g: F \to \K^k$ be a $\CM$ function.
Let $h:= \pairf{f, g}$; that is, the graph of $h$ is
\[
H := \set{\pair{\x, \y, \z}: \pair{\x, \y} \in F \et z = g\pair{\x,\y}}.
\]
Then, $h$ is $\CM$ and \adm.
\end{lem}
\begin{proof}
Since $g$ is continuous, $H$ is closed in $F \times \K^k$.
Since $ F$ is closed in $\K^{n + m}$, $H$~is closed in
$\K^{n + m + k}$.
\end{proof}

\begin{lem}\label{lem:product-adm}
For $i = 1, 2$, let $f_i: \K^{n_i} \leadsto \K^{m_i}$ be an \adm $\CN$
\pf, with graph~$F_i$.
The, the \pf $f_1 \times f_2: \K^{n_1 + n_2} \leadsto \K^{m_1 + m_2}$,
with graph $F_1 \times F_2$, is an \adm $\CN$ \pf.
\end{lem}

\begin{dfn}
Given a \pf $g: \K^n \leadsto \K^m$, we denote by
$\abs g$ the \pf $\abs g: \K^n \leadsto \K$, with graph
$\abs G :=
\set{\pair{x,t}: \exists \y\in \K^m\ \pair{\x,\y} \in F \et \abs{\y} = t }$.
\end{dfn}

\begin{dfn}
Given $C \subseteq \K^n$ and
$g: \K^n \leadsto \K$ \pf with graph~$G$,
and $\x \in C$, we say that $g$ reaches the minimum on $C$ at~$\x$,
if there exists $y \in g(\x)$ such that, for every $\pair{\x', y'} \in G$,
if $\x' \in C$, then $y \leq y'$; moreover, $y$ is the minimum of $g$ on~$C$.\\
We also define
\[
\inf_{\x \in C} g(\x) := \inf g(C)
\in \K \sqcup \set{\pm \infty}.
\]
\end{dfn}
Notice that $\inf_{\x \in C} g(\x) = + \infty$ iff $g(C) = \emptyset$.

\begin{lem}\label{lem:compact-adm}
Let $f: \K^n \leadsto \K$ be \adm,
and $C \subseteq \K^n$ be definable, 
such that $f(C)$ is non-empty.
\begin{enumerate}
\item If $C$ is \dcompact, then $\abs f$ reaches the minimum on~$C$;
\item however, $C$ \dcompact does not imply that $\abs f$
reaches the maximum on~$C$;
\item if $\inf_{x \in C} f(\x) \neq - \infty$, then $f$
reaches the minimum on~$C$ \rom(and similarly for the maximum\rom).
\end{enumerate}
\end{lem}
\begin{proof}
The graph $\abs F$ of $\abs f$ is closed in $C \times [0, + \infty)$,
and $C$ is \dcompact; hence, $\pi(\abs F)$ is closed in $[0, + \infty)$,
where $\pi: C \times [0, + \infty) \to [0, + \infty)$ is the projection
onto the second coordinate.
\end{proof}

\subsection{Preliminary results}\label{preliminary}

In this subsection we develop some preliminary results (corresponding to the results in \cite[\S 2]{wilkie99}) about o-minimal weak structures. We fix for the rest of the subsection a semi-closed o-minimal weak structure $\s$.


\begin{lem}\label{charb e' F sigma} Let $A\in\tts$.
Then the following are equivalent:

\begin{enumerate}
\item $A$ has empty interior.
\item $\cl{A}$ has empty interior.
\item $\cl{A}$ is meager.
\item $A$ is meager.
\end{enumerate}
\end{lem}

\proof

We first observe that, $\tts$ being semi-closed, every set in $\tts$ is an $\Fs$-set. In particular, by \ref{lem:meager-F}, the implications $(1\Rightarrow 4)$ and $(2\Rightarrow 3)$ are proven.

The implication $(1\Rightarrow 2)$ can be proven as in \cite[Lemma 2.7]{max},
where we conclude by using \ref{lem:Fubini-converse} instead of
Fubini's Theorem, and the previous observation.

The other implications are obvious.
\qed

With similar modification of Maxwell's proof, one can prove the following.
\begin{lem}
Let $A \in \tts_{n + 1}$ and $A \subset \K^n \times \K_+$.
If $A$ has empty interior, then $\cl A _0 \subset \K^n$ also has empty interior.
\end{lem}

To obtain versions of \cite[Theorems 2.3 - 2.6]{wilkie99} (whose proofs are in \cite[\S 4]{max}) which hold in definably complete Baire structures, we need to reprove some of Maxwell's results.

\begin{lem}\label{Max 3.1}
Let $C\in\tts_{n+1}$ have empty interior, and $B_0\subset\K^n$ be an open box. Then there exist $p\in\N$ and an open box $B\subseteq B_0$ such that for each $\x\in B$ the fiber $C_{\x}$ contains exactly $p$ points.
\end{lem}

\proof The proof can be easily adapted from \cite[Lemma 3.1]{max}, using \ref{lem:Fubini-converse}, together with the o-minimality of $\tts$ and Lemma \ref{charb e' F sigma}.
\qed

\begin{prop}\label{Max 3.2}
Let $f: \K^n \to \K$ such that $\Gamma(f)\in\tts_{n+1}$.
Let $\Df:=\set{\x\in A:\ f\mathrm{\ is\ discontinuous\ at\ }\x}$.
Then $\Df\in\tts$ and $\Df$ has empty interior.
\end{prop}

\begin{proof}
The graph of $f$ is in $\tts$, and therefore it is an $\Fs$;
hence, by Lemma~\ref{lem:discontinuity}, $\Df$~has empty interior.

It remains to show that $\Df$ is in $\tts$.
Up to a change of coordinate in the codomain, we can assume that $f$ is
bounded.
Let $F \subset \K^{n + 1}$ be the closure of $\Gamma(f)$.
Then, since $f$ is bounded, we have that
$\Df = \set{x \in \K^n: \card{F_x} \geq 2}$.
The latter set is in $\tts$.
\end{proof}

\bigskip

Having Proposition \ref{Max 3.2}, the proofs of the results in
\cite[\S 4]{max} go easily through in our case as well, up to some minor or
obvious modifications,%
\footnote{In the proof of \cite[Lemma 4.12]{max} the \lq\lq
unbounded case\rq\rq (iv) will be treated by using, instead of sequences,
the function $g$, defined in \cite[Lemma 3.2]{max} as
$\Gamma(g):=\set{(\x,y)|\ \exists z((\x,z) \in \Gamma(f) \land z y = 1)
\vee((\x,0) \in \Gamma(f)\land y = 1)}$,
where $f := \frac{\partial f^*}{\partial x_i}$.}
%
%
and so does the proof of the following version of \cite[Theorem 2.7]{wilkie99}.

\begin{thm}[Sard's Lemma for $\Cone$ functions in $\tts$]\label{sard for tts}
Suppose that $n\geq m\geq 1$ and that $f:U\to\K^m$ is a $\Cone$ function in $\tts$, where $U\subset\K^n$ is open. Then the set of singular values of $f$ is in $\tts$ and has empty interior \rom(hence, it is meager in~$\K^m$\rom).
\end{thm}

\begin{cor}\label{sard for pf}The above statement still holds if $f$ is a $\Cone$ admissible \pf whose graph $F$ is in $\tts$.
\end{cor}
\proof
By the Implicit Function Theorem and the definition of $F$, every point $(\x,\y)\in F$ has a definable neighbourhood $U=U_1\times U_2\subset\K^n\times\K^m$ such that $U\cap F$ is the graph of a $\Cone$ function $f_U:U_1\to U_2$. By reducing $U$, if necessary, we can ensure that $U_1$ and $U_2$ are in $\tts$ (in fact, we can assume they are boxes), so that $f_U\in\tts$. We can apply Theorem~\ref{sard for tts} to $f_U$ and obtain that the set of its singular values is meager. Now, since the set $\Sigma_f$ of the singular values of $f$ is given by $\bigcup_U\Sigma_{f_U}$, we can apply Corollary \ref{cor:local-meager} to the projection $\pi:\K^n\times\K^m\to\K^m$ onto the second factor and obtain that $\Sigma_f$ is meager. As in the proof of Theorem~\ref{sard for tts}, it is clear that $\Sigma_f\in\tts$. \qed

\bigskip

We now turn our attention to \cite[Corollary 2.9]{wilkie99}, which provides the main tool for the approximation of the boundary of the projection of a set in $\tts$. We need some preliminary lemmas. The following is our version of \cite[Theorem 2.8]{wilkie99}, and its proof does not present difficulties.

\begin{lem}\label{wilkie 2.8}
Suppose that $n>m\geq 1,\ \ov{a}\in\K^m,\ g:\K^n\leadsto\K^m$ is a $\Cone$
\adm \pf in~$\tts$,
$h:\K^n\to\K$ is a $\Cone$ function in $\tts$, and that $\av$ is a
regular value of $g$.
Then there are at most finitely many $b\in\K$ such that $(\av,b)$ is a singular value of the \adm \pf $\pairf{g,h}$.
\end{lem}

\begin{lem}\label{lem:restriction-adm}Let $f: \K^n \leadsto \K^m$ be a $\Cone$ \adm \pf, with graph~$F$.
Let $b \in \K$ and $1 \leq i \leq n$.
Define
\[
\hat F :=\set{\pair{\x,\y} \in \K^{n-1} \times \K^m:
\pair{x_1, \dotsc, x_{i-1}, b, x_{i+1}, \dotsc, x_n, y_1, \dotsc, y_m} \in F}.
\]

$\hat F$ is the graph of a \pf $\hat f$.
Assume that $\hat F$ has a finite number of \dconnected components.
Then, 
$\hat f$ is $\CN$ \adm.
Moreover, given $\overline a \in \K^m$, if $\pair{\overline a,b}$ is a regular
value for $\pairf{f, \pi_i}$ (where $\pi_i: \K^n \to \K$ is the projection onto
the $i$th coordinate), then $\overline a$ is a regular value for~$\hat f$.
\end{lem}

We can now prove the analogue of~\cite[Corollary~2.9]{wilkie99}.
\begin{prop}\label{prop:wilkie2.9}
Let $n, k \geq 1$,
$f = \pairf{f_1, \dotsc, f_k}: \K^{n + k} \to \K^k$ be an \adm $\Cone$ \pf
in~$\tts$, and $V := V(f)$.
Suppose further that $0$ is a regular value of~$f$, and that $U$ is an open
ball in~$\K^n$ with the property that the set $X := V \cap \pi^{-1}(U)$
is non-empty and bounded, where $\pi:= \Pi^{n+k}_n$.
Then either \rom(i\rom) $\pi(X) = U$, or \rom(ii\rom)
there exists $\eta > 0$ and distinct
$i_1, \dotsc, i_k \in \Nat$ with $1 \leq i_1, \dotsc, i_k \leq n + k$ such
that
$\det \Pa{ \frac {\partial(f_1, \dotsc, f_k)}
  {\partial(x_{i_1}, \dotsc,  x_{i_k})} }^2 (\x; 0)$
takes all values in the interval $[0, \eta]$ on $X$.
\end{prop}
\begin{proof}
The proof proceeds as in the original~\cite[Corollary~2.9]{wilkie99}.
Since $f$ is \adm, $V$~is closed in $K^{n + k}$.

In the case when $\pi(X)$ is finite,
as in Wilkie's proof, one shows that there exists
$Y$ a \dconnected component of $V$ contained in~$X$.
Since $V$ is closed and $X$ is bounded, $Y$~is \dcompact.
Let $1 \leq i_1 < \dots < i_k \leq n + k$, $j \neq i_1, \dotsc, i_k$, and
$z \in Y$ be a point where the map $\x \mapsto x_j$ is maximal.
This clearly implies that
$\det \Pa{ \frac {\partial(f_1, \dotsc, f_k)}
  {\partial(x_{i_1}, \dotsc, x_{i_k})} } (z;0) = 0$, and one
concludes this case as in~\cite{wilkie99}.

In the case when $\pi(X)$ is infinite, let $Y$ be a \dconnected component of~$X$ (which is definable, since it is an atom of the finite boolean algebra formed by all definable clopen subsets of $X$).
It follows that either
$\det\Pa{\frac {\partial(f_1, \dotsc, f_k)}
 {\partial(x_{n+1}, \dotsc, x_{n + k})} (z;0)} = 0$
for some $z \in Y$, or $\pi(Y) = U$.
Using Lemma \ref{wilkie 2.8}
and Lemma~\ref{lem:restriction-adm},
we conclude as in Wilkie's proof.
\end{proof}


\subsection{Weakening the DSF condition and the theorem of the boundary}\label{subs:boundary}

Let $\s$ be a semi-closed o-minimal weak structure.

\begin{dfn}\label{dfn:DCN}
$\s$~satisfies \DCN for all $N$ if for each $A \in \s_n$, there exist
$m \geq n$ and $r \geq 1$,
such that, for each~$N$, there exists a set
$S_N \subseteq \K^{m}$, which is a \emph{finite union} of sets, each of
which is an intersection of at most $r$ sets of the form $V(f_{N,i})$,
where each  $f_{N ,i}: \K^m \to \K$ is an \emph{\adm{}} $\CN$ \pf in~$\tts$,
and  $A = \Pi^m_n(S_N)$
%
\end{dfn}

In the above definition, note that:
\begin{enumerate}
\item
Each set $S_N$ is of the form $S_N = \bigcup_{0 \leq j < k_N} S_{N,j}$
(for some natural number~$k_N$),
where each set $S_{N,j}$ is of the form
\[
S_{N,j} = \bigcap_{0 \leq i < r}V(f_{N, rj + i}).
\]
\item
If $S_N$ is an intersection of $r$ sets, each of which is a finite union of
sets of the form $V(f_{N,i})$
(where $f_{N ,i}: \K^m \to \K$ are \adm $\CN$ \pfs in~$\tts$),
then $S_N$ can be rewritten in a way to satisfy the conditions
in the above definition (with the same $\pair{m,r}$, and using the same
\pfs $f_{N,i}$).
\item
$m$ and $r$  do \emph{not} depend on~$N$;
however, the number of sets forming the union
(and therefore the total number of \pfs $f_{N,i}$) might depend on~$N$.
\item
We only ask the \pfs $f_{N,i}$ to be in~$\tts$, not in~$\s$,
and only that they are \adm \pfs, instead of total functions.
Thus, the condition above is weaker than the one formulated in~\cite{KM99},
even for $\K = \Real$. Moreover, $\s$ satisfying \DCN for all $N$ does
 not imply that $\s$ is semi-closed.
\item DSF implies \DCN for all~$N$.
\item If each $f_{N,i}$ is a (total single-valued)  function, we can replace
the functions $f_{N,i}$ by a single function~$f_N$,
obtained from the $f_{N,i}$ using products and sums of squares;
this is the reason why in \cite{KM99} only one function~$f_N$ is used
(and in \cite{wilkie99} one $\Cinf$ function~$f$).
However, for general \adm \pfs, we can not conclude that $f_N$ is \adm.
\item
Let $\s$ be a semi-closed o-minimal weak structure satisfying \DCN for all~$N$. Then it is harmless to assume $\s$ to be closed: if $\s '$ is the collection of all closed sets in $\s$, then $\s '$ is a closed o-minimal weak structure satisfying \DCN for all~$N$ and moreover the Charbonnel closures of $\s$ and $\s '$ coincide.
\end{enumerate}

\begin{thm}[After~\cite{KM99}]\label{thm:KM}
Suppose that $\s$ is a semi-closed o-minimal weak structure satisfying \DCN for all~$N$.
Then $\tts$ is o-minimal, and the smallest structure containing~$\s$.%
\end{thm}

Before proving the above theorem, we will give a lemma which is useful
in applications.

\begin{lem}\label{lem:DCN-PBC}
Let $A_1, \dotsc, A_l \in \s_n$ satisfy the condition for $A$ in
Definition~\ref{dfn:DCN} \rom(we say that $A_i$ satisfy \DCN for all~$N$\rom).
Then, also every finite positive Boolean combination \rom(PBC\rom)
of $A_1, \dotsc, A_l$ satisfies \DCN for all~$N$.
Hence, if $\s'$ is a subset of~$\s$, such that:
\begin{itemize}
\item
every set in $\s'$ satisfies \DCN for all~$N$, and
\item
every set in $\s$ is a PBC of sets in~$\s'$,
\end{itemize}
then $\s$ satisfies \DCN for all~$N$, and therefore, by Theorem~\ref{thm:KM},
$\tts$ is an o-minimal structure.
\end{lem}
\begin{proof}
It is clear that it suffices to prove the following: for every
$A^1, A^2 \in \s_n$, if each $A^i$ satisfies \DCN for all~$N$,
then $A^1 \cup A^2$ and $A^1 \cap A^2$ also satisfy \DCN for all~$N$.

For $i = 1, 2$, let $\pair{m_i, r_i}$ be the DC-complexity of~$A^i$,
and let $m := \max(m_1, m_2)$, $r := \max(r_1, r_2)$.
For each~$N$ and $i = 1, 2$, let
\[
S^i_{N} := S^i_{N,1} \cup \dots \cup S^i_{N,k_{i,N}},
\]
such that each $S^i_{N,j} \in \s_{m_i}$ is
an intersection of $r_i$ sets of the form $V(g)$ for some \adm
$\CN$ \pf $g: \K^{m_i} \leadsto \K$ in~$\tts$, and
$A^i = \Pi^{m_i}_n S^i_{N}$.
Notice that we can always assume that $m_1 = m_2 = m$ and $r_1 = r_2 = r$.
In fact, for each $g$ as above, define $\tilde g: \K^m \to \K$ as
$\tilde g(x_1, \dots, x_m) := g(x_1, \dotsc, x_{m_i})$, and substitute
$V(g)$ with $V(\tilde g)$ everywhere.

For the union, notice that $S^1_{N} \cup S^2_{N} =  S^1_{N,1} \cup \dots \cup
S^1_{N,k_{1,N}} \cup  S^2_{N,1} \cup \dots \cup S^2_{N,k_{2,N}}$, and
$A^1 \cup A^2 = \Pi^m_n (S^1_{N} \cup S^2_{N})$.

For the intersection, let $\Lambda := \set{\pair{\x, \x'} \in \K^m \times \K^m:
x_i = x'_i, 1 \leq i \leq n}$,
Notice that
\[
A^1 \cap A^2 =
\Pi^{2m}_n\Pa{ (S^1_{N} \times S^2_{N}) \cap \Lambda }
= \bigcup_{\substack{j \leq k_{1,N},\\j' \leq k_{2,N}}} \Pi^{2m}_n
 \Pa{S^1_{N, j} \times S^2_{N, j'}) \cap \Lambda}.
\]
Hence, \wloG $k_{1,N} = k_{2,N} = 1$, and
$S^i_{N} = \bigcap_{j \leq r} V(f^i_{N,j})$, for some $f^i_{N,j} : \K^m
\leadsto \K$ \adm $\CN$ \pf in~$\tts$, $i = 1, 2$.
Thus, by distributivity, it suffices to treat the case when
$S^i_{N} = V(f^i_{N,1})$, $i = 1, 2$.
Define $\tilde F^1_{N} := \set{\pair{\x,\x',\y}:
\pair{\x, \y} \in \Gamma(f^1_{N,1})}$, and
$\tilde F^2_{N} := \K^m \times \Gamma(f^2_{N,1})$.
Notice that $\tilde F^i_{N}$ is the graph of an \adm $\CN$ \pf
$\tilde f^i_{N} : \K^m \times \K^m \leadsto \K$ in~$\tts$.
Moreover, $\Lambda = V(q)$ for some polynomial $q: \K^m \times \K^m \to \K$.
Finally,
$(S^1_{N} \times S^2_{N}) \cap \Lambda =
V(\tilde f^1_{N}) \cap V(\tilde f^2_{N}) \cap V(q)$, and therefore
$A^1 \cap A^2 = \Pi^m_n \Pa{V(\tilde f^1_{N}) \cap V(\tilde f^2_{N}) \cap V(q)}$.
\end{proof}

\begin{proviso*}
We fix for the rest of the section a closed o-minimal weak structure $\s$
satisfying \DCN for all~$N$.
\end{proviso*}

We will prove the following result, corresponding to \cite[Theorem 3.1]{wilkie99}.

\begin{thm}\label{boundary-DCN} Let $A\in\tts_n$ be closed. Then there exists a closed set $B\in\tts_n$ such that $B$ has empty interior and $\bd(A)\subseteq B$.
\end{thm}

Notice that, even without the \DCN hypothesis, the following is true:
if $A$ is closed, then $\bd(A)$ has empty interior.
The missing information is whether $\bd(A)$ is in $\tts$ or not.

We will follow the outline of \cite[\S 3]{wilkie99}, but we will use \cite{BS04} for some definitions and proofs. The two approaches are equivalent, but we find the latter easier to read.

\begin{dfn}\
\begin{itemize}
\item Given $\x \in \K^n$, let $\abs{\x} := \max\set{\abs{x_1}, \dotsc, \abs{x_n}}$,
and $\norm{\x} := \sqrt{x_1^2 + \dots + x_n^2}$.
Notice that $\x \mapsto \norm{\x}^2$ is a $\Cinf$ function, and so is the function $\x\mapsto\frac{1}{1+\|\x\|^2}$.
\item
Given $A \subseteq \K^n$ and $\eps \in \K_+$, define the {\em
$\eps$-neighborhood $A^\eps$} of $A $ as the set $\{x \in \K^n \mid \exists y
\in A \; |x-y| < \eps\}$.
\item \rom(The quantifier ``for all sufficiently small''\rom)
We write $\allsmall \eps \phi$ as a shorthand for
$(\exists \mu) (\forall \eps < \mu) \phi$,
where $\mu, \eps$ are always assumed to range in $\K_+$.
If $\veps=(\eps_1,\ldots,\eps_n)$, then $\allsmall\veps$ is an abbreviation for $\allsmall\eps_1\ldots\allsmall\eps_n$.
\item \rom(Sections\rom)
Given $S \subseteq \K^n \times \K_+^k$ and given $\eps_1,
\ldots, \eps_k \in \K_+$, we define $S_{\eps_1, \ldots, \eps_k}$ as the set
$\{x \in \K^n \mid
(x, \eps_1, \ldots, \eps_k) \in S\}$.
\item Let $A \subseteq \K^n, S \subseteq \K^{n}\times\K_+^k$.
$S$ approximates $A$ from below \rom($S\leq A$\rom) if
\[\allsmall \eps_0\allsmall \eps_1 \ldots \allsmall \eps_k
( S_{\eps_1, \ldots, \eps_k} \subseteq A^{\eps_0}).\]
\item Let $A \subseteq \K^n, S \subseteq \K^{n}\times\K_+^k$.
$S$ approximates $A$ from above on bounded sets \rom($S\geq A$\rom)
if \[\allsmall \eps_0 \allsmall \eps_1 \ldots \allsmall \eps_k
( A \cap B(0, 1/\eps_0) \subseteq (S_{\eps_1,\ldots, \eps_k})^{\eps_0}).\]
\end{itemize}
\end{dfn}

The reader can check the the following holds.
\begin{lem}\label{lem:appr-and}
For $i = 0, 1, 2$, let $\veps_i$ and $\x_i$ be tuples of length $n_i$.
Let $\phi_1(\x_0, \x_1)$ and $\phi_2(\x_0, \x_2)$ be formulae,
with the shown free variables \rom(and possibly bounded variables and
parameters\rom).
Define
\[
\phi (\x_0, \x_1, \x_2) := \phi_1(\x_0, \x_1) \et \phi_2(\x_0, \x_2).
\]
Then,
\begin{multline*}
\K \models \Bigl(
\Pa{\allsmall \veps_0 \allsmall \veps_1 \phi_1(\veps_0, \veps_1) } \et
\Pa{\allsmall \veps_0 \allsmall \veps_2 \phi_2(\veps_0, \veps_2) }\Bigr)
\rightarrow\\
\rightarrow \Pa{ \allsmall \veps_0 \allsmall \veps_1 \allsmall \veps_2
\phi(\veps_0, \veps_1, \veps_2)}.
\end{multline*}
\end{lem}

We will show that, to obtain Theorem~\ref{boundary-DCN}, it is enough to prove the following:

\begin{prop}\label{6.11}For each $n \in \Nat$, $A \in \tts_n$,
and each $N \geq 1$, the following holds:

\rom($\Phi_N$\rom):\ \
There exists $k \geq 1$ \rom(the {\em $N$-complexity} of~$A$\rom)
and a set $S \subseteq \K^n \times \K^k_+$
\rom (an {\em $\s(N)$-approximant} of~$A$\rom)
which is a finite union of sets of the form
\[
\set{ (\x, \veps) \in \K^n \times \K^k_+:
 \exists \y \in \K^{k - 1} f(\x,\y) \ni  \veps}
\]
\rom(the {\em $\s(N)$-constituents} of~$S$\rom),
where $f: \K^{n + k - 1} \leadsto \K^k$ is \emph{\adm{}},
$\CN$ and in~$\tts$, such that $S$ both approximates
$\bd(\overline A)$ from above on bounded sets and approximates $\overline A$
from below.

\end{prop}
A set $S$ of the above form is called an {\em $\s(N)$-set}.

\begin{rem}\label{rem:appr-dimension}
Let $0 \leq d < k$, and $f: \K^n \times \K^d \leadsto \K^k$ be
an \adm $\CN$ \pf in~$\tts$.
Let $S(f) := \set{\pair{\x, \veps} \in \K^n \times \K_+^k:
\exists \y \in \K^d\ f(\x, \y) \ni \veps}$.
Then, $S(f)$ is of the form
$\set{\pair{\x, \veps} \in \K^n \times \K_+^k:
\exists \z \in K^{k-1} \tilde f(\x, \z) \ni \veps}$, where
$\tilde f :\K^n \times \K^{k-1} \leadsto \K^k$ is an \adm $\CN$ \pf
in~$\tts$; in particular, $S(f)$ is an $\s(N)$-constituent.
The graph of $\tilde f$ is
\[
\tilde F := \set{\pair{\x,\z, \w} \in \K^n \times \K^{k-1} \times \K^k:
\pair{\x, z_1, \dotsc, z_d, \w} \in F}.
\]
\end{rem}

The following statement corresponds to~\cite[Lemma 6.7]{BS04} (which is a
remark at the end of \cite[\S4]{wilkie99}).
\begin{lem}
Given $N\geq 1$, every $\s(N)$-set has empty interior.
\end{lem}
\begin{proof}
It suffices to show that each $\s(N)$-constituent $S$ has empty interior.
$S$~is of the form $\im(g)\cap(\K^n\times\K_+^k)$, where
\[\begin{array}{cccc}
g: & \K^{n + k - 1} & \leadsto & \K^{n+k}\\
 & \pair{\x,\y} & \mapsto  & \pair{\x, f(\x,\y)},
\end{array}\]
for some $f: \K^{n+k-1} \leadsto \K^k$ \adm, $\CN$ and in $\tts$.
Since $g$ is in~$\tts$, Theorem~\ref{sard for pf} implies that
the image of $g$ has empty interior.
\end{proof}

The following two statements correspond to \cite[Lemmas 3.3 and 3.4]{wilkie99}, and their proofs do not present particular difficulties: it is enough to use Lemma \ref{lem:Fubini-converse} and Lemma \ref{charb e' F sigma} instead of Fubini's Theorem and \cite[Theorem~2.1]{wilkie99}.

\begin{lem}\label{wilkie3.3}
Let $A\in\tts_n, S\in\tts_{n+k}$. Suppose that $S$ has empty interior and is an {$\s(N)$-approximant} for $A$. Then so is the section $\ov{S}_{\ov{0}}=\{\x\in\K^n|\
(\x,\ov{0})\in\ov{S}\}\in\tts_n$.
\end{lem}

\begin{lem}\label{wilkie3.4}
If $A\in\tts_k$ has empty interior, then $\allsmall\veps\  \veps \notin A$.
\end{lem}

The remainder of this subsection is devoted to the proof of Proposition~\ref{6.11}.
Theorem~\ref{boundary-DCN} follows immediately from the proposition and
Lemma~\ref{wilkie3.3}.

The proof of Proposition~\ref{6.11} follows the pattern
of~\cite[\S10]{BS04}; however, we need to prove some more intermediate steps, due to the fact that we are dealing with several, not just one, \pfs in Definition \ref{dfn:DCN}.

\begin{lem}[Union]\label{lem:appr-union}
Let $N, r, n \geq 1$, $A_1, \dotsc, A_r$ be subsets of~$\K^n$,
and, for $i= 1, \dotsc, r$, let
$S_i \subseteq \K^{n} \times \K_+^{k_i}$ be an $\s(N)$-approximant for~$A_i$.
Then, $A := \bigcup_i A_i$ has an $\s(N)$-approximant.
\end{lem}
\begin{proof}
We may suppose that all the $A_i$ have the same $N$-complexity~$k$;
then, $\bigcup_i S_i$ is an $\s(N)$-approximant of~$A$.
\end{proof}

\begin{lem}\label{lem:appr-abs}
Let $f: \K^n \leadsto \K$ be an \adm $\CN$ \pf,
and define $S := \set{\pair{\x,t} \in \K^n \times \K_+: \abs{f(\x)} \ni t}$.
Then $S$ approximates $\bd[V(f)]$ from above on bounded sets.
\end{lem}
\begin{proof}
Fix $\eps > 0$, and let $V := V(f)$.
Let $X := \bd(V) \cap \cB{0; 1/\varepsilon}$,
and $Y_t := X \setminus \Pa{\abs f^{-1}(t)^\varepsilon}$.
Note that $X$ and $Y_t$ are \dcompact.
Let
\[
P := \set{t\in \K: t > 0 \et Y_t \neq \emptyset }.
\]
Assume for contradiction that the conclusion is false.
This implies that $P$ has arbitrarily small elements, if we chose
$\varepsilon$ small enough.
Let $\x \in \acc_{t \to 0^+} Y_t$
($\x$~exists, because each $Y_t$ is contained in the \dcompact set~$X$),
and $U := B(\x; \varepsilon/2)$.
Note that $V$ is closed (because $f$ is \adm), and that $\x \in \bd(V)$.

By shrinking $\varepsilon$ if necessary, we may assume that there exists
$\delta > 0$, such that
$F \cap \Pa{U \times (-\delta, \delta)}$ is the graph of a $\CN$ function
$g: U \to (-\delta, \delta)$, such that $g(\x) = 0$.
Since $\x \in \bd(V)$, $\abs g$~assumes a positive value $\gamma$
on~$U$.
Since $U$ is \dconnected and $g$ is continuous,
$\abs g$~assumes all values in the interval $[0,\gamma]$ in~$U$.
Choose $t_0 \in P$ such that $Y_{t_0} \cap U \neq \emptyset$,
and $t_0 < \gamma$.
Since $t_0 < \gamma$, $U \cap \abs g^{-1}(t_0) \neq \emptyset$; therefore,
$U \subseteq \abs g^{-1}(t_0)^\varepsilon$, and thus $Y_{t_0} \cap U = \emptyset$,
a contradiction.
\end{proof}

\begin{lem}[Zero-set of \pfs]\label{lem:appr-V}
If $f: \K^n \leadsto \K$ is \adm, $\CN$ and in $\tts$, then its
zero set $V(f)$ has an $\s(N)$-approximant $S \in \tts_{n + 2}$.
\end{lem}
\begin{proof}
Define the following 2 sets $S_+$ and $S_-$:
\[
S_{\pm}:= \set{\pair{\x, \varepsilon_1, \varepsilon_2} \in \K^n \times \K_+^2:
1 + \norm{\x}^2 \leq 1/\varepsilon_1 \ \&\ f(\x) \ni \pm\varepsilon_2},
\]
and $S := S_+ \cup S_-$.
By Lemma~\ref{lem:extension-adm}, $\pairf{\pm f, \phi}$ are $\CN$ and \adm,
where $\phi: \K^{n + 1} \to \K$,
$\pair{\x,y} \mapsto (1 + \norm{\x}^2 + y^2)^{-1}$ (and in~$\tts$).
Thus, $S$ is an $\s(N)$-set.

We prove that $S$ approximates $V(f)$ from below, namely
\[\allsmall \varepsilon_0 \allsmall \varepsilon_1 \allsmall \varepsilon_2\
S_{\varepsilon_1,\varepsilon_2} \subseteq V(f)^{\varepsilon_0}.
\]
Let $K := \set{\x \in \K^n : 1 + \norm{\x}^2 \leq 1 /\varepsilon_1}$,
and $H := K \setminus V(f)^{\varepsilon_0}$.
Note that $K$ and $H$ are \dcompact,
and $S_{\varepsilon_1,\varepsilon_2} \subseteq K$.

\begin{claim*}
$\abs f$ has a positive minimum on~$H$,
if $f(H)$ is non-empty.
\end{claim*}
If not, then, by Lemma~\ref{lem:compact-adm}, there exists $\x \in H$ such
that $\abs{f(\x)} \ni 0$; however, this means that $\x \in V(f) \cap H$,
contradicting the definition of~$H$.

Thus, if we choose $\varepsilon_2$ smaller than the minimum of $\abs f$ on $H$
(or arbitrarily if $H$ is empty), then
$S_{\varepsilon_1,\varepsilon_2} \cap H = \emptyset$, and therefore
$S_{\varepsilon_1,\varepsilon_2}
\subseteq K \cap V(f)^{\varepsilon_0} \subseteq V(f)^{\varepsilon_0}$.

To prove that $S$ approximates $\bd\Pa{V(f)}$ from above we proceed as
in~\cite[Lemma~10.3]{BS04}, using Lemma~\ref{lem:appr-abs} instead
of~\cite[Lemma~10.2]{BS04}.
\end{proof}

\begin{lem}[Projection]\label{lem:appr-projection}
Let $N \geq 1$.
If $A\subseteq \K^{n + 1}$ has an $\s(N + 1)$-approximant $S \subseteq \K^{n + 1}
\times \K^k_+$, then there is an $\s(N)$-approximant $S' \subseteq \K^{n} \times
\K^{k + 1}_+$ for $\Pi^{n + 1}_n A \subseteq \K^n$.
\end{lem}
The drop in regularity in the above lemma from $N + 1$ to $N$ is due to the
fact that the definition of $S'$ involves the derivatives of the
functions defining~$S$.

\begin{proof}
Define $S_{\varepsilon_1,\dotsc,\varepsilon_k}$,
$S_{\varepsilon_1,\dotsc,\varepsilon_k}[\varepsilon_{k+1}]$, and
$S'_{\varepsilon_1,\dotsc,\varepsilon_{k+1}}$
as in the proof of \cite[Lemma~10.6]{BS04}.
More precisely, $S = S^1 \cup \dots \cup S^l$, where each $S^i$ is an
$\s(N)$-constituent.
Hence, each $S^i$ is of the form
$S^{i} = F^i \cap (\K^{n+1} \times \K^{k-1} \times \K_+^{k})$,
where each $F^i$ is the graph of some \adm $\Continuous^{N+1}$ \pf
$f^i: \K^{n+1} \times \K^{k-1} \leadsto \K^k$.
Define, $S_{\veps}[\eps_{k+1}] := \bigcup_i S^i_{\veps}[\eps_{k+1}]$,
where $\veps := \pair{\eps_1, \dotsc, \eps_k}$, and
$S^i_{\veps}[\eps_{k+1}]$ is the set of points $\x$ in $(f^i)^{-1}(\veps)$,
such that one of the following conditions is satisfied for some
$1 \leq i_1 \leq \dots \leq i_k \leq n + k$:
\begin{itemize}
\item either $1 + \norm{\pair{x_{n+1}, \dotsc x_{n+k}}}^2 = 1/\eps_{k+1}$,
\item or
$\det\Pa{\frac{\partial f^i}{\partial(x_{i_1}, \dotsc, x_{i_k})}(\x; \eps)}^2
= \eps_{k+1}$.
\end{itemize}
Finally, $S'$ is the set whose sections $S'_{\veps, \eps_{k+1}} \subseteq \K^n$
are given by:
\[
S'_{\veps, \eps_{k+1}} := \Pi^{n+1}_n S_{\veps}[\eps_{k+1}].
\]

By lemmas~\ref{lem:differential-adm} and~\ref{lem:extension-adm},
$S'$ is an $\s(N)$-set.

The fact that $S'$ approximates $\overline{\Pi^{n+1}_n A}$ from below follows as
in~\cite[Lemma~10.6]{BS04}.

It remains to prove that $S'$ approximates $\bd(\cll{\Pi^{n+1}_n A})$
from above on bounded sets.

Using Lemma \ref{wilkie3.4}, it is easy to see that $\allsmall \veps$, $\veps$
is a regular value of each $f^1, \dotsc, f^k$.
Fix $\eps_0 > 0$.
Let $X := \bd(\cll A) \cap \cB{0; 1 / \varepsilon_0}$;
note that $X$ is \dcompact.
Let $\x \in X$,
and $U$ be the open ball  of center $\x$ and radius $\varepsilon_0$.
Reasoning as in \cite[Lemma 10.6]{BS04}, and using
Proposition~\ref{prop:wilkie2.9} instead of \cite[Lemma~10.4]{BS04},
we see that
$\allsmall \eps_1 \dotsc \allsmall \eps_{k + 1}$
$U \subseteq {S'}_{\eps_1, \dotsc, \eps_{k + 1}}^{\eps_0}$.
Using Lemma \ref{lem:exchange}, we deduce that
$\allsmall \varepsilon_1 \dotsc \allsmall \varepsilon_{k + 1}$
$X \subseteq {S'}_{\varepsilon_1, \dotsc, \varepsilon_{k + 1}}^{\varepsilon_0}$,
which is the conclusion.
\end{proof}

\begin{lem}[Product]\label{lem:appr-product}
Let $n_1, n_2, k_1, k_2, N \geq 1$.
For $i = 1, 2$,  let $A_i \in \tts_{n_i}$,
such that $A_i$ has empty interior \rom(in~$\K^{n_i}$\rom).
Assume that each $A_i$ has an $\s(N)$-approximant
$S^i \subset \K^{n_i} \times \K_+^{k_i}$.
Then, $A_1 \times A_2$ has an $\s(N)$-approximant
$S \subset \K^{n_1 + n_2} \times \K_+^{k_1 + k_2}$.
Moreover, up to permutation of variables, $S = S_1 \times S_2$.
\end{lem}
\begin{proof}
\Wlog, each $S^i$ has only one $\s(N)$-constituent, that is, it is of the form
\[
S^i := \set{\pair{\x, \veps} \in \K^{n_i} \times \K^{k_i}:
\exists \y\in \K^{k_i-1}\ f_i(\x,\y) \ni \veps},
\]
for some $\CN$ \adm \pf
$f_i: \K^{n_i} \times \K^{k_i -1} \leadsto \K^{k_i}$, $i = 1, 2$.
Define
\begin{multline*}
S := \set{\pair{\x_1, \x_2, \veps_1, \veps_2}  \in
\K^{n_1} \times \K^{n_2} \times \K_+^{k_1} \times \K_+^{k_2}:\\
\exists \y_1 \in \K^{k_1 -1}\ \exists \y_2 \in \K^{k_2 -1}\
f_1(\x_1, \y_1) \ni \veps_1 \et f_2(\x_2, \y_2) \ni \veps_2};
\end{multline*}
By Lemma~\ref{lem:product-adm} and Remark~\ref{rem:appr-dimension},
$S$~is an $\s(N)$-set in $\tts_{n_1 + n_2 + k_1 + k_2}$.
Since each $A_i$  has empty interior, also the $\cll{A_i}$ have empty
interiors; therefore, $\bd(\cll{A_i}) = \cll{A_i}$,
and we have $A_i \leq S^i$, and $S^i \leq A_i$.
By Lemma~\ref{lem:appr-and}, $S$~approximates $A_1 \times A_2$.
\end{proof}

\begin{lem}[Linear intersection]\label{lem:appr-intersection-linear}
Given $N, n, k \geq 1$, let $A \in \tts_n$ have an $\s(N)$-approximant
$S \subset \K^n \times \K_+^k$, and suppose $Y$ is an $(n - 1)$-dimensional
affine subset of~$\K^n$;
suppose further that $\inter{\cll A} \cap Y = \emptyset$.
Then, there is an $\s(N)$\hyph{}approximant $S' \subseteq \K^n \times
\K_+^{k+2}$ for $\cll A \cap Y$.
\end{lem}
\begin{proof}
The proof of \cite[Lemma~10.8]{BS04} goes through
(with $\s(N)$\hyph{}approximants replacing $M(\s)$\hyph{}approximants),
using Lemma~\ref{lem:extension-adm} to
ensure that the set $S'$ is indeed an $\s(N)$-set.
\end{proof}

\begin{lem}[Small intersection]\label{lem:appr-intersection-small}
Let $n, k_1, k_2, N \geq 1$; define $M := N + n$.
For $i = 1, 2$, let $A_i$ be closed sets in~$\tts_{n_i}$.
Assume that each $A_i$ has an $\s(M)$-approximant
$S^i \subset \K^n \times \K_+^{k_i}$.
Assume moreover that each $A_i$ has empty interior.
Then, $A := A_1 \cap A_2$ has an $\s(N)$-approximant
$S \subset \K^{n} \times \K_+^{3n + k_1 + k_2}$.
\end{lem}
\begin{proof}
$A = \Pi^{2n}_n \Pa{(A_1 \times A_2) \cap \Delta}$,
where $\Delta$ is the diagonal of $\K^n \times \K^n$.
By Lemma~\ref{lem:appr-product}, $A_1 \times A_2$ has an
$\s(M )$-approximant in $\s_{2n + k_1 + k_2}$.
By hypothesis, $A_1 \times A_2$ has empty interior, hence we can apply
Lemma~\ref{lem:appr-intersection-linear} $n$ times, and therefore
$(A_1 \times A_2) \cap \Delta$ has an $\s(M)$-approximant in
$\s_{2n + k_1 + k_2 + 2n}$.
Finally, by Lemma~\ref{lem:appr-projection}, $A$ has an
$\s(M - n)$-approximant in $\s_{4n + k_1 + k_2}$.
\end{proof}


\begin{lem}\label{lem:appr-Vz}
Let $f: \K^n \leadsto \K$ be an \adm $\CN$ \pf in~$\tts$.
Let $A := V(f) \times \set 0 \subset \K^{n+1}$.
Then, $A$ has a $\s(N)$-approximant in $\s_{n + 4}$.
\end{lem}
\begin{proof}
Define
\begin{multline*}
S := \set{\pair{\x, z, \eps_1, \eps_2, \eps_3} \in \K^{n+1} \times \K_+^3:\\
1 +  \norm{\x}^2 \leq 1/\eps_1\ \et\ \abs{z}^2 \leq \eps_2\ \et\ 
f(\x) + z \ni   \eps_3}.
\end{multline*}
Notice that $S$ is an $\s(N)$ set (with only one component): in fact,
$S = \set{\pair{\x, z ,\eps_1, \eps_2, \eps_3}: \exists y_1 y_2\
1/(1 + \norm{\x}^2) + y_1^2 = \eps_1 \et
z^2 + y_2^2 = \eps_2 \et
f(\x) + z \ni \eps_3
}$.
Notice also that $A$ has empty interior.
We claim that $A \leq S$ and $S \leq A$, proving the conclusion.
For fixed $t > 0$, let
$K(t) := \set{\x \in \K^n: 1 + \norm{\x}^2 \leq 1/t}$.

\begin{claim}
$S \leq A$.
\end{claim}
\Ie, $\allsmall \eps_0 \allsmall \veps$ $S_{\veps} \subseteq A^{\eps_0}$,
where $\veps := \pair{\eps_1, \eps_2, \eps_3}$.
Let $\pair{\x,z} \in S_{\veps}$.

Define $I := [-\sqrt{\eps_2}, \sqrt{\eps_2}]$
and $H := K(\eps_1) \setminus V(f)^{\eps_0/4}$.
$I$, $H$ and $K(\eps_1)$ are \dcompact,
and $S_{\veps} \subseteq K(\eps_1) \times I$.

We claim that $\abs f$ has a positive minimum on~$H$, if $f(H)$ is non-empty.
Otherwise, by Lemma~\ref{lem:compact-adm}, there exists $\x \in H$ such that
$f(\x) \ni 0$, contradicting the definition of~$H$.
Let $\delta > 0$ be such minimum (or $\delta = 1$ if $f(H)$ is empty).

If we choose $\eps_3$ smaller than~$\delta$,
then $\x \in K(\eps_1) \setminus H$,
and therefore $\x \in V(f)^{\eps_0/4}$.
Now choose $\eps_2$ smaller than ${\eps_0}^2 /4$,
and obtain
\[
\pair{\x,z} \in V(f)^{\eps_0/4} \times [-\eps_0/4, \eps_0/4] \subseteq A^{\eps_0}.
\]

\begin{claim}
$A \leq S$.
\end{claim}
\Ie, $\allsmall \eps_0 \allsmall \veps$
$A \cap B(0; 1/\eps_0) \subseteq (S_{\veps})^{\eps_0}$.
Fix $\eps_0 > 0$, and choose $1 > \delta_1 > 0$ such that
$B(0; 1/\eps_0) \subseteq K(\delta_1)$,
Let $\delta_2 := \eps_0/2$.
For any $\eps_2$ such that $0 < \eps_2 < \delta_2$,
let $\delta_3 := \eps_2/2$.
Finally, choose any $\eps_3$ such that $0 < \eps_3 < \delta_3$.
Let $\y := \pair{\x,z} \in A \cap B(0; 1/\eps_0)$.
We prove that, for $\eps_0$ and $\veps$ chosen as above,
$\y \in (S_{\veps})^{\eps_0}$.
First, notice $z = 0$ and $\x \in V(f)$.
Let $\w := \pair{\x, \eps_3}$.
Notice that $\dist(\y, \y') = \eps_3 < \eps_0$,
and that $\w \in S_{\eps}$, and therefore $\y \in (S_{\eps})^{\eps_0}$.
Hence,
\[
\forall \eps_0 \exists \delta_1 \forall \eps_1 < \delta_1
\exists \delta_2  \forall \eps_2 < \delta_2
\exists \delta_3 \forall \eps_3 < \delta_3\
\Pa{A \cap B(0; 1/\eps_0) \subseteq (S_{\veps})^{\eps_0}}.
\qedhere
\]
\end{proof}

\begin{proof}[Proof of \ref{6.11}]
First, we prove the case when $A \in \s_n$.
Fix $N \geq 1$.
Let $M$ be large enough (how large will be clear from the rest of the proof).

By hypothesis, there exists $m \geq n$ and $r \geq 1$,
such that $A = \Pi^m_n(S_M)$, for some $S_M \subseteq \K^m$ of the form
$S_M = \bigcup_{0 \leq j < k_M} S_{M,j}$
where each set $S_{M, j}$ is of the form
\[
S_{M,j} = \bigcap_{0 \leq i < r}V(f_{M, i, j}),
\]
and each $f_{M,i,j}: \K^m \leadsto \K$ is a $\CM$ \adm \pf in~$\tts$.

Let $A_j := \Pi^m_n(S_{M,j})$.
If we prove that each $A_j$ satisfies~$(\Phi_N)$,
then, by Lemma~\ref{lem:appr-union}, $A$ also satisfies~$(\Phi_N)$.
Therefore, \wloG, $k_M = 1$,
\ie $S_M = \bigcap_{0 \leq i < r} V_{M,i}$, where
$V_{M,i} := V(f_{M, i})$
(where each $f_{M,i}: \K^m \leadsto \K$ is a $\CM$ \adm \pf in~$\tts$).
By Lemma~\ref{lem:appr-V}, each $V_{M,i}$ satisfies $(\Phi_M)$.
We need to prove that $S_M$ satisfies $(\Phi_{M'})$ (for a suitable~$M'$).
If all the $V_{M,i}$ were with empty interior, we could apply
Lemma~\ref{lem:appr-intersection-small}.
Otherwise, for every $i$, define
$W_{M,i} := V_{M,i} \times \set 0 \subset \K^{m+1}$.
By Lemma~\ref{lem:appr-Vz}, each $W_{M,i}$ has an $\s(M)$-approximant in
$\s_{m + 4}$; moreover, each $W_{M,i}$ has empty interior,
and therefore, by Lemma~\ref{lem:appr-intersection-small},
$W_M := \bigcap_i W_{M,i}$ has an $\s(M - (r - 1)(m + 1))$-approximant in
$\s_{(3\cdot 2^r-2)m + 4\cdot 2^r - 5}$.
Since $S_M = \Pi^{m+1}_m(W)$,
$S_M$~has an $\s(M - r m + m - r)$-approximant in $\s_{(3\cdot 2^r-2)m + 4\cdot 2^r - 5}$.
Finally, by Lemma~\ref{lem:appr-projection}, $A$~has an
$\s(M - r m + n - r )$-approximant in $\s_{(3\cdot 2^r-2)m + 4\cdot 2^r - 5}$.

The general case $A \in \tts_n$  can be proved as in \cite[\S10]{BS04},
using Lemma~\ref{lem:appr-projection} instead of \cite[Lemma~10.6]{BS04},
and Lemma~\ref{lem:appr-intersection-linear} instead of Lemma~10.8.
\end{proof}


\subsection{Cell decomposition}\label{subsec:cell dec}

We can proceed to prove Theorem~\ref{wilkie1.8} by a cell decomposition argument: for every $A\in\tts_n$, the ambient space $\K^n$ can be partitioned into finitely many sets $A_1,\ldots,A_N\in\tts_n$, such that $A$ (and hence its complement) is the union of some of the $A_i$s. We follow the outline of \cite[\S 4]{wilkie99}. The reader can refer to \cite[Definitions 4.1 and 4.3]{wilkie99}, where we replace $\R$ by $\K$.

Our aim is now to prove the analogue of the $\tts$-cell Decomposition Theorem 4.5 in \cite{wilkie99}.
Once established this result, we see that Theorem~\ref{wilkie1.8} follows easily, as explained in the remarks preceding the proof of \cite[Theorem~4.5]{wilkie99}.

There are three points in Wilkie's proof that do use
some reasoning whose translation in our context is not readily apparent.
We will examine them.

\begin{oclaim}
Assume that $A \in \tts_n$ has empty interior.
For each $i \geq 1$, consider the set
\[
A_i := \bigset{ \x \in C : \exists y_1, \dotsc, y_i\ (y_1 < \dots < y_i
 \wedge \bigwedge_{j = 1} ^i (\x, y_j) \in A)
}.
\]
Then each set $A_i$ lies in $\tts_n$, and $A_N$ has empty interior in $\K^n$
for some $N \geq 1$.
\end{oclaim}
\begin{proof}
We proceed as in \cite[Lemma~7.8]{BS04}.
The definition of $A_i$ implies immediately that $A_i \in \tts_n$.

Let $N := \gamma(A) + 1$, and fix $\x \in C$.
Note that if the fibre $A_x$ has cardinality greater or equal to~$N$, then it
has non-empty interior.

Since $A$ has empty interior, it is meager.
Therefore, by Lemma \ref{lem:Fubini-converse},
the set of those points $\x \in C$ such that $A_{\x}$
has non-empty interior is meager.
Thus, $A_N$ is meager, and hence it has empty interior.
\end{proof}

\begin{oclaim}
Let $C'$ be an open $\tts$-cell compatible with $\cl{A_1}$, \dots, $\cl{A_N}$,
$\tilde H$, $\tilde H_f$ and $\tilde H_g$, and such that
$C' \cap \cl{A_1} \neq \emptyset$.
Choose $k < N$ maximal such that $C' \cap \cl{A_k} \neq \emptyset$.
Then $C' \subseteq A_k$.
\end{oclaim}
\begin{proof}
As in Wilkie's proof, we conclude that $\tilde H$, $\tilde H_f$ and
$\tilde H_g$  are disjoint from~$C'$.

We have $C' \subseteq \cl{A_k}$, because $C' \cap \cl{A_k} \neq \emptyset$.
Consider a point $\x \in C'$.
Let $M$ be the cardinality of the fibre $A_{\x}$;
note that, by definition of~$k$, $M \leq k$.
Let $y_0 := f(\x)$, $y_{M + 1} := g(\x)$, and, for $1 \leq i \leq M$,
$y_i$ be the $i$-th point of $\K$ such that $(\x,y_i) \in A$.

If $1 \leq i \leq M$, since $\x \notin \tilde H$, we may find open
neighbourhoods $V_i$ of $\x$ in $\K^{n}$ and $J_i$ of $y$ in~$\K$,
such that for each $\x' \in V_i$ there is at most one $y' \in J_i$ such that
$(\x',y') \in A$.

Similarly, if $i = 0$ or $i = M + 1$ then, since $\x \notin \tilde H_f \cup
\tilde H_g$, we may choose $V_i$ and $J_i$ such that $(V_i \times J_i) \cap A =
\emptyset$.
Let $T := \set{y \in \K: (\x,y) \in A\ \&\ y \notin \bigcup_i J_i}$,
and $T' := \set{\x} \times T$.
Note that $T'$ is a compact subset of~$C$, that $A$ is a closed subset of $C$
disjoint from $T'$.
Hence, the distance between $T'$ and $A$ is some positive number $d > 0$.
Let $U := \bigcap_{i = 0}^{M + 1} V_i \cap \set{\x' \in C': \abs{\x' - \x} < d}$.

Therefore, for every $\x' \in U$,
\begin{equation}\label{eq:card}
\card{(\set{\x'} \times \K) \cap A} \leq
\card{(\set{\x} \times \K) \cap A} = M.
\end{equation}
We conclude as in \cite{wilkie99}:
as $\x \in \cl{A_k}$, we may choose $\x' \in U \cap A_k$ here,
from which it follows (using the maximality of $k$) that $M = k$.
Hence $\x \in A_k$ and the claim is justified.
\end{proof}

Thus, for each $i = 1, \dotsc, k$, we may define the function $f_i: C' \to K$
is $\tts$ by $f_i(\x) = y$ iff $y$ is the $i$-th point of $\K$ such that
$(\x,y_i) \in A$.

\begin{oclaim}
Each function $f_i$ is continuous.
\end{oclaim}
\begin{proof}
Let $\x \in C'$.
Let $U$, $V_i$ and $J_i$ be defined as in the proof of the previous claim, for
$i = 1, \dotsc, k$.
Let $\x' \in U$.
Note that, since we have equality in~\eqref{eq:card}, then, for every $i = 1,
\dotsc, k$, there is exactly one $y'_i \in J_i$ such that $(\x',y'_i) \in A$.
Note also that $y'_i = f_i(\x')$.
Fix $i$ such that $1 \leq i \leq k$, and fix $J$ neighbourhood of $y_i = f_i(\x)$.
In the construction of $V_i$ and $J_i$, we could have chosen $J_i$ such that
$J_i \subseteq J$, and then found a corresponding~$V_i$.
Proceeding in the construction, we see that, for every $J$ neighbourhood of
$f_i(\x')$, we can find $U$ neighbourhood of $\x$ such
that $f_i(U) \subseteq J$, which is equivalent to the definition of $f_i$
being continuous at~$\x$.
Since $\x \in U$ is arbitrary, the claim is proved.
\end{proof}

\section{Pfaffian functions}\label{section:pfaffian}

Khovanskii's results in \cite{khovanskii} show that any expansion of the real field with a Pfaffian chain of functions satisfies the hypotheses of \cite[Theorem~1.9]{wilkie99}.
Let $\K$ be a definably complete Baire structure. In this section we give an analog of Khovanskii's results, thus providing an example of a class of structures to which Theorem~\ref{o-minimality theorem} applies.

\begin{dfn}\label{pfaff def}
Let $f_1,\ldots,f_s:\K^n\to\K$ be definable and $\Cone$. We say that $(f_1,\ldots,f_n)$ is a {\em Pfaffian chain} if $\frac{\partial f_i}{\partial x_j}\in\K[\x,f_1,\ldots,f_i]$ for $i=1,\ldots,s$ and $j=1,\ldots,n$. A definable map $F=(F_1,\ldots,F_m):\K^n\to\K^m$ is {\em Pfaffian} if $F_1,\ldots,F_m\in\K[\x,f_1,\ldots,f_s]$ for some Pfaffian chain $(f_1,\ldots,f_s)$.

Consider polynomials $p_{ij}\in\K[\x,y_1,\ldots,y_i],\ q_k\in\K[\x,y_1,\ldots,y_s]$ such that
\begin{tabular}{lllr}
$\frac{\partial f_i}{\partial x_j}(\x)$&$=$&$p_{ij}(\x,f_1(\x),\ldots,f_i(\x))$&$\ \ \ i\leq s,j\leq n$\\
$F_k(\x)$&$=$&$q_k(\x,f_1(\x),\ldots,f_s(\x))$&$\ \ \ k\leq m$
\end{tabular}

The {\em complexity} of $F$ is the sequence of integers $(n,m,s,\deg q_k,\deg p_{ij}:\ i\leq s,\ j\leq n,\ k\leq m)$.
\end{dfn}

We prove the following.

\begin{thm}\label{pfaff is omin}
Let $\K$ be a definably complete Baire structure. Let $\mathcal F$ be a family of $\K$-definable Pfaffian chains and let ${\K}_{\mathcal F}$ be the reduct of $\K$ generated by $+,\cdot$ and $\mathcal F$. Then ${\K}_{\mathcal F}$ is o-minimal.
\end{thm}

\begin{cor}\label{axiom exp} Let $\R_{\exp}$ be the real ordered field with the exponential function. Then the following statements axiomatize a recursive subtheory of $Th(\R_{\exp})$ which is o-minimal.

\begin{itemize}
\item Axioms of ordered field.
\item Axioms ensuring that the models are definably complete and Baire.
\item $\forall x\ \exp'(x)=\exp(x)\ \&\ \exp(0)=1$.
\end{itemize}
\end{cor}

Analogous statements hold, for example, for the structures 
$\langle\R;\, +$, $\cdot$, $0$, $1$, $\exp$, $\sin\rest[0,1]\rangle$ and  $\langle\R;\ +\cdot,0,1,x^{\alpha}\rangle$ ($\alpha$~a real number).

\bigskip

To obtain Theorem~\ref{pfaff is omin} it is enough to show that ${\K}_{\mathcal F}$ satisfies the hypotheses of Theorem~\ref{o-minimality theorem}. Hence, it suffices to prove the following version of Khovanskii's Theorems (see \cite[Theorems~1 and~2]{khovanskii}):

\pagebreak[2]
\begin{thm}\label{khov1e2}\ 
\begin{enumerate}
\item Suppose $F:\K^n\to\K^n$ is Pfaffian. Then the number of regular zeroes
of $F$ is finite and can be bounded 
by a function of the complexity of~$F$.
\item Suppose $F:\K^n\to\K^m$ is Pfaffian. Then the number of definably
connected components of $F^{-1}(0)$ is finite and can be bounded
by a function of complexity of $F$.
\end{enumerate}\end{thm}

The fact that the bounds in the above theorem depend only on the complexity
imply, in particular, that they do not depends on the coefficients of the
polynomials in the Pfaffian chain,
or on other parameters in the definition of~$F$.
Moreover, the reader can verify that the explicit bounds given
in~\cite{khovanskii} continue to work in this context.

Before proceeding with the proof, we need to develop a version of Sard's Lemma
holding true in this context.

\subsection{The Sard property and Noetherian Differential Rings}\label{subs:sard}

In Subsection \ref{preliminary} we have seen that a strong version of Sard's Lemma holds true for
functions definable in an o-minimal weak structure (see \cite[Theorems 2.7 and 2.8]{wilkie99}).
In this subsection we will show a version of Sard's Lemma for functions belonging to a Noetherian differential ring (Theorem~\ref{sard for noeth rings}). The proof of the two mentioned statements are completely different from one another;
in particular, the proof of \ref{sard for noeth rings} is a quite simple modification of the classical argument for Sard's Lemma.

\begin{notation} Fix $n\in\Nat\setminus\set{0}$ and a definably connected
definable open set $U\subseteq\K^n$.
Let $\Cinf(U,\K)$ be the ring of definable $\Cinf$ functions from $U$ to~$\K$.
\end{notation}

\begin{dfn}\label{noetherian rings}
A ring $M$ with the following properties

\begin{itemize}
\item $M\subseteq \Cinf(U,\K)$;
\item $M$ is Noetherian;
\item $M$ is closed under partial differentiation;
\item $M \supseteq \K[x_1,\ldots,x_n]$.
\end{itemize}

is called a {\em Noetherian differential ring}.

If $G := (g_1,\ldots,g_k) \in M^k$, we denote by $V(G)$ the set of
zeroes of~$G$,
and by $V^{\reg}(G)$ the set of \emph{regular} zeroes of~$G$.
\end{dfn}

Generalities on Noetherian differential rings of functions over definably
complete structures can be found in \cite{servi-articolo}. In particular, we
will need the following result, which states that in a Noetherian differential
ring there are no flat functions.

\begin{prop}\label{noetherian} Let $M\subseteq \Cinf(U,\K)$ be a
Noetherian differential ring and let $0\not\equiv g\in M$.
Then for every $x\in U$ such that $g(x)=0$, there exist $k\in\Nat$ and a derivative $\theta$ of order $k$ such that $\theta g(x)\not=0$.
\end{prop}

Fix a Noetherian differential ring $M\subseteq \Cinf(U,\K)$.

\begin{rem} For $g_1,\ldots,g_k\in M$, the set
$V:=V^{\reg}(g_1,\ldots,g_k)$ is in~$\Fs$; in fact consider the following
closed definable subset of $U \times \K$:
\[
C:=\bigcup_{E(x)}\set{(x,y)\in U\times\K :
\ \bigwedge_{i=1}^k g_i(x)=0\land
  \det(E(x))y-1=0},
\]
where $E(x)$ ranges over all maximal rank minors of the Jacobian matrix of
$(g_1,\ldots,g_k)$ in~$x$.
Now, $V=\Pi^{n+1}_n(C)$; since $C$ is an $\Fs$  of $\K^{n + 1}$
and $\Pi^{n+1}_n$ is continuous, $V$~is also an~$\Fs$.
\end{rem}

In this subsection we prove the following version of Sard's Lemma:

\begin{thm}\label{sard for noeth rings}
Fix $k,m\in\Nat$, $k\leq n$. Let
\begin{itemize}
\item $H=(h_1,\ldots,h_{n-k})\in M^{n-k}$ and $V:=V^{\reg}(H)\not=\emptyset$;
\item $F=(F_1,\ldots,F_m)\in M^m$ and $f:=F\rest V:V\to\K^m$;
\item $\Delta_f \subseteq V$ be the set of singular points of~$f$, and
$\Sigma_f := f(\Delta_f)$
be the set of singular values of~$f$.
\end{itemize}
Then, $f:V\to\K^m$ has the Sard property, \ie $\Sigma_f$ is a meager set
\rom(in~$\K^m$\rom).
\end{thm}

\proof
We proceed by induction on $\dim V$ and $m$.
If $m = 0$, there are no singular points. If $\dim V = 0$, then $V$ is discrete. In particular, for every $a \in\Delta_f$ there exists $U_a$ neighbourhood of $a$
such that $\Delta_f \cap U_a = \set a$. Hence we can apply Corollary~\ref{cor:local-meager} and we are done.

Consider now the general case.

\begin{claim}\label{cl:V=K}
We can restrict to the case $V=\K^{k}$.
\end{claim}
By Corollary \ref{cor:local-meager}, it suffice to prove that for every ${a} \in \Delta_f$ there exists a neighbourhood $U_{{a}}$ of ${a}$ such that $f(U_{{a}}
\cap \Delta_f)$ is meager. Fix ${a} \in \Delta_f$.
Using the Implicit Function Theorem, it is easy to check that there is a
neighbourhood $U_{{a}}$ of ${a}$ and a definable diffeomorphism $\Phi:\K^k\to
V\cap U_{{a}}$ such that $H\circ\Phi\equiv 0$ and each $F_i\circ\Phi$ belong to a Noetherian differential ring $M'\subseteq \Cinf(\K^k,\K)$ (see \cite{servi-articolo} for the details).
Hence Claim~\ref{cl:V=K} is proved and we may assume that $f:\K^k\to\K^m$,
and $f\in M\subseteq \Cinf(\K^k,\K)$.

Let $X_0 := \set {{a} \in \Delta_f: D f ({a}) \neq 0}$, where $D f$ is the Jacobian matrix of~$f$.
We first prove that $f(X_0)$ is meager.

Again by Corollary~\ref{cor:local-meager}, it suffice to prove that for every
${a} \in X_0$ there exists a neighbourhood $U_{{a}}$ of ${a}$ such that $f(U_{{a}}
\cap X_0)$ is meager.

Fix ${a}\in X_0$.

\begin{claim}\label{cl:first-coordinate}
We may assume that $f(x)=(x_1,f_2(x),\ldots,f_m(x))$.
\end{claim}
In fact, since $D f({a}) \neq 0$, w.l.o.g. we can assume that $\partial {f_1}({a}) /
\partial x_1 \neq 0$ and ${a} = 0$.

Consider definable neighbourhoods $O$ and $\tilde O\subset\K^k$ of $0$, where the following map is a diffeomorphism:
\[
\begin{array}{cccc}
G:&O&\to&\tilde O\\
 &x&\mapsto&({f_1}(x),x_2,\ldots,x_{k}).
\end{array}
\]
Let $\Delta$ be the determinant of the Jacobian of $G$ and let $\hat{M}:=\set{g\circ G^{-1}|\ g\in M}\subset \Cinf(\tilde O,\K)$; then the ring $\tilde{M}:=\hat{M}[\Delta^{-1}]$
is clearly Noetherian and differentially closed; define $\tilde{f}:=f\circ
G^{-1}\in\tilde{M}$. Since $G$ is a diffeomorphism, it is enough to prove the
statement for $\tilde{M}$ and $\tilde{f}$, and Claim~\ref{cl:first-coordinate}
is proved.


For every $t\in\K$, consider the Noetherian differential ring
\[
N_t:=\set{g_t:=g(t,x_2,\ldots,x_k)|\ g\in M}\subset
\Cinf(\tilde O \cap\K^{k-1},\K).
\]
Let $f_t:\K^{k-1}\to\K^{m-1}$ be the map $((f_2)_t,\ldots,(f_m)_t)$.
By inductive hypothesis, the set  $\Sigma_{f_t}$ is meager in $\K^{m-1}$.
Moreover, $f(X_0\cap \tilde O) \cap (\set{t} \times \K^ {m - 1} ) \subseteq
\set{t} \times \Sigma_{f_t}$.
Hence $f(X_0\cap \tilde O)\subseteq D:=\set{(t,y)\in\K\times\K^{k-1}|\
  y\in\Sigma_{f_t}}$.
By what we have just observed, $T(D):= \set{t \in \K: D_t
  \text{ is meager}}$ is residual, because $D_t = \Sigma_{f_t}$,
hence by Lemma~\ref{lem:Fubini-converse},
$D$~is meager.
It follows by Corollary~\ref{cor:local-meager} that $f(X_0)$ is meager.

\bigskip

Now, let ${a}\in \Delta_f$ such that $D f({a}) = 0$, and let $P$ be the
least natural number such that there exists $i \leq m$ and a
derivative $\theta$ of order $P$ such that, if $g_{\theta} := \theta f_i$, then
$g_{\theta}({a}) = 0$ and $D g_{\theta}({a}) \neq 0$.
Such a $P$ exists by Proposition~\ref{noetherian}. Let $W_{\theta}:=V^{\reg}(g_{\theta})\subset \K^k$ (notice that the inclusion is strict, hence $\dim W_{\theta}<k$). Then there is a definable open neighbourhood $O$ of ${a}$ such that
\[
\Delta_f\cap O\subseteq\bigcup_{\mathrm{ord}(\theta)\leq P}W_{\theta}.
\]
Hence it is enough to prove that $f(\Delta_f\cap W_{\theta})$ is meager.
Let $h_{\theta} := f \rest W_{\theta}$. By inductive hypothesis, $\Sigma_{h_{\theta}}$ is meager.
Note that if $x \in W_{\theta}$ is a singular point for $f$, then $x$ is also a
singular point for~$h_{\theta}$; that is, $\Delta_f \cap W_{\theta} \subseteq \Delta_{h_{\theta}}$,
and we are done.

\qed


\begin{cor}\label{morse}
Let $F \in M^k$ and $G \in M$.
Define $X := \Vreg(F) \subseteq U$, and, for every $\av \in \K^n$,
$g_{\av} : X \to \K$ as $g_{\av}(\x) := G(\x) + \sum a_i x_i$.
Then, the set $A = \{(a_1,\ldots,a_n) \in \K^n: g_{\av}$
is not a Morse function\footnote{A definable $\Ctwo$ function $f$, from a
$\Ctwo$ $\K$-manifold to~$\K$, is a Morse function if, as in the classical
definition, every singular point of $f$ is nondegenerate.}
on $X\}$ is meager.
\end{cor}


\begin{proof}
We proceed as in~\cite{gull-pollack}.
\begin{claim}\label{cl:morse-open}
The lemma is true if $k = 0$, \ie if $X = U$.
\end{claim}
In fact, $\av \in A$ iff $-\av$ is critical values of $\nabla G$,
and we can apply Theorem~\ref{sard for noeth rings}.

By the Implicit Function Theorem, around every point $p \in X$
there exists an open definable neighbourhood~$U_p$,
such that the restriction of some $n - k$ of the coordinate functions
on~$\K^n$ (\wloG, the first $n - k$), constitute a coordinate system in~$U_p$;
let $V_p := \Pi^n_{n-k}(U_p)$ and $\phi_p : V_p \to U_p$ be the inverse map of
$\Pi^n_{n-k} \rest U_p$.
Let $\tilde M$ be the ring of functions on $V_p$ if the form $h \circ \phi$,
where $h \in M$: notice that $\tilde M$ is contained in some Noetherian
differential ring $M_p$ (see \cite{servi-articolo}).
Let $A_p = \set{(a_1,\ldots,a_n) \in \K^n: g_{\av}
\text{ is not a Morse function on } V_p}$.
Proceeding as in~\cite{gull-pollack}, using Lemma~\ref{lem:Fubini-converse}
instead of Fubini's theorem, and Claim~\ref{cl:morse-open} applied to
functions in the ring $M_p$, we see that $A_p$ is meager for every $p \in X$.
Since $A = \bigcup_{p \in X} A_p$, Corollary~\ref{cor:local-meager} implies
that $A$ is meager.
\end{proof}

\begin{rem}\label{pfaff is noetherian}
Note that if $(f_1,\ldots,f_s)$ is a Pfaffian chain, then the ring $\K[\x,f_1,\ldots,f_s]$ is a Noetherian differential ring. In particular, Theorem~\ref{sard for noeth rings} holds for functions in this ring.
\end{rem}

\subsection{Proof of Theorem \ref{khov1e2}}
\setcounter{oclaim}{0}
We will follow the outline of \cite{marker}.

We argue by induction on the length $s$ of the Pfaffian chain. If $s=0$ then $F$ is a polynomial map and the bound is
given by \cite[Proposition 11.5.4]{BCR98}.

Let $s>0$.
\begin{inductivehypothesis}\label{inductive hyp}
We suppose that, for all Pfaffian chains of length $\leq s-1$, the two statements of Theorem~\ref{khov1e2} hold true.
\end{inductivehypothesis}

We first prove the first statement of Theorem~\ref{khov1e2}.
Let $F:\K^n\to\K^n$ be Pfaffian with respect to a Pfaffian chain $(f_1,\ldots,f_s)$, with $F=(F_1,\ldots,F_n)$ and $F_i(\x)=q_i(\x,f_1(\x),\ldots,f_s(\x))$.

\begin{lem}\label{reduction}
There are Pfaffian maps $H:\K^{n+1}\to\K^n$ and $G:\K^{n+1}\to\K$ such that
\begin{enumerate}
\item $H$ has length $s-1$ and $G$ has length $s$.
\item $V(G)=V^{\reg}(G)$.
\item If $\ov a\in V^{\reg}(F) \subseteq \K^n$, then $\exists b \in \K$ such
that $(\ov a,b) \in V^{\reg}(H,G)\subseteq \K^{n+1}$.
\end{enumerate}
\end{lem}

\begin{proof}
Define $H_i(\x,y):=q_i(\x,f_1(\x),\ldots,f_{s-1}(\x),y)\ (i=1,\ldots,n)$ and $G(\x,y):=y-f_s(\x)$.
\end{proof}

Hence it is enough to bound the cardinality of $V^{\reg}(H,G)$.

\begin{dfn}
A definable continuous function $f: \K^d \to \K^{d'}$ is \emph{proper}
if the pre-image of every \dcompact set is \dcompact.
\end{dfn}
\begin{rem}
A definable continuous function $f: \K^d \to \K^{d'}$ is proper
iff $\lim_{\abs x \to \infty} \abs{f(x)} = + \infty$.
\end{rem}

\begin{lem}\label{proper}
We may assume that $H$ is proper.
\end{lem}
\begin{proof} Suppose $H$ is not proper. For all $r\in\K$, we define a proper Pfaffian map $Q^r:\K^{n+2}\to\K^{n+1}$ such that:
\begin{enumerate}
\item the length of $Q^r$ is $s-1$ and its complexity does not depend on $r$;
\item for all $(\ov a,b)\in V^{\reg}(H,G)$, there exist $r\in\K$ and $c\in\K$ such that $(\ov a,b,c)\in V^{reg}(Q^r,G)$.
\end{enumerate}
It follows that, if $\forall r\ \card{V^{\reg}(Q^r,G)}<N$, then $\card{V^{\reg}(H,G)}<N$.
The components of $Q^r$ as defined as follows: $Q_0^r(\x,y,z)=\sum_{i=1}^nx_i^2+y^2+z^2-r^2$; $Q^r_i(\x,y,z)=H_i(\x,y)$ for $i=1,\ldots,n$.
\end{proof}

\begin{lem}\label{regular}
We may assume that $V(H)=V^{\reg}(H)$ and that $V(H,G)=V^{\reg}(H,G)$.
\end{lem}

\begin{proof}
Suppose this is not the case.
For every $\bv \in \K^n$, we consider the Pfaffian proper map
$H_{\bv} := H - \bv$.
Let $B$ be the set of all $\bv \in \K^n$ such that
$V(H_{\bv}) = V^{\reg}(H_{\bv})$ and $V(H_{\bv},G)=V^{\reg}(H_{\bv},G)$.
By Theorem~\ref{sard for noeth rings}, $B$~is a co-meager subset of~$\K^n$.
Note that $H_{\bv}$ has length $s - 1$ and same complexity as~$H$.
Suppose we did prove that, for every $\bv \in B$,
$\card{V^{\reg}(H_{\bv},G)} \leq N$.
The set of all $\av\in\K^n$ such that $\card{V^{\reg}(H_{\av},G)}\geq N + 1$
is open (by the Implicit Function Theorem, applied to $H$ restricted to the
manifold $V(G)=V^{\reg}(G)$) and disjoint from~$B$, and therefore empty.
\end{proof}


We have thus reduced our problem to the following situation: $\Gamma := V^{\reg}(H) = V(H)\subseteq\K^{n+1}$ is a smooth \dcompact
Pfaffian curve of length $s-1$ and
$G:\K^{n+1}\to\K$ is
a Pfaffian map of length $s$ such that $G\rest \Gamma$ has only regular zeroes. We need to bound the number of such zeroes.

\begin{dfn} \label{arc}An {\em arc} of a non singular curve $\Gamma$ is the image of a differentiable
function $\phi \colon I \to \Gamma$ such that $I \subseteq \K$ is an interval and
$\phi'(t)$ is nonzero for all $t\in I$. The function $\phi$ is called a
{\em parametrization} of the arc. When no confusion is possible we use the word
``arc'' both for $\phi$ and its image.
\end{dfn}

\begin{lem} \label{atlas}
$\Gamma$ is the union of finitely many arcs.
\end{lem}
\begin{proof}
By the inductive hypothesis and Theorem~\ref{o-minimality theorem}.
\end{proof}

\begin{dfn}
Given a $\Cone$ function $f: \K^d \to \K^d$, let $J(f): \K^d \to \K$
be the determinant of the Jacobian matrix of~$f$.
\end{dfn}

\begin{dfn}
Let $\xi_H$ be the unique vector field on $\K^{n+1}$ such that for
every smooth definable function $g\colon \K^{n+1} \to \K$ we have $\xi_H(\x)
\cdot \nabla g(\x) = J(H,g)(\x)$. Note that $\xi_H$ is tangent to $\Gamma$ and is never zero
on $\Gamma$. We say that the arc $\phi\colon I \to \Gamma$ is {\em orientation
preserving} if $\phi'(t) \cdot \xi_H(\phi(t)) > 0$ for every $t\in I$. Note that if $\phi\colon (a,b) \to \Gamma$ is not orientation preserving, then
its reverse arc $-\phi (t) = \phi(b-t+a)$ is orientation preserving.
\end{dfn}

\begin{dfn} We say that two points $\x,\y \in V(H,G)$ are {\em consecutive} if
there are an orientation preserving arc $\phi \colon I \to \Gamma = V(H)$ and
$t_1<t_2$ in $I$ such that $\x = \phi(t_1), \y = \phi(t_2)$ and $\phi(t) \notin
V(G)$ for every $t \in (t_1, t_2)$. \end{dfn}

\begin{lem} \label{signs}
Let $\x,\y$ be consecutive points in $V(H,G)$. Then $J(H, G)$ assumes opposite
signs at $\x,\y$. So in particular $\x \neq \y$. \end{lem}
\begin{proof}
We are going to use the elementary fact that if a function $h \colon I \to \K$
defined on an interval $I \subseteq \K$ has two consecutive zeros $t_1<t_2$ in
$I$, and has nonzero derivative at these points, then $h'(t_1)$ and $h'(t_2)$
have opposite signs.

To reduce to this situation consider an orientation preserving arc $\phi \colon I \to \Gamma$ with $\x
= \phi(t_1) \in V(G), \y = \phi(t_2) \in V(G)$ and $\phi(t) \notin V(G)$ for every
$t \in (t_1, t_2)$. The derivative $\frac {\de (G \circ \phi) (t)} {\de t}$ equals
$\phi'(t) \cdot \nabla G (\phi (t))$, which has the same sign as $\xi_H(\phi (t))
\cdot \nabla G (\phi (t)) = J(H,G)(\phi (t))$ (this is nonzero since $J(H,G) \neq
0$ on $V(G)$). So if $J(H,G)$ assumes the same sign at $\x,\y$, then $G \circ
\phi \colon I \to \K$ would contradict the elementary fact stated above.
\end{proof}

\begin{lem} \label{next}
For each $\x \in V(H,G)$, there is $\y\in V(H,G)$ such that $\x,\y$ are
consecutive. \end{lem}
\begin{proof} Let $\Gamma$ be the union of the orientation preserving arcs $\phi_0,\ldots,\phi_k$. We can assume that this family of arcs is essential, \ie no arc $\phi_i$ is contained in the union of the remaining arcs. Suppose $\phi_0$
contains~$\x$. If this arc does not contain a consecutive point to $\x$, then
it cannot contain any points of $V(G)$ coming after $\x$. Let $\phi_1$ be the arc such that $\lim_{t \to \sup
I} \phi_0(t)\in \phi_1$ and $\lim_{t \to \inf
I} \phi_1(t)\in \phi_0$. There is only one such arc, because $\Gamma$ is a smooth curve and otherwise the Implicit Function Theorem would be violated. We prolong the arc $\phi_0$ by
attaching $\phi_1$ to it. Suppose that the arc $\phi_1$ contains no consecutive points to $\x$. If $\lim_{t \to \sup
I} \phi_1(t)\in \phi_0$ and $\lim_{t \to \inf
I} \phi_0(t)\in \phi_1$, then the orientation reversing arc $-\phi_0$ must contain a consecutive point to $\x$, or else
$\x$ would be consecutive to itself, contradicting Lemma \ref{signs}. Otherwise, let $\phi_2$ be the unique arc which contains $\lim_{t \to \sup
I} \phi_1(t)$. Notice that, again by the Implicit Function Theorem, it is not possible that $\lim_{t \to \sup
I} \phi_2(t)\in \phi_1$. We carry on attaching arcs with this procedure, until we either find a consecutive point to $\x$ or we find an arc $\phi_i$ such that $\lim_{t \to \sup
I} \phi_i(t)\in \phi_0$. In this case, by the argument above, the arc $\phi_0$ must contain a consecutive point to $\x$.
\end{proof}

\begin{lem}
There is a Pfaffian function $\hat{J}\colon \K^{n+1} \to \K$ of length $s-1$
which coincides with $J(H,G)$ on $V(G)$.
\end{lem}
\begin{proof}
Let $\hat{J}(\x,y)$ be such that $\hat{J}(\x, f_s(\x)) = J(H,G)(\x)$. 
\end{proof}

Define $\jhat$ to be the restriction of $\hat J$ to~$\Gamma$.
Note that $\jhat$ assumes opposite signs at two consecutive points $\x,\y$ of
$V(H,G)$.

\begin{lem}
$V(H,G)$ is finite, and we can compute a bound $N$ on its cardinality in terms
of the complexity of $H, G$.
\end{lem}
\begin{proof}
Let $\eps > 0$ be the minimum of the absolute value of $\jhat$ on the closed and
bounded set $V(H,G)$. 
Then $\jhat$ assumes every value between $-\eps$ and $+\eps$ between any two 
consecutive points $\x,\y$ of $V(H,G)$. 
By Theorem~\ref{sard for noeth rings} $\jhat$ has a regular value 
$t\in (-\eps, +\eps)$. 
Since $\hat J$ has length $\leq s - 1$, using the inductive hypothesis we can
compute a finite bound on the cardinality of $\jhat^{-1}(t)$. 
This is also a bound on $V(H,G)$ since we can associate
injectively to each $\x \in V(H,G)$ a point of $\jhat{-1}(t)$ lying in the arc
between $\x$ and the point consecutive to~$\x$ (which exists by Lemma~\ref{next}).
\end{proof}

Combining all the lemmas, we obtain a proof of the first statement of Theorem~\ref{khov1e2}. We now prove the second statement.

Let $F:\K^n\to\K^m$ be Pfaffian with respect to a Pfaffian chain $(f_1,\ldots,f_s)$.

We need some preliminary results.

\begin{dfn}
Let $\cof(\K)$ be the cofinality of $\K$. A {\em sequence} is a map $x:\cof(\K)\to\K^m$. If $\xseq$ is a sequence,
we say that $x_k \to l$ if for every neighbourhood $V$ of $l$ there exists $\mu < \cof(\K)$ such that $x_k \in V$ for every $k > \mu$.
We call $\xseq$ {\em infinitesimal} if $x_k \to 0$.
\end{dfn}

\begin{lem}\label{cont comp}
Let $F: \K^n \to \K$ be definable, continuous, proper and nonnegative, and $M \in \Nat$.
Suppose there is an infinitesimal nonnegative sequence $\epseq$ such that for every $k<\cof(\K)$, $F^{-1}(\eps_k)$ has less than $M$ def-connected components. Then $F^{-1}(0)$ has less than $M$ def-connected components.
\end{lem}

\proof
Since $F$ is proper, $F^{-1}(0)$ and $F^{-1}(\eps_k)$ are \dcompact.
Assume, for contradiction, that there exists a partition $\set{C_0, \dotsc, C_M}$ of $F^{-1}(0)$ into non-empty definable clopen subsets.
Let
\[\begin{aligned}
\delta &:= \frac{1}{3}\min_{i \neq j} d (C_i, C_j),\\
W &:=  \set{\x \in \K^n: d (\x, F^{-1}(0)) < \delta},\\
B_i & := \set{\x \in \K^n: d (\x, C_i) < \delta},\\
J_i & := F(B_i).
\end{aligned}\]
Note that $\delta > 0$, that the $B_i$s are open (in $K^n$) and disjoint, that $B_i \cap F^{-1}(0) = C_i$, and that $W = \bigsqcup_i B_i$.
We note that each $J_i$ is def-connected: consider
w.l.o.g. $i = 0$.
Let $\varepsilon \in J_0$, and let $\y \in B_0$ such that $F(\y) = \varepsilon$.
Let $\x \in C_0$ such that $d (\x,\y) < \delta$.
Note that the segment $[\x,\y]$ is contained in~$B_0$.
Since $[\x,\y]$ is def-connected, $[0,\varepsilon] = F\Pa{[\x,\y]}$ is also def-connected, and therefore $J_0$ is def-connected.

Let
\[ \begin{aligned}
\theta_i &:= \sup J_i\\
\eta_1 & := \min_i \theta_i.
\end{aligned}\]
We claim that
there exists $\eta_2 > 0$ such that $F^{-1}\Pa{[0,\eta_2)} \subseteq W$.
Let $D := F^{-1}\Pa{[0,1]} \setminus W$.
Note that $D$ is \dcompact, because $F$ is proper.
If $D = \emptyset$, we can define $\eta_2 = 1$.
Otherwise, $F(D)$ is \dcompact and non-empty.
Let $\eta_2 := \min D$.
Since $F^{-1}(0) \cap D = \emptyset$, we have that $\eta_2 > 0$.
Let $F(\x) < \eta_2$.
Then, $\x \notin D$, and therefore $\x \in W$.

Define $\eta = \min(\eta_1, \eta_2)$.
Therefore, for every $\varepsilon < \eta$, we have
\[\begin{aligned}
 F^{-1}(\varepsilon) &\subseteq \bigsqcup_i B_i\\
\varepsilon & \in \bigcap_i F(B_i).
\end{aligned}\]
Let $k < \cof(\K)$ such that $\eps_k < \eta$.
Since $F^{-1}(\eps_k)$ has at most $M$ def-connected components, we deduce that $F^{-1}(\eps_k) \cap B_i = \emptyset$ for at least one~$i$.
However, this contradicts $\eps_k \in \bigcap_i F(B_i)$.
\qed

We turn to the proof of the second statement of Theorem~\ref{khov1e2}.


\begin{oclaim}\label{cl:proper}
We may assume that $F$ is proper (the preimage of a \dcompact is \dcompact).
\end{oclaim}
\begin{proof}[Proof of Claim~\ref{cl:proper}]
For every $r\in\K_+$, define the proper map
\[
G_r(x_1,\ldots,x_{n+1})=(F(x_1,\ldots,x_n,x_1^2+\ldots x_{n+1}^2-r^2)).
\]
If we find a bound for the number of connected components of $G_r^{-1}(0)$ not depending on the parameter $r$, then the same number will be a bound valid for $F$.
\end{proof}

\begin{oclaim}\label{cl:m=1}
We may assume $m=1$ and $F\geq 0$.
\end{oclaim}
\begin{proof}[Proof of Claim~\ref{cl:m=1}]
We may replace $F$ by $\sum F_i^2$.
\end{proof}

\begin{oclaim}\label{cl:0reg}
We may assume that $0$ is a regular value for $F$.
\end{oclaim}

\begin{proof}[Proof of Claim~\ref{cl:0reg}]
Consider the function $F_{\eps} := F-\eps$,
for $\eps\in\K_+$. It follows from Theorem~\ref{sard for noeth rings} that the set of critical values of $F$ is meager, hence we can find an infinitesimal sequence $(\eps_n)_{n<\cof(\K)}$ such that $\eps_n$ is a regular value for $F$. If we find a bound which works for $F_{\eps_n}$, then by Lemma \ref{cont comp}, the same bound will work for $F$.
\end{proof}

\begin{oclaim}\label{cl:morse}
We may assume that $x_n$ is a Morse function of $V(F)$.
\end{oclaim}
\begin{proof}[Proof of Claim~\ref{cl:morse}]
By Corollary~\ref{morse}, we can choose $(a_1,\ldots,a_n)\in\K^n$ such that $a_n\not=0$ and $\sum a_ix_i$ is a Morse function on $V(F)$. Define $G(x_1,\ldots,x_n)=F(x_1,\ldots,x_{n-1},x_n-\frac{1}{a_n}\sum_{i=1}^{n-1}a_ix_i)$.
Then $G$ is proper, $0$ is a regular value of $G$, $x_n$ is a Morse function on $V(G)$ and a bound on the number of connected components of $V(G)$ will also work for $V(F)$.
\end{proof}

Once these four claims are established,
{
note that every non-empty clopen definable subset $C$ of $V(F)$ is
d-compact, and hence the function $x_n$ has at least one critical point
on~$C$%
\footnote{If $\dim V(F) > 0$, then $x_n$ has actually at least two critical
  points on~$C$.}; it follows by a standard argument that
the number of def-connected components of $V(F)$ finite and is bounded by
the number of critical points of $x_n$ on $V(F)$, if the latter is also finite.

We can then proceed as in \cite{marker}:
}%
a calculation shows that the critical points of $x_n$ on $V(F)$ are regular zeroes of the map $(F,\frac{\partial F}{\partial x_1},\ldots,\frac{\partial F}{\partial x_{n-1}})$, a bound on whose number is given by the first statement in Theorem~\ref{khov1e2}. This concludes the proof.

\section{Relative Pfaffian closure}\label{section:speissegger}

A consequence of Wilkie's Theorem of the Complement 
\cite[Theorem~1.9]{wilkie99} is that the structure generated by the real
ordered field together with {\em all} Pfaffian chains is o-minimal. 

Here we prove that the Pfaffian closure of an o-minimal structure {\em inside a
  definably complete Baire structure} is o-minimal. 
This will be obtained by proving that such a structure satisfies the hypotheses of Theorem~\ref{thm:KM}, and this is the reason why we proved such a general statement (of which \cite[Theorem 1.8]{wilkie99} and \cite[Theorem 1]{KM99} are special cases, even if in \ref{thm:KM} one considers only expansions of the real field).

\subsection{Preliminary results on o-minimal structures}
We will need the following results about o-minimal structures.

Let $\F$ be an o-minimal structure expanding a (real closed) field.
In this subsection, by ``definable'' we will mean ``definable with parameters
from~$\F$'', and, by ``cell'', ``cell definable in $\F$''.

\begin{prop}\label{prop:DM-closed}
For every $N \geq 1$ and every $Y \subseteq \F^{n}$ closed and definable
there exists $h: \F^{n} \to [0,1]$ definable and~$\CN$,
such that $Y = V(h)$. In particular, $\F$ is generated by its $\CN$ definable functions.
Moreover, if $Z$ is a closed definable subset of $\F^n$ disjoint from $Y$, then
we can also require that $Z = V(1-h)$.
\end{prop}
\begin{proof}
We can use \cite[Corollary~C.12]{DM96}, since the proof works also for o-minimal structures expanding any real closed field, not just $\R$.
\end{proof}

\begin{lem}\label{lem:closed-projection}
For every $X \subseteq \F^n$ definable there exists $Y \subseteq
\F^{n + 1}$, also definable, such that $Y$ is closed and $X = \pi(Y)$.
If moreover $X$ is a $\CN$ cell, then $Y$ can be also chosen to be a
\rom(closed\rom) $\CN$ cell of the same dimension as~$X$.
\end{lem}
\begin{proof}
Since $X$ is a finite union of $\CN$ cells, and projection commutes with
topological closure, it suffices to deal with the case when $X$ is a $\CN$ cell.
If~$X$ is closed, define $Y := X \times \set{0}$.
Otherwise, $\partial X$ is a closed non-empty set; let $h: \F^n \to [0,1]$
be definable and $\CN$ such that $\partial X = V(h)$.
Define
\[
Y := \pi^{-1}(\cl X) \cap \set{(\x,z) \in \F^{n+1}: z \cdot h(\x) = 1}.
\]
It is easy to see that $Y$ is a cell satisfying the conclusion.%
\end{proof}

\begin{lem}\label{retraction}
Let $Y \subseteq \F^n$ be a \emph{closed} $\CN$ cell.
Then, there exists a definable $\CN$ retraction $r : \F^n \to Y$.
\end{lem}
\begin{proof}
After a permutation of variables, \wloG $Y = \Gamma(f)$, for some definable
$\CN$ function $f: W \to \F^{n-d}$, where $W$ is an open cell in $\F^{d}$.
Let $U := W \times \F^{n-d}$ and define 
$r_0 : U \to Y$, $r_0(\z,\y) := (\z, f(\z))$. 
Notice that $U$ is an open neighbourhood of $Y$ and $r_0$ is a retraction.
Let $V$ be an open definable subset  of $\F^{n}$, such that $Y\subseteq V$ and
$\cl V  \subseteq U$.
By Proposition~\ref{prop:DM-closed}
there exists $h: \F^n \to [0, 1]$ definable and $\CN$ such that $Y = h^{-1}(1)$
and $\F^n \setminus V = h^{-1}(0)$.

Since $Y$ is a $\CN$ cell, there exists $\phi: Y \to \F^d$ definable $\CN$
diffeomorphism, with $d := \dim Y$.
\Wlog, we can assume that $\phi(0) = 0$.
For every $t \in \F$ and $\x \in Y$, define
\[
t * \x := \phi^{-1}\Pa{t \cdot \phi(\x)} \in Y.
\]
Define\vspace{-2ex}
\[
r(\x) :=
\begin{cases}
 0 & \text{ if } \x \notin U;\\
h(\x) * r_0(\x) & \text{ if } \x \in U.
\end{cases}
\qedhere\]
\end{proof}

\subsection{Expansions of o-minimal structures by total smooth functions}
In this subsection we generalize Theorem~\ref{o-minimality theorem} to the situation where $\K_{\cal F}$ expands an o-minimal structure. More precisely, let $\K$ be a definably complete Baire structure, $\Kt$~be an o-minimal reduct of~$\K$, expanding the field structure, and
$\Ffam$ be a family of total $\Cinf$ functions definable in $\K$. We assume
that $\Ffam$ is closed under permutation of variables, contains the
coordinate functions $(x_1, \dotsc, x_n) \mapsto x_i$, and that if $f \in
\Ffam$, then $(\x,y) \mapsto f(\x)$ is also in~$\Ffam$. Let $\Kt(\Ffam)$ be the reduct of $\K$ generated by $\Kt$ and $\Ffam$. We give necessary and sufficient conditions for $\Kt(\Ffam)$ to be an
o-minimal structure.

\begin{dfn}Let $\Gfam_0$ be the set of  all total continuous functions definable
in~$\Kt$, and $\Gfam$ be the set of functions of the form $h \circ f$, for
some $f : \K^n \to \K^m$ in $\Ffam^m$ and some $h: \K^m \to \K$ in~$\Gfam_0$
(notice that $\Gfam_0 \subseteq \Gfam$).

For every $n \in \Nat$, let $\s_n$ be the family of subsets of $\K^n$
of the form $V(g)$, for some $g : \K^n \to \K$ in $\Gfam$,
and let $\s := (\s_n)_{n \in \Nat}$.
\end{dfn}

\begin{thm}\label{thm:o-minimal-expansions}
$\Kt(\Ffam)$ is o-minimal iff, for every $X$ in~$\s$, $\gamma(X) < \infty$.
\end{thm}

\begin{proof}
Notice that $\s$ is a closed weak structure.
It is obvious that every set in $\s$ is definable in $\Kt(\Ffam)$.
Conversely, since $\Kt$ is o-minimal, Prop.~\ref{prop:DM-closed}
and the fact that $\Gfam_0 \subseteq \Gfam$ imply that
the structure generated by $\s$ expands~$\Kt$;
since moreover $\Ffam \subseteq \Gfam$,
$\s$~generates $\Kt(\Ffam)$.

Hence, by Theorem~\ref{thm:KM}, it suffices to show that $\s$ satisfies \DCN
for all~$N$.  That is, let $n \in \Nat$ and fix $A \in \s_n$.
It is enough to prove the following:\\
(*) There exists $m \geq n$, such that, for every $N \in \Nat$, $A$ is of the
form $\pi(V(g_N))$ for some $g_N: \K^m \to \K$ in $\Gfam$ and~$\CN$.

Let $g \in \Gfam$ such that $A = V(g)$.
Hence, $g = h \circ f$, for some $f : \K^n \to \K^m$ in $\Ffam^m$
and some $h: \K^m \to \K$ in $\Gfam_0$.
Let $h_N: \K^{m} \to \K$ be $\CN$ and definable in $\Kt$,
such that $V(h) = V(h_N)$ (the existence of $h_N$~is given by Prop.~\ref{prop:DM-closed}),
and define $g_N := h_N \circ f: \K^n \to \K$.
Note that $g_N$ is $\CN$ and in~$\Gfam$.  Note moreover that
\[
A = V(g) = f^{-1}\Pa{V(h)} = f^{-1}\Pa{V(h_N} = V(g_N),
\]
and we are done (in fact, we see that we can take $m = n$ in (*)).
\end{proof}




\subsection{Speissegger's theorem}\label{subsec:VRL-main}

We proceed to define a notion of relative Pfaffian closure. 
We recall that Speissegger's results in \cite{speissegger} concern expansions
of the real field.
Let $\Rz$ be an o-minimal expansion of the real field.
Let $U \subseteq \R^n$ be an open subset definable in~$\Rz$, and
$\omega$ be an $\Rz$-definable $\Cone$-form on~$U$ which is never~$0$.
A \intro{leaf} with data $(U, \omega)$ is a closed connected real submanifold
of $U$ of dimension $n-1$ that is orthogonal to $\omega$ at every point.
A \intro{Rolle leaf} (RL) is a leaf $L$ which moreover satisfies the condition:
if $\gamma:[0,1] \to U$ is a $\Cone$ curve with  end-points in~$L$, 
then $\gamma$ is orthogonal to $\omega$ in at least one point. 
Speissegger proved that if we add to $\Rz$ all Rolle leaves with data
definable in $\Rz$, then we still get an o-minimal structure.

We now generalize Speissegger's theorem to o-minimal structure outside the
real line.
We remark that the first results in this direction are due to Fratarcangeli in~\cite{frat}. 
However his definitions and methods are substantially different from ours (he follows \cite{speissegger} whereas we follow~\cite{KM99}) and the results he obtains are a special case of the main theorem in this section.
We will use a definition of ``Rolle leaves'' (which we call Virtual Rolle
Leaves) which is more complicated than the one in \cite{fr}, but which will
allow us to give in Section~\ref{sec:effective} 
effective bounds on a series of topological invariants of sets definable in the Pfaffian closure of an o-minimal
expansion of the real field, thus answering a question of Fratarcangeli
\cite[p.6]{fr}.

\begin{proviso*}
Let $\Kt$ be an o-minimal structure (expanding a field), and $\K$ be an
expansion of $\Kt$ that is definably complete and Baire.
For the rest of this section,
by ``connected'' we will mean ``definably connected'' (in~$\K$),
by ``connected component'' we will mean ``definably connected component'', and
by ``cell'' we will mean ``cell definable (with parameters) in~$\Kt$''.
\end{proviso*}
$\K$-manifolds (which we will simply call manifolds) have already been defined
(Def.~\ref{dfn:manifold}).

\begin{dfn}
Let $\omega = a_1 d x_1 + \dots + a_n d x_n$ be a definable $\Cone$
differential form, defined on some definable open subset 
$U \subseteq \K^n$, such that $\omega \neq 0$ on all~$U$.
A \intro{multi-leaf} with data $(U, \omega)$ is a
is a $\Cone$ manifold~$M$ contained in~$U$ and closed in~$U$, 
of dimension~$n-1$, such that
$M$ is orthogonal to $\omega$ at all of its points
(\ie, $T_aM = \ker(\omega(a))$, for every $a \in M$).
\end{dfn}
Compare the above with the definition of $\Kt$-leaf in \cite[5.2]{fr}, where
he asks that $M$ is connected.

We must now face the problem of generalizing the notion of Rolle leaf
to the context of definably complete Baire structures.

We let an \intro{arc} be a definable $\Cone$ 
map $\gamma: [0,1] \to \K^n$, such that  $\gamma'$ is always non-zero.

The most natural notion of generalized Rolle leaf would be the
following (\cf \cite[Remark p. 33]{fr}):

\begin{dfn}\label{ARL} An \intro{alternate Rolle leaf} (ARL) is a
connected multi-leaf $L$  with data $(U, \omega)$ which moreover satisfies the condition:
if $\gamma:[0,1] \to U$ is an arc 
in~$L$, then $\gamma$ is orthogonal to $\omega$ in at least one point.
\end{dfn}
Unfortunately it is not clear whether in definably complete Baire structures
definable $\Cone$ connected manifolds of dimension one are parametraizable as a finite union of arcs.
This fact creates an impediment, as will be clear later (see
Subsection~\ref{subsec:VRL-proof}), and forces us to modify this definition.

One could think of replacing the use of arcs with the use of connected manifolds of dimension 1 (see the definition of
Rolle leaf according to Fratarcangeli \cite[1.5]{fr}.
 
   The drawback of this choice is that it
   is not possible to express with a first-order formula the fact that
   a set is definably connected. However, for the application we have in mind
   (see Section \ref{sec:effective}) we need the definition to be first-order (in a sense which will be made precise later).

   Hence we will introduce the notion of \intro{Virtual Rolle
   Leaf} (VRL, see Definition \ref{VRL}), which has the advantage of being first
 order (as ARL is) and at the same time of involving the notion of
 manifold of dimension one, rather than that of arc. 

\bigskip

We are now ready to define the notion of relative Pfaffian closure:

\begin{dfn}\label{def:PK}
Inductive definition: for every $n \in \Nat$, let
$\K_{n+1}$ be the expansion of $\K_n$ to a language $L_{n+1}$ with a new
predicate for every VRL with $\K_n$-definable data.
Let $L^*=\bigcup_n L_n$ and define the {\bf relative Virtual
   Pfaffian closure of $\Kt$
  inside~$\K$}, denoted by $\mathcal V\PK $, as the $L^*$-expansion of $\K_0$ where every predicate is interpreted as
the corresponding Rolle leaf.
\end{dfn}

 Our aim is to prove the following version of Speissegger's Theorem:

\begin{thm}\label{speissegger's thm} Let $\K$ be a definably
  complete Baire structure and $\Kt$ be an o-minimal reduct of
  $\K$. Then $\mathcal V\PK$ is o-minimal.
\end{thm}

\subsection{Virtual Rolle Leaves}\label{subsec:VRL}

We will now give the precise definition of Virtual Rolle leaf. 
The idea is the following: unlike the definition of ARL, where we considered
all arcs $X$, in the definition of VRL we we consider closed manifolds $X$ of
dimension 1 (not necessarily connected) such that $X$ does not have compact
connected components. 
We want to find a first order condition on $X$ that implies a bound on
the number of connected components of $X$: 
every component has two ``end-points at infinity'' (see
Definition~\ref{virtual bd} below); hence, if we ask that $X$ has at most $2k$
end-points at infinity, we obtain that $X$ has at most $k$ connected
components. It remains to express the requirement that $X$ have no compact
components in a first order way: this is done by asking the existence of a
definable $\Cone$ function without critical points on~$X$.

Finally, the Rolle condition for a leaf $L$ is expressed by asking that for
any $X$ as above that intersects $L$ in a number of points which is greater
than the number of connected components of~$X$, there is a point
where $X$ is orthogonal to the 1-form defining~$L$.

\smallskip 

\begin{dfn}
A \intro{weak cell} of dimension $d$ is a $\Kz$-definable set $U\subseteq \K^n$ which
is diffeomorphic, via a $\Kz$-definable map $\phi_U$, to~$\K^d$.
For every $0 < t \in \K$, we define 
$U_t := \phi_U^{-1}(\set{x \in \K^d: \norm x = t})$. 
\end{dfn}
We consider the diffeomorphism $\phi_U$ as part of a weak cell: 
the same subset $U$ of $\K^m$ two different choices of diffeomorphisms
should be considered two different weak cells.
Notice that, for $0 < t \in \K$, $U_t$ is a compact manifold of
dimension $d - 1$.

\begin{dfn}
Let $U \subseteq \K^n$ be a weak cell.  We say that $X \subseteq U$ is a
\intro{twine} in~$U$ if $X$ is a $1$-dimensional $\Cone$ manifold, such that
$X$ is closed in~$U$.  We say that $X \subseteq U$ is a \intro{good twine}
in~$U$ if $X$ is a twine in~$U$ and moreover there exists a definable $\Cone$
function $\rho: X \to \K$ without critical points.
\end{dfn}

\begin{rem}\label{rem:clopen-algebra}
Let $X \subseteq \K^n$ be definable.
We denote by $\Bool(X)$ the Boolean algebra of definable clopen subsets of $X$.
$\Bool(X)$ is finite iff $\cc(X)$ (the number of connected
components of~$X$) is finite, and in that case
each connected component of $X$ is definable 
and an atom of $\Bool(X)$, 
and moreover $\card{\Bool(X)} = 2^{\cc(X)}$.

Moreover, for every $n \in \Nat$, the following are equivalent:
\begin{enumerate}
\item $\Bool(X) \leq 2^n$;
\item $\cc(X) \leq n$;
\item if $Y_1, \dotsc, Y_{n+1}$ are disjoint element of $\Bool(X)$,
then at least one of them is empty.
\end{enumerate}
\end{rem}

\begin{rem}\label{rem:twine-compact}
Let $U$ be a weak cell and $X$ be a twine in~$U$.
Let $\emptyset \neq Y \in \Bool(X)$.
Then, $Y$ is a also twine.
If moreover $X$ is good, then $Y$ is also good and not compact.
In particular, if $X$ is a good twine and $\cc(X) < \infty$, 
then  no connected component 
of $X$ is compact.
\end{rem}

\begin{dfn}\label{virtual bd}
Let $U$ be a weak cell and $X$ be twine in~$U$.
For each $0 < t \in \K$, let
$X_t := \set {x \in X \cap U_t:
X \text{ is transversal to } U_t \text{ at } x}$.
We denote by
\[
\vb_U(X) := \limsup_{t \to + \infty} \card{X_t} \in \N \cup \set{\infty},
\]
the \intro{virtual boundary} of~$X$.
\end{dfn}
Notice that $X_t$ is a $0$-dimensional manifold,
and hence $\card{X_t} = \cc(X_t)$.
Notice also that, unlike the number of connected components,
$\vb_U(X)$ can be defined with a first order formula.

\begin{lem}
Let $U$ be a weak cell and $X \subseteq U$ be a good twine in~$U$.
Assume that $X = X_1 \sqcup X_2$, where $\emptyset \neq X_i \in \Bool(X)$, 
$i = 1, 2$.
Then 
$\vb_U(X) = \vb_U(X_1) + \vb_U(X_2)$. 
\end{lem}

\begin{lem}
Assume that $\K$ is o-minimal.
Let $U$ be a weak cell and $X \subseteq U$ be a good 
twine in~$U$.
If $X$ is a connected, then $\vb_U(X) = 2$.
More generally, $\vb_U(X) = 2 \cdot \cc(X)$.
\end{lem}
\begin{proof}
It suffices to do the case when $X$ is connected.
Since $\K$ is o-minimal, $X$ is then
the image of some definable $\Cone$ function $\gamma: (0,1) \to \K$.
The conclusion follows from the o-minimality of~$\K$.
\end{proof}

\begin{lem}\label{lem:twine-components}
Let $U$ be a weak cell and $X$ be a good twine in~$U$.
If $X$ is non-empty, then $\vb_U(X) \geq 1$, and if moreover $\K$ is an
expansion of $\R$, then $\vb_U(X) \geq 2$.

Moreover, if $\vb_U(X)$ is finite, then $\cc(X) \leq \vb_U(X)$, and in
particular $X$ has a finite number of connected component, and each
component of $X$ is not compact.
\end{lem}

\begin{proof}
By Remark~\ref{rem:clopen-algebra},
it suffices to show that if $X$ is non-empty, then $\vb_U(X) \geq 1$;
the remainder follows from Remark~\ref{rem:twine-compact}.
Assume, for contradiction, that $\vb_U(X) = 0$.
Since $\vb_U(X) = 0$, there exists $R > 0$ such that, for every $t > R$,
$X$ meets $U_t$ only non-transversally.
Define $r : U \to \K$, $r(x) := \abs{\phi_U(x)}$,
let $U_{> R} := r^{-1}(R, + \infty) = \set{x \in U: \abs{\phi_U(x)} > R}$,
$Y := X \cap U_{> R}$ and $s := r \rest Y$; notice that $Y$ is open in~$X$.
Notice also that $s$ has only critical points (on~$Y$), and therefore $s$ is
locally constant.
Hence, for every $t > R$, $s^{-1}(t)$ is clopen in~$Y$, and therefore it is
open in~$X$.
By Remark~\ref{rem:twine-compact}, $X$ is not compact; thus, there exists 
$t_0 > R$ such that $Z := X \cap U_t = Y \cap U_t$ is non-empty, and, by what
we  have said before, it is clopen in~$X$.
Hence, by Remark~\ref{rem:twine-compact}, $Z$~is not compact.
However, $Z$ is closed and bounded in $\K^n$, contradiction.

The case $\K$ expanding $\R$ follows from the fact that each connected
component of $X$ (not necessarily definable!) is the image of some $\Cone$
function $f: (0,1) \to \R$, and some standard analysis.
\end{proof}

\begin{rem}
Let $X \subset \K^n$ be a weak cell of dimension~$1$.
Then, $X$ is a good twine in itself, and $\vb_X(X) = 2$.
\end{rem}

\begin{dfn}\label{VRL} A \intro{virtual Rolle leaf} (VRL) is a multi-leaf $L$
with   data $(U, \omega)$ which satisfies the following condition: for every
$n\in\N$,\ for every $V \subseteq U \times \K^n$ weak cell and every 
$X$ good twine in~$V$,
if $\card{X \cap (L\times\K^n)} > \vb_V(X)$, then $X$ is orthogonal to
$\pi^*(\omega)$ in at least one point, where $\pi^*(\omega)$ is the $1$-form
on $U \times \K^n$ induced by $\omega$ via the projection $\pi: U \times \K^n
\to U$.
\end{dfn}

With the notation of the above definition, if $\vb_V(X)$ is infinite, then the
premise is false, and therefore the condition is automatically satisfied
(for the given~$X$). Therefore, to verify whether $L$ is a VRL, we need to
check only the good twines $X$ such that $\vb_V(X)$ is finite.
Moreover, by Lemma~\ref{lem:twine-components}, such a good twine $X$ satisfies
$\cc(X) \leq \vb_V(X)$.
Therefore, if $X_1, \dotsc, X_m$ are the components of $X$, 
and $\card{X \cap  L} >  \vb_V(X)$, then for at least one $i$ we have
$\card{X_i \cap L} > 1$.

\subsection{O-minimality of \texorpdfstring{$\mathcal V\PK$}{VP(K0,K)}}\label{subsec:VRL-proof}

In this subsection we prove Theorem \ref{speissegger's thm}.
For this subsection, a Rolle leaf will be a virtual Rolle leaf.

\begin{dfn} Let $\RK=\set{(\RK)_n|\ n\in\N}$ be such that $(\RK)_n$ consists of
all the finite unions of sets
$
A \cap L_1 \cap \dots \cap L_k,
$
which we call {\em basic Rolle sets},
where $A \subseteq \K^n$ is $\Kt$-definable, and each $L_i$ is a Rolle leaf
with data $(U_i, \omega _i)$ in~$\Kt$.
\end{dfn}

We will show:

\begin{prop}\label{RK} $\RK$ is a semi-closed o-minimal weak structure,
satisfying \DCN for all~$N$.
\end{prop}

Since $\RK$ generates $\K_1$ in Def.~\ref{def:PK}, this, together with Theorem~\ref{thm:KM}, shows that
$\K_1$ is o-minimal; by applying inductively the same result to each $\K_n$,
we obtain a proof of \ref{speissegger's thm}.

We will prove Proposition~\ref{RK} via a series of lemmas.

\begin{lem}
$\RK$ is a weak structure.
\end{lem}
\begin{proof}
As in~\cite[Lemma~3]{KM99}. Notice that $\RK$ is closed under cartesian products by definition of VRL.
\end{proof}

\begin{rem}
Every basic Rolle set is the projection of another basic Rolle set, such that
all the open sets $U_i$ in the data are the same open set~$U$.
\end{rem}
\begin{proof}
As in~\cite[\P3.4]{KM99}.
\end{proof}

\begin{prop}\label{prop:fr}

Let $\Omega = (\omega_1, \dotsc, \omega_q)$ be a tuple of $\Kt$-definable
nonsingular 1-forms defined on some common open subset $U$ of~$\K^n$, and let
$A$ be a $\Kt$-definable subset of~$U$.
Then, there is a natural number $N$ such that, whenever $L_i$ is a VRL  of
$\omega_i = 0$ for each $i = 1, \dotsc, q$, then $A \cap L_1 \dots \cap L_q$
is the union of fewer than $N$ connected manifolds.
Moreover, $N$ does not depends on the parameters used in defining
$\Omega$, $U$, and~$A$ (and on the choices of the leaves~$L_i$).

\end{prop}
The proof of this proposition is in Subsection \ref{subs:frat}.
Notice the similarity with \cite[Theorem~1.7]{fr}.


\begin{prop}\label{prop:fr-54}
Let $U$ be a $\Kt$-definable open subset of~$\K^n$, and $\omega$ be a
$\Kt$-definable 1-form on~$U$, such that $\omega \neq 0$ on all~$U$.
Let $L$ be a multi-leaf  with data ($U, \omega)$.
Let $C$ be a definable connected $\Cone$ manifold of dimension at most
$n-1$ contained in~$U$, such that $C$ is orthogonal to $\omega$ at all of its
points.
Then, either $C$ is contained in $L$, or $C$ is disjoint from~$L$.
\end{prop}
\begin{proof}
\cite[Lemma~5.4]{fr}.
\end{proof}

\begin{lem}\label{lem:KM-WS6}
$\RK$ is semi-closed. 
\end{lem}
\begin{proof}
We use:
\begin{enumerate}
\item union commutes with projection;
\item the class of projections (from various~$\K^n$) of closed sets in $\RK$
is closed under intersections.
\end{enumerate}
It suffices to prove that any Rolle leaf $L\subseteq \K^n$ is the projection of a closed set
in~$\RK$. Let $(U,\omega)$ be the data (definable in~$\Kt$) of~$L$.
Do a $\Cone$ cell decomposition of~$U$.
It suffices to prove that, for each cell $E_i$ in the decomposition, $L \cap
E_i$ is the projection of a closed set in~$\RK$.

Let $U_i$ be an open cell of the decomposition, and $L_i := L \cap U_i$.
Consider a $\Cone$ closed cell $D_i \subseteq \K^{n + 1}$ such that $U_i =
\pi(D_i)$, and $D_i$ is of dimension~$n$
($D_i$~exists by Lemma~\ref{lem:closed-projection}).
Let $\tilde\omega := \pi^*(\omega)$, the 1-form on $U \times \K$ induced
by~$\omega$, and $\tilde L := \pi^{-1}(L)$; note that $\tilde L$ is a Rolle
leaf,  with data $(U\times \K,\tilde \omega)$.
Define $C_i := D_i \cap \tilde L$; $C_i$ is a basic Rolle set, closed in
$\K^{n + 1}$, and $\pi(C_i) = L_i$.

If instead $E_i$ is a $\Cone$ cell in the decomposition of dimension less
than~$n$, consider the $\Kt$-definable set of points in $E_i$ whose tangent
space (w.r.t. the manifold~$E_i$) is contained in $\ker(\omega)$.
Decompose this again into $\Kt$-definable connected submanifolds.
By Prop.~\ref{prop:fr-54},
any of these is either disjoint or contained in~$L$.
Hence, $L \cap E_i$ is a finite union of sets definable in $\Kt$, and hence is
itself definable in~$\Kt$, and thus projection of a closed set (in~$\Kt)$.
\end{proof}

\begin{lem}\label{lem:KM-WS5}
$\RK$ is an o-minimal weak structure.
\end{lem}
\begin{proof}
The conclusion  can be easily obtained from Proposition~\ref{prop:fr},
reasoning as in \cite[Corollary~2.7]{speissegger}.
\end{proof}




Hence, we can conclude that $\RK$ is a semi-closed o-minimal weak structure.
The last step is proving that $\RK$ satisfies \DCN for all~$N$, and hence is
an o-minimal \emph{structure}.
Notice that the following lemma does not imply Lemma~\ref{lem:KM-WS6}, because
the \DCN condition does not imply that a weak structure is semi-closed.

\begin{lem}\label{lem:KM-DCN}
$\RK$ satisfies \DCN for all~$N$.
\end{lem}
\begin{proof}
Proceed as in the proof of the preceding lemma.
Using Lemma \ref{lem:DCN-PBC}, we are reduced to prove:

(*)\ \ If $U \subseteq \K^n$ is open and definable in~$\Kt$,
$\omega$ is a $\Cone$ form, 
also definable in $\Kt$, and $L$ is a Rolle leaf with data $(U, \omega)$,
then there exists a natural number~$r$,
such that, for every $N \geq 1$, there is a set $S \subseteq \K^{n+1}$,
such that $S$ is a finite union of sets,
each of whose is an intersection of at most $r$ sets
of the form $V(f_{N,i})$, where each
$f_{N,i}: \K^{n + 1} \leadsto \K$ is a $\CN$ \adm \pf in $\RKt$,
$i = 1, \dotsc, l$, and  $L = \pi\Pa{S}$, where $\pi := \Pi^{n+1}_n$.
\smallskip

%

Note that the above claim is the \DCN hypothesis for~$L$, with $m = n + 1$.
By inspecting the following proof, 
the reader can easily verify that $r$ indeed does not depend on~$N$.

Fix~$N$.
Do a decomposition of $U$ into $\CN$ cells~$E_i$, such that on each open cell
$\omega$ is a $\CN$ form.
It suffices to prove (*) for each $L \cap E_i$.

\case 1
If $E_i$ is a cell of dimension less than~$n$, then, as in the proof of the
previous lemma, $E_i \cap L$ is definable in $\Kt$, and hence is the
projection of $V(f_N)$, for some $\CN$ function $f_N : \K^{n + 1} \to \K$
definable in~$\Kt$.

\case 2
If $E_i$ is an open cell, by Proposition~\ref{prop:fr},
$L \cap E_i$ has a finite number of connected components~$L_1, \dotsc, L_k$;
moreover, since $L$ is a manifold of dimension $n - 1$ and $E_i$ is open, each
$L_j$ is also a manifold of dimension $n - 1$, and moreover  it is
a Rolle leaf with data $(E_i, \omega)$;
hence, by substituting $U$ with~$E_i$,
\wloG we can assume that $U$ is an open cell.

Let $\omega := a_1 dx_1 + \dots + a_n dx_n$,
and $V_j := \set{\x \in U: a_j(\x) \neq 0}$, $j = 1, \dotsc, n$.
Note that $V_j$ is open and definable in~$\Kt$.
Decompose again $U$ into $\CN$ cells, in a way compatible with each~$V_j$.
For the non-open cells, proceed as in Case~1.
For the open ones, do the same trick as before, and
reduce to the case $a_n(\x)$ never~$0$ on~$U$, and therefore
we can assume that $a_n$ is the constant function~$1$.

Hence, $L$ is a closed subset of~$U$, and satisfies all conditions for being
the graph of an \adm $\CN$ \pf $l: \K^{n-1} \leadsto \K$,
except that $L$ might not be closed in~$\K^n$.
If $U = \K^n$, we can easily conclude as in \cite[Lemma~6]{KM99}.
Otherwise, we have more work to do.


Let $\theta := \Pi^n_{n-1}$,
$U' := \theta(U)$ be the basis of the cell~$U$,
$\phi': \K^{n-1} \stackrel{\sim}{\longrightarrow} U'$ and
$\phi: \K^n \stackrel{\sim}{\longrightarrow} U$ be $\Kt$-definable
$\CN$ diffeomorphisms, such that $\phi'\circ \theta = \theta \circ  \phi$.
Let $\tilde L := \phi^{-1}(L)$, and $\tilde \omega := \phi^*(\omega)$.
Then, $\tilde L$ is a Rolle leaf, with data $(\K^n, \tilde\omega)$.
Moreover, $\tilde L$ is the graph of a $\CN$ \adm \pf
$\tilde l: \K^{n-1} \leadsto \K$
(in fact, $\tilde L$ is closed in~$\K^n$).

Define $\tilde g(x_1, \dotsc, x_n) := \tilde l (x_1,\dotsc, x_{n - 1}) -  x_n$,
$\tilde g: \K^n \leadsto \K$.
By Lemma \ref{lem:difference-adm}, $\tilde g$ is \adm;
notice that $\tilde L = V(\tilde g)$.

We would like to pullback $\tilde g$ via $\phi$;
the problem is that $\tilde g \circ \phi^{-1}$ is not defined on all~$\K^n$.

Let $D \subseteq \K^{n+1}$ be a \emph{closed} $\CN$ cell,
such that $\pi(D) = U$, 
$f_{N,1}: \K^{n+1} \to \K$ be a $\CN$ and $\Kt$-definable function,
such that $D = V(f_{N,1})$,
and $r: \K^{n+1} \to D$ be a $\Kt$-definable $\CN$ retraction
($D$, $f_{N,1}$ and $r$ exists by Lemmas \ref{lem:closed-projection} and \ref{retraction}).
Let $f_{N,2} := \tilde g \circ \phi^{-1} \circ \pi \circ r$.
Notice that $\phi^{-1} \circ \pi \circ r$ is a total $\CN$ function,
and therefore, by Lemmas~\ref{lem:composition-adm} and~\ref{lem:KM-WS5},
$f_{N,2}$~is \adm; $f_{N,1}$ is also obviously \adm.
It is also clear that $L = \pi\Pa{ V(f_{N,1}) \cap V(f_{N,2}) }$.
\end{proof}

\begin{rem}\label{KM sbagliano}
In the proof of \cite[Lemma 6]{KM99} there is a gap, in that
$f_{N,2}$~might not be a total function:
this is \emph{the} reason why we had to work with \adm \pfs
instead of total functions.
It is still true that the proof contained in \cite{KM99} implies the
o-minimality of the closure under {\em total} $\Cinf$
$R$-Pfaffian functions of an o-minimal expansion $R$ of the real field
(cf.\ Theorem~\ref{thm:o-minimal-expansions}).
\end{rem}


However, the correspondences under consideration are single valued, due to the Rolle condition.
\begin{lem}
Let $\omega := a_1 \de x_1 + \dots + a_n \de x_n$ be a $\CN$ $1$-form
on~$\K^n$, such that $a_n \equiv 1$, and let $F$ be a Rolle leaf for~$\omega$.
Then, $F$ is the graph of a $\Continuous^{N+1}$ partial function $f: U \to \K$,
with open domain $U \subseteq \K^{n-1}$.
\end{lem}
\begin{proof}
The fact that $f$ is an \adm $\Continuous^{N+1}$ \pf is clear.
It remains to prove that $f$ is single-valued.
If not, there exist $\x \in \K^{n-1}$ and $y_1 < y_2 \in \K$,
such that, for $i = 1,2$, $p_i := \pair{\x, y_i} \in F$.
Let $J$ be the ``vertical'' segment with endpoints $p_1$ and~$p_2$.
By the Rolle condition, there exists $q \in J$ such that $J$ is orthogonal
to $\omega$ at~$q$.
Since $J$ is vertical, this means that $\omega(q)$ is ``horizontal'',
contradicting the fact that $a_n \equiv 1$.
\end{proof}

However, as we said before, the partial function $f_{N,2}$ might not be
total, as the following example shows.
\begin{example}
Let $f: \Real \hookrightarrow \Real$
be the partial function $f(x) := 1/x$, defined
on~$\Real_+$, and $F \subset \Real^2$ be the graph of~$f$.
Let $\omega(x,y) := y^2 \de x + \de y$ be a $1$-form defined on~$\Real^2$.
Then, $F$~is a $\Cinf$ Rolle leaf of $\omega = 0$.
In fact, $f$~solves the differential equation $f' = -f^2$, and therefore we
can apply \cite[Example~1.3]{speissegger}.
\end{example}

\subsection{Proof of Proposition \ref{prop:fr}}\label{subs:frat}

We will assume familiarity with~\cite{fr}.
Some important but easy observations are the following ones:
\begin{asparaitem}
\item \cite[Lemma 5.9]{fr} does not require that the manifolds $L_i$ are
connected, and therefore can be applied to $L_i$ multi-leaves.

\item \cite[Prop.~5.10]{fr} does not use neither the conditions that the $L_i$
are connected nor the Rolle condition, and remains true for $L_i$ multi-leaves.

\item \cite[Prop.~5.7]{fr} can be used in the following form:

\begin{prop}
Let  $U$ and $V$ be definable open subsets of $\K^n$, 
and let $\sigma: V \to U$ be a definable diffeomorphism.
Let $\omega$ be a definable $1$-form on $U$, and $L$ be a multi-leaf with data
$(U, \omega)$.
Then, $\sigma^{-1}(L)$ is a multi-leaf with data 
$(V, \sigma^*(\omega))$.\\
If $L$ is a VRL and $\sigma$ is $\Kz$-definable, 
then $\sigma^{-1}(L)$ is a VRL.
\end{prop}
\end{asparaitem}
Hence, the Rolle condition is used directly only at the end of the proof, on
\cite[p.~39]{fr}. We will show how to use the Virtual Rolle condition.

The proof will proceed by induction on~$q$.
If $q = 0$, the conclusion follows from o-minimality of $\Kt$; hence, we can
assume $q \geq 1$.

For the inductive step, we assume that we have already proved the conclusion for
$q - 1$: that is, we assume that we have proved the
result for every $(q - 1)$-tuple $\Omega'$ of $\Kt$-definable nonsingular
1-forms defined on some open set $U'$ of $\K^{n'}$, for every $\Kt$-definable
set $A' \subseteq U'$, and for every
corresponding $(q-1)$-tuple of VRL with data $(U', \Omega')$.

Fix $U$, $\Omega$, $L_1, \dotsc, L_q$, and $A$ as in the assumption of the
theorem.
Let $d := \dim(A)$.
We prove the conclusion by a further induction on~$d$.

As  in the proof of~\cite[Theorem~1.7]{fr}, we can reduce to the case when
$A$ is a $\Kt$-definable $\Cone$-cell of dimension $d \geq q$, contained in~$U$,
and $\Omega$ is transverse to~$A$; that is, for every $a \in A$,
the projections of (the vector fields associated to)
$\omega_1, \dotsc, \omega_q$ on $T_a(A)$ are linearly independent.
Notice that ``$\Omega$ transverse to $A$'' is equivalent to
``the projection on $T(A)$ of the $q$-form $\omega_1 \wedge \dots \wedge
\omega_q$ is never null''.

If $d > q$, we can conclude by induction on $d$ as in \cite[p~.39,
``\textsc{Case} $d > q$'']{fr};  as we noticed before, the
Rolle condition is not used in \cite[5.10]{fr}, and therefore we can use it in
our situation.

Hence, it remains to treat the case $d = q$.

If $d = q$, we treat first as a way of exemplification the case $d = q = 1$.
Then, $A$ is a good twine in itself, thus, by the Rolle condition, and the
fact that $\omega_1$ is transverse to~$A$,
$\card{A \cap L_1} \leq \vb_A(A) = 2 $, and we are done.

In general, if $d = q$,
define $L' := A \cap L_1 \dotsc L_{q - 1}$ (or $L' := A$ if $q = 1$).
Notice that $L'$ is a twine in~$A$.
Let $\omega' := \omega'_1 \wedge \dotsc \wedge \omega'_{q-1}$, where each
$\omega'_i$ is the projection of $\omega_i$ onto (the tangent space of)~$A$.
Notice that $\omega'$ is a non-singular $(q-1)$-form on~$A$. 
If we identify $\omega'$ with the corresponding vector field on $A$,
then $\omega'$ is always tangent to ~$L'$.
Notice also that $A \cap L_1 \cap \dotsc \cap L_q$ is a $0$-dimensional
manifold, and therefore 
$\cc(A \cap L_1 \cap \dotsc \cap L_q) = \card{L' \cap L_q}$.

We have to further decompose $A$ in order to transform $L'$ into a \emph{good}
twine. 
Fix a map $p: A \to \K$, such that $p$ is $\Kz$-definable, is $\Cone$, and
has no critical points on~$A$.
For every $x \in A$, let $c(x)$ be the gradient vector of~$p$ at~$x$
(by definition, $c(x)$~is tangent to~$A$).

Define $A_{crit}$ to be the set of points in $A$ such that
$\omega'$ is orthogonal to~$c$, and $A_{reg} := A \setminus A_1$.
After a further cell decomposition, \wloG we can assume that either $A =
A_{reg}$, or $A = A_{crit}$.


If $A = A_{reg}$,  let $\rho$ be the restriction of $p$ to $L'$.
Notice that, by definition of $A_{reg}$, $\rho$ is a definable $\Cone$ 
function without critical points, and hence $L'$ is a good twine in~$A$.
Fix a $\Kt$-definable diffeomorphism $\phi_A$ between $A$ and $\K^d$, 
and define $A_t$ accordingly.
By induction on~$q$, there is $N \in \Nat$ such that
$L' \cap A_t \cap \set{x \in A: \omega' \text{ is not orthogonal to } A_t}$ 
has at most $N$ connected components, where $N$ does not depend on~$t$.
Hence, by definition, $\vb_A(L') \leq N$.
Thus, since $L_q$ is a VRL and $\Omega$ is transverse to~$A$, 
$\card{L' \cap L_q} \leq  N$, and we are done.

If instead $A = A_{crit}$, for every $t \in \K$ let 
$B(t) := \set{x \in A: p(x) = t}$: each $B(t)$ is a $\Kz$-definable set of
dimension $d - 1$.
By induction on~$q$, $L'$ has a uniformly bounded number of connected components $M_1, \dotsc, M_r$.
Moreover, $p$ is constant on each $M_i$, and therefore for each $i \leq r$
there exists $t_i \in \K$ such that  $M_i \subseteq B(t_i)$.
Thus, $M_i \cap L_q \subseteq L_1 \cap \dots \cap L_q \cap B(t_i)$, and
therefore 
\[
A \cap L_1 \cap \dots \cap L_q  \subseteq \bigcup_{i = 1} ^r
B(t_i) \cap L_1 \cap \dots \cap L_q .
\]
By induction on~$d$, there exists a uniform (independent from~$t$)  
bound $r'$ for
$\cc(B(t_i) \cap L_1 \cap \dots \cap L_q)$, and therefore
$\cc(A \cap L_1 \cap \dots \cap L_q ) \leq r r'$.

\subsection{Variants of the Rolle Property}\label{subsec:variants}

In this subsection we compare different notions of Rolle leaves: the original
definition of Rolle leaf (RL), due to Speissegger, which makes sense only for
expansions of the real field was given at the beginning of 
Section \ref{subsec:VRL-main}. 
Alternate Rolle leaves (ARL) and Virtual Rolle leaves (VRL) were defined in Definitions \ref{ARL} and \ref{VRL} respectively.

\begin{dfn}\label{FRL}
   A \textbf{Rolle leaf according to Fratarcangeli} (FRL) is a
   connected multi-leaf $L$ with data $(U, \omega)$, which moreover
   satisfies the condition: for every $m\in\N$, if $X \subset U\times\K^m$ is a definable connected $\Cone$
submanifold of~$U\times\K^m$ of dimension 1, and $X$ intersects $L$
   in at least two points, then $X$ is orthogonal to $\omega$ in at least one
   point (compare with \cite[1.5]{fr}).
   \end{dfn}

\begin{prop}
Let $\K$ be an expansion of the real field. Then every RL is a VRL.
\end{prop}

In particular, we recover Speissegger's theorem is a special case of ours.

\begin{proof}

Let $L \subset \R^n$ be a RL with data $(U, \omega)$.
Let $V \subseteq U$ be a weak cell and $X$ be a good twine in~$V$.
Assume that $\card{X \cap L} > \vb_V(X) =: m$ (the case when $V\subseteq U\times\R^k$ can be treated similarly).
We must show that $X$ is orthogonal to $\omega$ in at least one point; assume,
for contradiction, that this is not the case.
Let $X_i$ be a connected component of~$X$ (notice that $X_i$ is not
necessarily definable).
Since $X$ is a good twine, $X_i$ is not compact; moreover,
$X$ has at most $m$ connected components.
Hence, $X_i$ intersects $L$ in at least two points, for some connected 
component $X_i$.
Thus, since $L$ is a RL and $X_i$ is arc-connected, $X_i$ is orthogonal to
$\omega$ in at least one point, contradiction. 
\end{proof}

\begin{prop}
 Let $\K$ be definably complete. Then every FRL is a VRL.
\end{prop}
In particular, we recover Fratarcangeli's theorem is a special case of ours.

\begin{proof}
Let $L$ be a FRL with data $(U, \omega)$.
Let $V \subseteq U$ be a weak cell and $X \subseteq V$ be a good twine in~$V$,
such that $\card{X \cap L} > \vb_V(X) =: m$ (for simplicity, we are dealing
with the case $n = 0$ in Definition~\ref{VRL}).
By Lemma~\ref{lem:twine-components}, $X$ has at most $m$ connected
components; therefore, there exists $Y$ component of $X$ such that
$\card{Y \cap X} \geq 2$.
Thus, since $L$ is a FRL, $Y$ it orthogonal to $\omega$ at some point.
\end{proof}

There is the following question left.
Let $\K$ be definably complete and Baire.
Let $F: \K^n \to \K$ be a Pfaffian function
(\eg, $F$ is a definable $\Cinf$ function satisfies $d F / d x_i = g_i(x,
F(x))$, for some $\Cinf$ $\Kz$-definable functions $g_i: \K^n \to \K$.
Let $\Kz(F)$ be the expansion of $\Kz$ by $F$.
Is $\Kz(F)$ o-minimal?
Let $C$ be the graph of~$F$, and $\omega$ be the 1-form on $U := \K^{n+1}$
$g_1 d x_1 + \dotsc g_n d x_n - dy$.
Notice that $C$ is a connected Leaf with data $(U, \omega)$.
The question has positive answer if $C$ is either a FRL or a VRL.
We don't know it either is true, but, since being a VRL is a first-order
condition, we can add either the condition ``$C$ is VRL'' to the axioms of
$\Kz(F)$, or we can add the condition ``every graph of a Pfaffian function is
a VRL'' to the axioms of~$\K$.
In both ways, we obtain an axiomatization of $\Kz(F)$ that ensures o-minimality.


\section{Effective bounds}\label{sec:effective}
In this section we apply our results to derive uniform and effective bounds on
some topological invariants (\eg the number of connected components) of sets
definable in the Pfaffian closure of an o-minimal expansion of the real field.

Let $T_0$ be a recursively axiomatized (not necessarily complete) o-minimal theory (if $T_0$ is not recursively
axiomatized, then the effective results of these section are still valid with respect to an oracle
for $T_0$). 
Let $\Rz$ be an o-minimal expansion of the real field, which is a model of
$T_0$ and let $\mathcal P(\Rz)$ be the Pfaffian closure of $\Rz$ (in the sense
of~\cite{speissegger}).

 \begin{dfn}
Let $X \subseteq \R^n$ be definable in~$\mathcal P(\R_0)$.
We call the \textbf{topological complexity of~$X$} (denoted by
$\tc(X)$)
the least $N \in \Nat$, such that there exist:
\begin{enumerate}
\item a simplicial complex $Z$ composed by less than $N$ simplexes,
each of dimension less than~$N$;
\item and a $\mathcal P(\R_0)$-definable homeomorphism to~$f: X \approx \abs{Z}$.
\end{enumerate}
\end{dfn}

Note that, since $\mathcal P(\R_0)$ is o-minimal, the topological complexity is a well defined natural number.

\smallskip

Let $X$ be defined by a formula $\varphi$, where some of the variables are evaluated as a suitable tuple of
parameters. This definition will involve a finite number of Rolle leaves $L_1,\ldots,L_k$. As one can see from the
inductive definition of Pfaffian closure, every leaf $L_i$ will have data $(U_i,\omega_i)$ definable (by a formula
$\phi_i$, where some of the variables are evaluated as a suitable tuple of
parameters) in terms of a finite number of Rolle leaves
$L_{i,1},\ldots,L_{i,n_i}$ of lower complexity (\ie appearing at some earlier
stage of the inductive construction). Hence, to the set $X$ (or better, to its definition $\varphi)$ we can associate a
finite sequence $F_1=L_1,\ldots,F_k=L_k,F_{k+1}=L_{1,1},\ldots,F_{k+1+n_1}=L_{1,n_1},\ldots,F_m$ of Rolle leaves, which
are involved in its definition. The aim of the following
definition is to code the set $X$ by this sequence of leaves (\cf~\cite{gab,fr}).

\begin{dfn}\label{format}
 Let $\mathcal L_P$ be the language of $\R_0$ to which we adjoin a countable set of new predicates
$\set{P_1,\ldots,P_n,\ldots}$. 
A \intro{format} of a definable set $X$ is the following finite sequence of
$\mathcal L_P$-formulae (without parameters): $(\varphi,\mathbf P,\Phi)$, where

\begin{itemize}
 \item for a suitable choice of parameters $\av$, the set $X$ is defined by $\varphi(\cdot,\av)$;
 \item $\mathbf P=(P_1,\ldots,P_m)$ and every $P_i$ represents a Rolle leaf $F_i$ involved in this definition of~$X$;
 \item $\Phi=(\phi_1,\ldots,\phi_m)$ and, for a suitable choice of parameters 
$\av_i$, the formula
$\phi_i(\cdot,\av_i)$ defines the graph of $\omega_i$ on $U_i$, where $(U_i,\omega_i)$ is the data of the leaf $F_i$.
\end{itemize}

\end{dfn}
We did not allow the parameters in the definition. In particular, every other Rolle leaf with the same data $(U,
\omega)$ has the
same format.

\begin{exa}
Let $X = L_1\cup L_2$, where $L_i$ are Rolle leaves with data
$(U_i,\omega_i)$. 
Let $L_{3}$ be a Rolle leaf with $\R_0$-definable data $(U_{3},\omega_{3})$. 
Suppose $(U_1,\omega_1)$ are $\langle\R_0,L_{3}\rangle$-definable and
$(U_2,\omega_2)$ are $\R_0$-definable. 
Let the graphs of $\omega_1,\omega_2,\omega_{3}$ be defined by formulas
$\phi_1(\av_1,\x,\y),\phi_2(\av_2,\x,\y),\phi_3(\av_3,\x,\y)$
respectively, where $\av_1,\av_2,\av_3$ are tuples of parameters. 
Then a format for $X$ is given by the
$L_0\cup\set{P_1,P_2,P_3}$-formulas $(\varphi,\mathbf P,\Phi)$, where $\varphi=P_1\vee P_2;\ \mathbf P=(P_1,P_2,P_3);\
\Phi=(\phi_1,\phi_2,\phi_3)$.

\end{exa}

\smallskip

The next definition requires the notion of Rolle leaf to be first order. This is the reason why we introduced Virtual
Rolle Leaves: the property of being a VRL is type-definable, \ie it can be expressed by a countable (recursive)
conjunction of first order formulae.

\begin{dfn}
Let $X$ be a definable set and $\theta=(\varphi,\mathbf P,\Phi)$ be a format 
for~$X$. 
Let $T_{\theta}$ be the first order theory (in the language of $\R_0$ adjoined
with the predicates $P_1,\ldots,P_m$) with the following \emph{recursive} 
(but not necessarily complete) axiomatization:
\begin{itemize}
\item Axioms of $T_0$;
\item Axioms of Definably Complete Baire Structure;
\item $\phi_i$ defines the graph of a non-singular $\Cone$ 1-form $\omega_i$ on some definable open subset $U_i$;
\item $P_i$ is a VRL with data $(U_i,\omega_i)$.
\end{itemize}

\end{dfn}

We now show the existence of a bound on the topological complexity of~$X$, 
which depends (recursively) only on a format for~$X$.

\begin{thm}
 There is a recursive function $\eta$ which, given a set $X$ definable in $\mathcal P(\R_0)$ and a format $\theta$ for~$X$, returns a natural number $\eta(\theta)$ which is an upper bound on the topological complexity of~$X$.
\end{thm}

\begin{proof}
Let $X$ be a definable set and $\theta=(\varphi,\mathbf P,\Phi)$ be a format
for $X$. Note that, by Theorem~\ref{speissegger's thm}, 
the theory $T_{\theta}$ is o-minimal. In particular, there is a natural
number $N$ such that $\tc(X)<N$. 
Moreover, $T_{\theta}$ is recursively enumerable hence we can recursively
enumerate all the formulas which (for every choice of the parameters) are
provable in this theory. 
Take the first formula in this enumeration which defines a homeomorphism
between the set defined by $\varphi$ and some simplicial complex~$Z$. 
Define $\eta(\theta)$ as the number of complexes which form~$Z$.
\end{proof}

\begin{cor}
There are recursive bounds on the following topological invariants of sets 
definable in $\mathcal P(\R_0)$: number of connected components, sum of the
Betti numbers, number of generators of the fundamental group.
\end{cor}
\begin{proof}
Let $X$ be a set definable in $\mathcal P(\R_0)$ and
$N$ be the recursive bound on $\tc(X)$ given by the above theorem.
Let $F$ be a simplicial complex with at most $N$ simplexes, each of them of
dimension at most~$N$, such that $\abs F$
is homeomorphic to $X$ ($F$ exists by definition of $\tc(X)$).
Clearly, the number of connected components of $\abs F$ (and hence of $X$) 
is at most~$N$.

If $F$ were a \emph{closed} complex, by classical algebraic topology theory,
$N$ would be also a bound for the other mentioned topological invariants.
Otherwise, let $F'$ be the barycentric subdivision of~$F$.
By \cite[Lemma7.1]{EW}, there exists a closed simplicial complex $C$
which is also a sub-complex of~$F'$, such that
$\abs{F'}$ (and hence~$X$) is homotopic to~$\abs C$.
Since $C$ is a closed complex, the number $m$ of simplexes of $C$ gives an
upper bound to the sum of the Betti numbers of $\abs C$ (and hence of~$X$),
and to the number of generators of $\pi_1(X, x_0)$ (for any $x_0 \in X$), and
$m$ is bounded by a recursive function of~$N$.
\end{proof}

We can also obtain bounds on the Hausdorff measure of definable sets.
\begin{lem}\label{lem:bound-measure}
Let $(C_i)_{i \in I}$ be a collection \rom(not necessarily definable\rom)
of sets definable in $\mathcal P(\Rz)$, all  with the same format~$\Phi$, 
such that each is contained in $B(0; 1)$,
and of dimension at most~$d$.
Then there exists a uniform bound on their $d$-dimensional Hausdorff
measure~$\mathcal H^d(C_{i})$.
\end{lem}
\begin{proof}
By a compactness argument, \cite{speissegger} implies that
there exists a uniform bound on  $\gamma(C_{i})$
(where $\gamma$ is as in \cite{wilkie99}).
We then conclude using the Cauchy-Crofton formula: see~\cite{dries03} for
details.
\end{proof}

\section{Conclusion}
\label{sec:conclusion}
We conclude with some open problems.

\begin{open problem}
Let $T$ be an o-minimal theory (expanding \RCF).
Let $\exp$ be a new unary function symbol, and $T(\exp)$ be the following
expansion of~$T$:
\begin{itemize}
\item $T(\exp)$ is definably complete and Baire;
\item $\exp(0) = 1$;
\item $\exp' = \exp$;
\item the graph of $\exp$ is a VRL.
\end{itemize}
Is $T(\exp)$ consistent?
Notice  that, by Theorem~\ref{speissegger's thm}, if consistent,
$T(\exp)$ is o-minimal.
Moreover, any o-minimal structure is either power bounded,
or already defines an exponential function \cite{miller96}
(and therefore in the latter case it is already a model of $T(\exp)$).
Whether $\RCF(\exp)$ is complete or not is not known, but it is surely
consistent, since if $T$ has an Archimedean model, then $T(\exp)$ is
consistent. 
Notice also that there are real closed fields which do not have expansions to
models of $\RCF(\exp)$ \cite{kks}; however, any real closed field has an
elementary extension which admits such an expansion.
\end{open problem}

\begin{open problem}
Is $\RCF(\exp)$ complete?
Assume that Schanuel's Conjecture holds.
Then we can combine our results with those in \cite{mac-wilkie} and obtain the
following result:
if $\K$ is a model of $\RCF(\exp)$, such that every unary function definable
via $\exp\rest(0,1)$ has rational exponent\footnote{If $\K$ is an expansion of
an ordered field, we say that a definable function $g:\K\to\K$ has rational
exponent if there exists $q\in\Q$ such that $\lim_{x\to +\infty}g(x)x^q$ is
finite and nonzero.}, then $\K$ is elementarily equivalent to $\Real(\exp)$.
An analogous result has been obtained in \cite{jones--servi08a} for the
expansion of the real field with a power function $x^{\alpha}$ ($\alpha\in\R$)
and, for $\alpha$ sufficiently generic, it is not necessary to assume
Schanuel's or any other unproven conjecture.
\end{open problem}

\begin{open problem}
Let $\F$ be an o-minimal structure.
For every definable (with parameters) continuous function $g: \F \to \F$, let
$G$ be a new unary function symbol,
and $T'$ be the following expansion of the elementary diagram of $\F$ with the
new symbols:
\begin{itemize}
\item $T'$ is definably complete and Baire;
\item $G$ is a $\Cone$ function, and $G' = g$;
\item the graph of $G$ is a VRL.
\end{itemize}
Is $T'$ consistent?
Again, notice that,
by Theorem~\ref{speissegger's thm}, if consistent, $T'$~is o-minimal.
Moreover, if $\F$ expands the real line, then $T'$ is consistent.
A positive answer to the above question would allow to define an integral for
functions defined in o-minimal structures outside the real line
(at the price of enlarging the structure).
\end{open problem}



\begin{thebibliography}{99}

\bibitem[BO01]{BO01}
A. Berarducci, M. Otero.
\newblock{\em Intersection theory for o-minimal manifolds.}
Annals of Pure and Applied Logic 107 (2001), 87--119.


\bibitem[BS04]{BS04}
A. Berarducci and T. Servi, {\em
An effective version of Wilkie's theorem of the complement and some effective
o-minimality results},
APAL 125(2004) 43-74.

\bibitem[BCR98]{BCR98} J. Bochnak, M. Coste, M-F. Roy, {\em Real Algebraic Geometry}, Springer, 1998.

\bibitem[DMS10]{DMS08} A. Dolich, C. Miller, C. Steinhorn, {\em Structures
  having o-minimal open core}, 
Trans. Amer. Math. Soc.  362  (2010),  no. 3, 1371--1411.


\bibitem[Dries98a]{vddries} L. van den Dries, {\em Dense pairs of o-minimal
  structures}, Fundamenta Mathematicae 157 (1998).

\bibitem[Dries98b]{dries98}
L. van den Dries, {\em
Tame Topology and O-minimal Structures},
Cambridge University Press (1998).

\bibitem[Dries03]{dries03}
L. v.d.Dries, {\em Limit sets on o-minimal structures},
O-minimal Structures, Lisbon 2003,
Proceedings of a Summer School by
the European Research and Training Network Raag,
pp.~172--215.

\bibitem[DM96]{DM96}
L. van den Dries and C. Miller, {\em
Geometric categories and o-minimal structures},
Duke Math. Journal 84:2 (1996) 497--540.


\bibitem[EW08]{EW}
M. J. Edmundo,   A.~Woerheide,
{\em Comparison theorems for o-minimal singular (co)homology},
 Trans. Amer. Math. Soc.  360  (2008),  no. 9, 4889--4912.


\bibitem[Fratarc06]{fr}
S. Fratarcangeli, {Rolle leaves and o-minimal structures},
Ph.D. Thesis, April 2006.

\bibitem[Fratarc08]{frat} S. Fratarcangeli, {\em A first-order version of Pfaffian Closure}, Fund. Math. 198 (2008), 229-254.

\bibitem[GV04]{gab}
A. Gabrielov and N.~Vorobjov,
 {\em Complexity of computations with Pfaffian and Noetherian  functions},
 in {Normal forms, bifurcations and finiteness problems in
              differential equations},
 {NATO Sci.\ Ser.~II Math.\ Phys.\ Chem.},
 {vol.~137},
{pp.~211--250},
{Kluwer Acad.\ Publ.},
{Dordrecht},
{2004}.


\bibitem[GP74]{gull-pollack}
V Guillemin, A. Pollack, {\em Differential Topology}, Prentice-Hall, Englewood Cliffs, 1974.



\bibitem[JS08]{jones--servi08a} G. Jones, T. Servi, {\em On the decidability
  of the real field with a generic power function}, Submitted, 2009.



\bibitem[KM99]{KM99} M. Karpinski, A. Macintyre, {\em A generalization of Wilkie's theorem of the complement, and an application to Pfaffian closure},  Selecta Math. (N.S.)  5  (1999),  no. 4, 507--516.


\bibitem[Kelley55]{Kelley55}
J. Kelley.
General Topology.
Springer, 1955.

\bibitem[Khov91]{khovanskii} A.G. Khovanskii, {\em Fewnomials}, Translations of Mathematical Monographs, 88.
American Mathematical Society, Providence, RI, 1991. viii+139 pp.

\bibitem[KKS97]{kks} F.-V. Kuhlmann, S.~Kuhlmann,  S.~Shelah.
{\em Exponentiation in power series fields.}
Proc. Amer. Math. Soc.  125  (1997),  no. 11, 3177--3183.

\bibitem[MW96]{mac-wilkie} A. Macintyre, A. Wilkie,
{\em On the decidability of the real exponential field}, in: Kreiseliana, A. K.
Peters, Wellesley, MA, 1996, pp. 441--467.

\bibitem[Marker97]{marker} D. Marker, {\em Khovanskii's theorem.  Algebraic model theory (Toronto, ON, 1996)},  181--193,
NATO Adv. Sci. Inst. Ser. C Math. Phys. Sci., 496, Kluwer Acad. Publ., Dordrecht, 1997.

\bibitem[Maxwell98]{max} S. Maxwell, {\em A general model completeness result for
expansions of the real ordered field}, Ann. Pure Appl. Logic 95 (1998) 185--
227.

\bibitem[Miller96]{miller96} C. Miller.
{\em A growth dichotomy for o-minimal expansions of ordered fields.}
 Logic: from foundations to applications (Staffordshire, 1993),
 385--399, Oxford Sci. Publ., Oxford Univ. Press, New York,  1996.

\bibitem[Miller01]{miller} C. Miller.
\newblock{\em Expansions of dense linear orders with the intermediate value property.}
\newblock J. Symbolic Logic 66 (2001), no. 4, 1783--1790.


\bibitem[Oxtoby80]{Oxtoby80}
J. C. Oxtoby,
{\em Measure and Category.
A Survey of the Analogies between Topological and Measure Spaces},
\newblock Springer, 1980.

\bibitem[Servi07]{servi-articolo}
T. Servi.
\newblock{\em Noetherian Varieties in Definably Complete Structures,} Logic and Analysis 1 (2008), 187--204.



\bibitem[Speiss99]{speissegger} P. Speissegger, {\em The Pfaffian closure of an o-minimal structure}  J. Reine Angew. Math.  508  (1999), 189--211.

\bibitem[Wilkie99]{wilkie99} A. Wilkie, {\em A theorem of the complement and some new
o-minimal structures}, Selecta Math. (N.S.) 5 (1999).
397--421.


\end{thebibliography}
\end{document}